\documentclass[reqno,11pt]{amsart}

\usepackage{amsmath}
\usepackage{amsfonts, amssymb,amsthm,amscd,stmaryrd}
\usepackage{amsthm}
\usepackage{mathrsfs}
\usepackage{amsfonts}
\usepackage{parskip} 
\usepackage[hidelinks]{hyperref}
\usepackage[margin=3.57cm]{geometry}
\usepackage{enumerate,color,graphicx,
}
\usepackage{psfrag}  
\usepackage{url}
\usepackage{mdwlist}   
\usepackage{mathrsfs}
\usepackage{etoolbox}

\makeatletter
\patchcmd{\@part}{\Huge}{\normalsize}{}{} 
\makeatother

\numberwithin{equation}{section}

\allowdisplaybreaks

\newtheorem{theorem}{Theorem}[section]

\newtheorem{proposition}[theorem]{Proposition}
\newtheorem{lemma}[theorem]{Lemma}
\newtheorem{conjecture}{Conjecture}

\newtheorem{definition}[theorem]{Definition}
\newtheorem{corollary}[theorem]{Corollary}

\newcommand{\ep}{\varepsilon}

\newcommand{\R}{\mathbb{R}}
\newcommand{\N}{\mathbb{N}}
\newcommand{\Z}{\mathbb{Z}}

\newcommand{\eps}{\ep}

\newcommand{\latt}{\Z^2}
\newcommand{\lattv}{\Z^2}

\newcommand{\SHSW}{\aleph \! \! \! \smallsetminus \! \! \! \aleph  \! \daleth\mathsf{HSW}}
\renewcommand{\SHSW}{\raisebox{-0.04 em}{\includegraphics{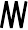}}\mathsf{HSW}}

\newcommand{\northeast}{\mathsf{NE}}

\newcommand{\southeast}{\mathsf{SE}}

\newcommand{\bridgeprod}{\mathsf{B}\Pi}

\newcommand\La{\Lambda}

\newcommand{\bbU}{\mathbb{U}}
\newcommand{\bbUi}{{\mathbb U}^\infty}

\newcommand{\ueight}{{\mathbb U}^8}
\newcommand{\originwalk}{o}

\newcommand\tall{t}

\theoremstyle{remark}

\newcounter{mycount}

\newcommand{\remove}[1]{}
\newcommand{\add}[1]{{\color{black}#1}}

\def\mik{1}
\newcommand\cpsfrag[2]{\ifnum\mik=1\psfrag{#1}{#2}\fi}

\newcommand\ga{\gamma}     

\newcommand{\hidden}[1]{}

\newcommand{\joinloc}{\mathsf{Join}}

\newcommand{\SAW}{\mathsf{SAW}}
\newcommand{\SAB}{\mathsf{SAB}}
\newcommand{\HSW}{\mathsf{HSW}}

\newcommand{\CHSW}{\mathsf{HorConfHSW}}
\newcommand{\FBHSW}{\mathsf{FewBranchHSW}}
\newcommand{\SAP}{\mathsf{SAP}}
\newcommand{\WSAP}{\mathsf{WSAP}}
\newcommand{\FSAP}{\mathsf{RSAP}}
\newcommand{\MJ}{\mathsf{MJ}}
\newcommand{\Strip}{\mathsf{Strip}}
\newcommand{\Rect}{\mathsf{Rect}}
\newcommand{\ajll}{\mathsf{list}}
\newcommand{\AJLL}{\mathsf{JoinLocLists}}

\newcommand{\bbN}{\mathbb{N}}
\newcommand{\bbZ}{\mathbb{Z}}
\newcommand{\bbR}{\mathbb{R}}

\newcommand{\bbH}{\mathbb{H}}
\newcommand{\calR}{\mathcal{R}} 
\newcommand{\calS}{\mathcal{S}} 
\newcommand{\calT}{\mathcal{T}} 

\newcommand{\emm}{m_0}
\newcommand{\EMM}{M_0}
\newcommand{\concat}{\mathcal{C}}
\newcommand{\eightLa}{\Lambda^8}
\newcommand{\eightLambda}{\Lambda^8}
\newcommand{\eightOm}{\bbU^8}
\newcommand{\eightOmega}{\bbU^8}
\newcommand{\inftyOm}{\bbU^\infty}

\newcommand{\inftyLa}{\Lambda^\infty}
\newcommand{\inftyLambda}{\Lambda^\infty}

\title[Bounding the number of self-avoiding walks]{Bounding the number of self-avoiding walks: Hammersley-Welsh with polygon insertion}
\date{\today}
\author[H.~Duminil-Copin]{Hugo Duminil-Copin}
\address{Institut des Hautes \'Etudes Scientifiques and University of Geneva,
	Le Bois-Marie 35 route de Chartres,
	91440 Bures-sur-Yvette, France.
	Member of SwissMAP.}
\email{duminil@ihes.ch, hugo.duminil@unige.ch}

\author[S.~Ganguly]{Shirshendu Ganguly}
\address{Department of Statistics, 
  U.C. Berkeley,
  Berkeley, CA, 94720-3840, U.S.A.}
\email{sganguly@berkeley.edu}

\author[A.~Hammond]{Alan Hammond}
\address{Departments of Mathematics and Statistics, 
  U.C. Berkeley,
  Berkeley, CA, 94720-3840, U.S.A.}
\email{alanmh@stat.berkeley.edu}
\keywords{}

\author[I.~Manolescu]{Ioan Manolescu}
\address{D\'epartement de math\'ematiques,
Universit\'e de Fribourg,
Chemin du Mus\'ee 23,
CH-1700 Fribourg, Switzerland.
Member of SwissMAP.}
\email{ioan.manolescu@unifr.ch}

\keywords{}

\thanks{2010 Mathematics Subject Classification. Primary:  60K35.  Secondary: 60D05}

\begin{document}

\maketitle

\begin{abstract}
        Let $c_n = c_n(d)$ denote the number of self-avoiding walks of length~$n$ starting at the origin in the Euclidean nearest-neighbour lattice~$\Z^d$.
        Let $\mu = \lim_n c_n^{1/n}$ denote the connective constant of~$\Z^d$.
        In 1962, Hammersley and Welsh~\cite{HamWel62} proved that, for each $d \geq 2$, there exists a constant $C > 0$
        such that $c_n \leq \exp(C n^{1/2}) \mu^n$ for all $n \in \N$. 
        While it is anticipated that $c_n \mu^{-n}$ has a power-law growth in $n$,
        the best known upper bound in dimension two has remained of the form $n^{1/2}$ inside the exponential.
      
		The natural first improvement to demand for a given planar lattice is a bound of the form 
		$c_n \leq \exp (C n^{1/2 - \eps})\mu^n$, where $\mu$ denotes the connective constant of the lattice in question. 
		We derive a bound of this form for two such lattices, for an explicit choice of $\eps > 0$  in each case. 
		For the hexagonal lattice $\bbH$, the bound is proved for all $n \in \N$; 
		while for the Euclidean lattice $\Z^2$, it is proved for a set of $n \in \N$ of limit supremum density equal to one.

        A power-law upper bound on $c_n \mu^{-n}$ for $\bbH$
        is also proved, contingent on a non-quantitative assertion concerning this lattice's connective constant.
\end{abstract}
\vspace{6mm}

\tableofcontents

\section{Introduction}

\subsection{Definitions and results}

We will denote by $\N$ the set of {\em positive} integers.
For $u \in \R^2$, let $\| u \|$ denote the Euclidean norm of $u$. The two-dimensional nearest-neighbour Euclidean lattice $\latt = (\lattv,E(\lattv))$
has origin $O = (0,0)$.  

A {\em walk} of length $n  \in \N \cup \{ 0 \}$ is a map $\gamma:\{0,\dots,n \} \to \lattv$ 
such that $\| \gamma_{i + 1} - \gamma_i \| = 1$ for each $i \in \{0,\dots,n-1\}$.
An injective walk is called {\em self-avoiding}. 
Write $\SAW_n$ for the set of all self-avoiding walks of length $n$ with $\gamma_0 = O$.  
The book~\cite{MadSla13} and lecture notes~\cite{BauDumGoo12} offer introductions to the topic of self-avoiding walk.

It follows from a simple sub-multiplicativity result that there exists a constant $\mu$, called the {\em connective constant}, such that 
\begin{align*}
	\mu^n \leq | \SAW_n |\leq \exp (o(n))\, \mu^n. 
\end{align*}
The question that we address here concerns the upper bound. 
Until recently, the best known upper bound on $|\SAW_n|$ was given by the celebrated work of Hammersley and Welsh~\cite{HamWel62}. 

\begin{theorem}[Hammersley-Welsh bound]\label{thm:HW}
For $d \geq 2$,	there exists $C > 0$ such that, for all $n \in \N$,
	\begin{align}\label{eq:HW}
		 |\SAW_n| \leq \exp (C n^{1/2})\,\mu^n. 
	\end{align}
\end{theorem}

Our definition of $\SAW_n$ is made for $d=2$ because this article's results concern this dimension, but the model's definition in higher dimensions is no different.
Theorem~\ref{thm:HW}'s proof depends on  unfolding self-avoiding walks so that certain special self-avoiding walks known as bridges result. 
 Kesten substantially improved Theorem~\ref{thm:HW} for dimensions $d\ge3$ in~\cite{Kes64}.  Recently,~\cite{Hut18}
has observed that, in any dimension at least two, the sub-ballisticity of self-avoiding walk~\cite{DumHam13} entails that the constant $C > 0$ in Theorem~\ref{thm:HW} may be chosen arbitrarily. 
In sufficiently high dimension~$d$, Slade~\cite{Slade89}  has proved the far stronger inference that $\vert \SAW_n \vert \sim c_d \, \mu^n$ for some $c_d > 0$; when $d \geq 5$, this result is due to Hara and Slade~\cite{HarSla92b}.

Our aim is to improve the exponent $1/2$ in the power of $n$ in the exponential in Theorem~\ref{thm:HW} in the two-dimensional case. We mention first that it is expected that 
\begin{align*}
	 |\SAW_n| = n^{11/{32} + o(1)}\mu^n \, , 
\end{align*}
where the exponent $11/32$ is predicted~\cite{Nie82} for any planar lattice, in contrast to the lattice-dependent value of $\mu$.  See~\cite{LawSchWer04} for discussions of this prediction and the conjectural ${\rm SLE}_{8/3}$ scaling limit of planar self-avoiding walk.

The improved bound that we present requires a power-law lower bound on the correction to exponential growth  for self-avoiding polygons 
 of a given length, where such polygons  are in essence self-avoiding walks that return to their starting points.  
This bound is available only subsequentially for $\Z^2$. Here then is the first of our main results.

\begin{theorem}[Improved Hammersley-Welsh on $\Z^2$]\label{thm:HW1}
For any $\eps <\frac1{466}$, there exist infinitely many values of $n \in \N$ such that 
	\begin{align}\label{eq:HWimp1}
		|\SAW_n| \leq \exp (n^{1/2 - \eps})\, \mu^n \, . 
	\end{align}
	Indeed, there are infinitely many $j \in \N$
	for which all $n \in \N \cap [ j,j^{4/3} ]$ satisfy~(\ref{eq:HWimp1}).
\end{theorem}

The technique of proof of Theorem~\ref{thm:HW1} is rather robust, and it is plausible that a similar result may be obtained  for many planar lattices. 
For a particular lattice, special structure permits the derivation of a stronger result.
 This lattice is the hexagonal lattice $\bbH$, which is dual to the triangular lattice $\mathbb T = \bbZ+{\rm e}^{{\rm i}\pi/3}\bbZ$: see the later  Figure~\ref{fig:triangle} for a depiction. 
On $\bbH$, an analysis of a discretely holomorphic observable, which was exploited in~\cite{DumSmi12} and is discussed in Section~\ref{sec:observable},  leads to a bound of the form~(\ref{eq:HWimp1}) for all~$n$. Let $\SAW_n(\bbH)$ denote the set of self-avoiding walks of length $n$ starting from a given vertex in $\bbH$.
The hexagonal lattice's connective constant $\mu(\bbH)$ was proved to equal $\sqrt{2 + \sqrt 2}$ in~\cite{DumSmi12}. 

\begin{theorem}[Improved Hammersley-Welsh on $\bbH$]\label{thm:HW2}
	Let $\eps \in ( 0, \frac1{42})$. Then, for any $n \in \N$ high enough, 
	\begin{align}\label{eq:HWimp2}
		|\SAW_n(\bbH)| \leq \exp (n^{1/2 - \eps} ) \, \mu(\bbH)^n \, . 
	\end{align}
\end{theorem}

Beyond these two bounds, we present a third theorem, in which 
 a polynomial upper bound on the hexagonal lattice's normalized walk count $\mu(\bbH)^{-n} |\SAW_n(\bbH)|$ is obtained subject to a qualitative conjecture concerning numbers of self-avoiding walks. In order to pose the conjecture, 
consider the universal cover $\bbUi$ of $\mathbb H$ with singularity at the origin; this informal description will be made precise in Section~\ref{s.cover}.
Let $\pi_\infty$ be the canonical projection of $\bbUi$ to $\mathbb H$. 
For $k\geq 0$, let $\La_k$ be the ball of radius $k$ for the graph distance on $\mathbb T$, centred at the origin; we will view $\La_k$ as a set of faces of $\bbH$; note that $\La_0$ contains one face, whose centre is the origin. 
Let $\bbUi_k:=\bbUi\setminus \pi_\infty^{-1}(\Lambda_k)$,  this being  the universal cover of $\mathbb H \setminus \La_k$; note that $\bbUi \subset \bbUi_0$.

\begin{figure}
	\begin{center}
	   	\includegraphics[scale = 0.5]{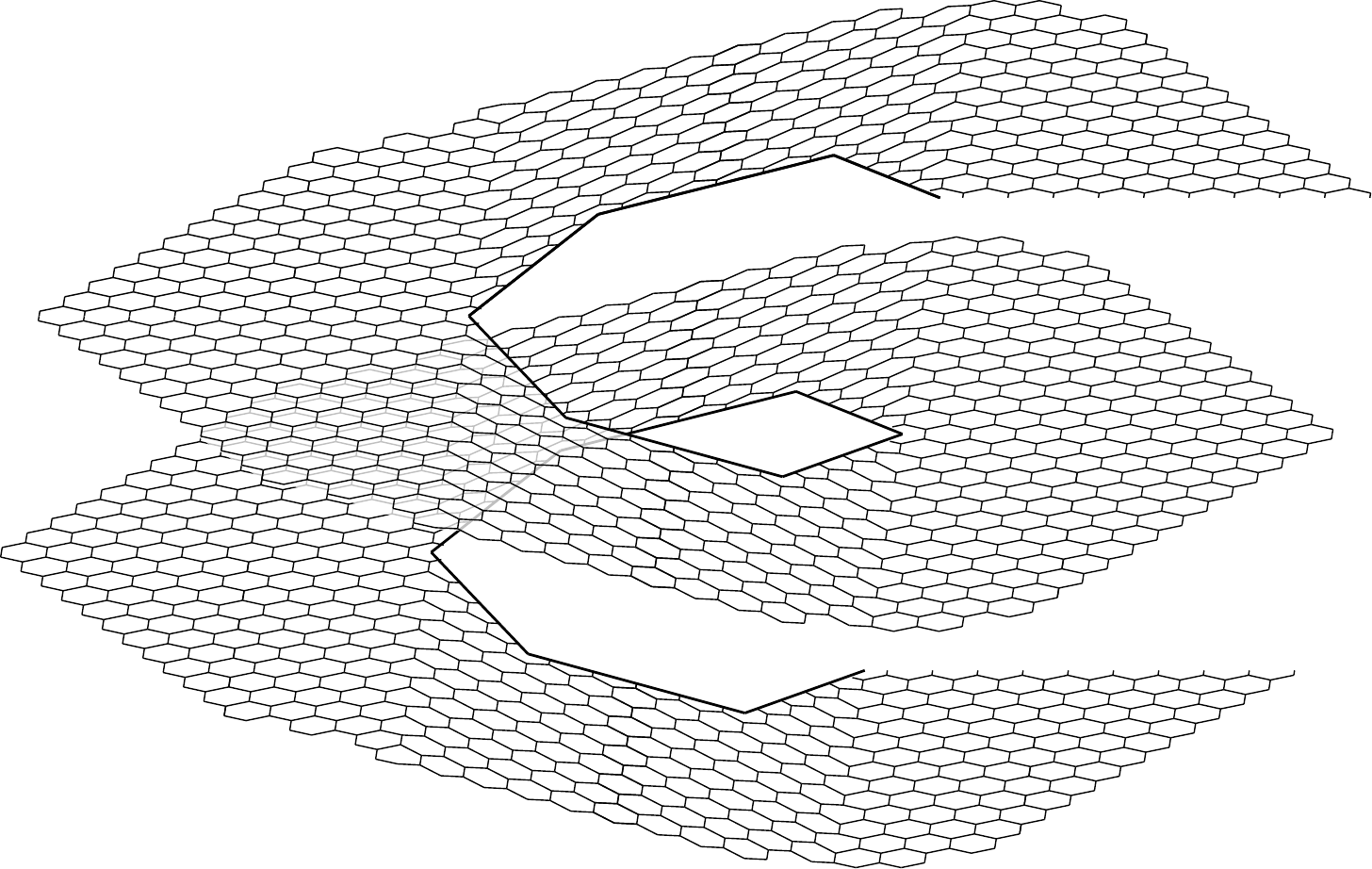}
		\caption{A finite section of some $\bbUi_k$ is  embedded in $\bbR^3$.}
	\label{fig:bbU}
	\end{center}
\end{figure}

Let $k \in \N$. For $v \in \bbUi_k$, set $\SAW_n(\bbUi_k,v)$ equal to the set of self-avoiding walks in $\bbUi_k$ of length $n$ that start at $v$; absence of translation invariance entails the specification of~$v$ in this notation since it is the cardinality of these sets that will concern us. When a self-avoiding walk is split into two, the pieces are also self-avoiding. Thus, the   sequence $\big\{ \sup_{v \in \bbUi_k} \big\vert \SAW_n(\bbUi_k,v) \big\vert : n \in \N \big\}$ is submultiplicative, so that Fekete's lemma permits us to define 
   \begin{equation}\label{e.uik}
    \mu(\bbUi_k) \, :=  \, \lim_{n\rightarrow\infty} \, \sup_{v \in \bbUi_k} \big\vert \SAW_n(\bbUi_k,v) \big\vert^{1/n} \, .
    \end{equation}

\begin{conjecture}\label{conj:1}
	There exists $k \in \N \cup \{0 \}$ for which $\mu(\bbUi_k)$ is equal to $\mu(\bbH)$.
\end{conjecture}

The conjecture is in our view likely to be valid even when $k = 0$.
In principle, however, a typical self-avoiding walk on $\bbUi_k$ may wind continually around the axis of the universal cover, resulting in  $\mu(\bbUi_k) > \mu(\bbH)$.


Here is our third main result.

\begin{theorem}\label{thm:a}
Assume Conjecture~\ref{conj:1}.
    There exist positive constants $C_0$ and $C$ such that, for $n \in \N$,
    \begin{equation}\label{eq:a}
    	\big\vert  \SAW_n(\bbH) \big\vert \, \leq \, C_0 \, n^C \mu(\bbH)^n \, .
	\end{equation}
\end{theorem}


Although this theorem is a conditional result, it is plausible that it considerably simplifies the task of proving~\eqref{eq:a}. 
Indeed, the task of deriving a power-law upper bound for $\mu(\bbH)^{-n} c_n$ has been reduced to  proving that the number of self-avoiding walks on $\bbU_k$ is bounded above by $\mu(\bbH)^{n+o(n)}$ with a non-quantitative  $o(n)$ bound.


The introduction continues by presenting some fundamental notation and concepts. It ends by explaining the paper's structure. 

\subsection{Notation}\label{s.notation}

For $a < b$, write $\llbracket a,b \rrbracket := [a,b]\cap\bbZ$.

Let $(e_1,e_2)$ denote the standard unit vectors that generate $\Z^2$.
For a point $z \in \bbZ^2$, we write $x(z)$ and $y(z)$ for its first and second coordinates. 
For a finite set of points $A \subset \bbZ^2$, write  
$y_{\max}(A) = \max\{y(a) :\, a \in A\}$ and $y_{\min}(A) = \min\{y(a) :\, a \in A\}$; the same applies to the $x$-coordinate.
For walks $\gamma$,  $y_{\max}(\gamma)$ and $y_{\min}(\gamma)$ refer to the definitions above, with $\gamma$ being the set of points visited by $\gamma$. 
\newcommand{\height}{\mathsf{h}}
Then $\height(\gamma) = y_{\max}(\gamma) - y_{\min}(\gamma)$ is called the {\em height} of $\gamma$. 
(The word `height' will also refer to the $y$-coordinate of a point in $\Z^2$.)

For a walk $\gamma \in \SAW_n$ and $0 \leq a \leq b \leq n$, $\gamma_{[a,b]}$ is the walk $(\gamma_a,\gamma_{a+1},\dots, \gamma_b)$. We define the {\em concatenation} $\gamma\circ\gamma'$ of two walks $\gamma$ and $\gamma'$ starting from $O$ of respective lengths $n$ and~$m$ by the formula
$$(\gamma\circ\gamma')_k=\begin{cases} \, \gamma_k&\text{ for }k\le n\\
\, \gamma_n+\gamma'_{k-n}&\text{ for }k\ge n.\end{cases}$$
The length of a self-avoiding walk $\gamma$ will be denoted by $\vert \gamma \vert$.
\subsection{Multi-valued maps}\label{s.mvp}
Our main arguments will be presented using a simple combinatorial inequality.

Let $A$ and $B$ be finite sets, and 
let $\mathcal{P}(B)$ denote the power set of~$B$.
A {\em multi-valued map} from $A$ to $B$ is a function $\Phi: A \to \mathcal{P}(B)$. An arrow is a pair $(a,b) \in A \times B$ for which $b \in \Phi(a)$; such an arrow is said to be outgoing from $a$ and incoming to $b$. We consider multi-valued maps in order to find upper bounds on $\vert A \vert$, and for this, we need the next lemma, which gives upper (and lower) bounds on the number of incoming (and outgoing) arrows. 
\begin{lemma}\label{l.mvm}
Let $\Phi: A \to \mathcal{P}(B)$.
Set $\emm$ to be the minimum over $a \in A$ of the number of arrows outgoing from $a$, and $\EMM$ to be the maximum over $b \in B$ of the number of arrows incoming to~$b$. Then
  $\vert B \vert \geq \emm \EMM^{-1} \vert 
A \vert$.
\end{lemma}
\noindent{\bf Proof.} The quantities $\EMM \vert B \vert$ and $\emm \vert A \vert$
are upper and lower bounds on the total number of arrows. \qed

When there is merely one outgoing arrow from each element of $A$, we call the multi-valued map $\Phi$
{\em degenerate}; in this case, we may instead view $\Phi$ as a function from $A$ to $B$.

\subsection{Self-avoiding bridges and polygons}\label{s.bridge}\label{s.def}
A self-avoiding walk in $\Z^2$ of length~$n$  that starts at $O$ and satisfies $0 < y(\gamma_k) \leq y(\gamma_n)$ for $1 \leq k\leq n$ is called a bridge.
Let $\SAB_n$ denote the set of self-avoiding bridges of length $n$.

When bridges are concatenated, the outcome is also a bridge. The resulting supermultiplicativity leads to
$$
	|\SAB_n| \leq \mu_b^n \, , 
$$	
where $\mu_b = \lim_n |\SAB_n|^{1/n} \leq \mu$.

Self-avoiding polygons play an essential role in the proofs of the strengthened Hammersley-Welsh bounds, Theorems~\ref{thm:HW1} and~\ref{thm:HW2}. 

A self-avoiding walk $\gamma: \{ 0, \dots, n \} \to \Z^2$ is called {\em closing} if $\| \gamma_n - \gamma_0 \| = 1$. A closing walk has odd length.
A self-avoiding polygon is formed from such a closing walk $\gamma$  by taking the union of the set of nearest-neighbour edges that interpolate consecutive endpoints of~$\gamma$ with the further edge 
$\{\gamma_n,\gamma_0 \}$. This polygon has length $n+1$, which is even. A self-avoiding polygon of length four is called a plaquette; it consists of the four edges that border a certain unit square.

Let $\SAP_n$ denote the set of equivalence classes of self-avoiding polygons of length $n \in 2\N$, where two such polygons are identified if there exists a vector in $\Z^2$
that translates one to the other. 

An element of $\SAP_n$ corresponds to $2n$ closing walks of length $n-1$ that start at the origin: any one of $n$ edges may be removed from a length~$n$ polygon, and a choice of two orientations then made for the resulting closing walk.

The next lemma follows from an argument of Kesten~\cite{Kes63}.
\begin{lemma}\label{l.sapsab}
For any $n \in \N$,
$$
 \big\vert \SAP_{2n + 2} \big\vert \geq \frac{\big\vert \SAB_n \big\vert^2}{4(2n+1)n(n+1)^3}\, .
$$
\end{lemma}
\noindent{\bf Proof.}
The final inequality of the proof of~\cite[Theorem~$2.9$]{BauDumGoo12} asserts, when $d=2$, that the number of closing walks in~$\Z^2$ of length~$2n+1$ that begin at the origin is at least  $n^{-1} (n+1)^{-2} (2n+1)^{-1} \vert \SAB_n \vert^2$. The number of elements of $\SAP_{2n+2}$ is thus seen to be at least the product of $2^{-1}(2n+2)^{-1}$ and the latter quantity. \qed

\noindent{\em Two remarks: 1.} We study self-avoiding walks, but walks in which the constraint of self-avoidance is violated play a vital role in the proofs of Theorems~\ref{thm:HW1} and~\ref{thm:HW2}. Walks, bridges and polygons are understood to be self-avoiding: the  violation of this constraint  will be noted with emphasis. 

{\em 2.}
We will make such assertions as `a walk and a polygon are disjoint' and `a polygon intersects a given set of vertices'. These assertions are abuses of notation, but their meaning is, we hope, clear.

\subsection{Structure of the paper}
The celebrated Hammersley-Welsh unfolding argument proves Theorem~\ref{thm:HW}. 
Our first two principal results depend on developing this argument, and we begin Section~\ref{sec:HW} by reviewing it. We then explain roughly how we will modify the argument with the use of polygon insertion.

This heuristic discussion will make clear that certain basic tools are needed: polygon abundance,  and Madras' join technique for polygons in~$\Z^2$. These tools are reviewed or developed in  
Section~\ref{sec:SAP}. By polygon abundance, we mean a power-law lower bound on the normalized number of certain useful polygons of given length~$n$. This bound is obtained subsequentially in~$n$ for~$\Z^2$. 
A stronger bound valid for the hexagonal lattice and for {\em all} lengths has a proof that involves the parafermionic observable. Section~\ref{sec:observable} contains this derivation. 

Sections~\ref{sec:proofsone} and~\ref{sec:proofstwo} offer the rigorous counterpart to the heuristics of Section~\ref{sec:HW}:
 assertions of polygon abundance will be used  to prove the square lattice Theorem~\ref{thm:HW1}. Section~\ref{sec:proofsone} begins by explaining how the task will be carried out over the two sections. Section~\ref{sec:proofstwo} ends with some brief comments concerning the prospect for obtaining  improvements to Theorem~\ref{thm:HW1} by making changes to this theorem's proof. Certain perturbations of the argument for Theorem~\ref{thm:HW1} are needed alongside the stronger polygon abundance estimate available for the hexagonal lattice to yield Theorem~\ref{thm:HW2}. The construction and changes are described in Appendix~\ref{sec:madras_join_hex}.  
 
The main part of the paper ends with two sections treating the hexagonal lattice, the first of which we have mentioned.
The final Section~\ref{s.thm:a}  contains the proof of the strong but conditional Theorem~\ref{thm:a}.

\medskip

\noindent{\bf Acknowledgments.} 
The first author is supported by the ERC grant CriBLaM, an IDEX grant from Paris-Saclay and the NCCR SwissMAP, 
the second by a Miller Research Fellowship, the third by NSF grant DMS-$1512908$ and the fourth by the the NCCR SwissMAP.

The authors thank Alexander Glazman and Matan Harel for several valuable conversations. 
The first author thanks Tony Guttmann, Nick Beaton and Iwan Jensen for useful discussions regarding Conjecture~\ref{conj:1}.

Finally, we thank a referee who read the manuscript meticulously and offered many useful comments.

\section{The Hammersley-Welsh argument and its prospective modification}\label{sec:HW}

In two subsections, we recount this argument and, in a third, we lay out a plan to modify it.

\subsection{Bridge lists}\label{sec:bridge}
We define the bridge-product space 
$\bridgeprod_{\ell,j}$, for $\ell,j \in \N$,
to be the set of ordered lists $\big(b_1,\dots,b_j\big)$ of length~$j$ whose elements are bridges whose respective lengths~$\ell_k$ 
satisfy $\sum_{k=1}^{j} \ell_k= \ell$ and whose respective heights $h_k$ form a strictly decreasing sequence.
\remove{The latter, height, condition will not be invoked during the proofs of our main results; it is needed for a bound specific to the classical derivation of Theorem~\ref{thm:HW}. }

\begin{lemma}\label{l.bridgeproduct} We have the following two properties.
    \begin{enumerate}
    \item
   If $j \leq \big( 2\ell \big)^{1/2}$, then
    $ \big\vert \bridgeprod_{\ell,j} \big\vert \leq \ell^j \mu^\ell$, and if
     $j>  2 \ell^{1/2}$, then $ \bridgeprod_{\ell,j} = \emptyset$.
    \item
    There exists $C > 0$ such that,
    for $\ell \in \N$,
    $$
     \bigg\vert \bigcup_{j=1}^\ell \bridgeprod_{\ell,j} \bigg\vert \, \leq \,  \exp ( C \ell^{1/2}) \mu^\ell \, .
    $$
    \end{enumerate}
\end{lemma}

\noindent{\em Remark.} The proofs demonstrate the statements with the change $\mu \to \mu_b$. However, the results under review leading to Theorem~\ref{thm:HW}
demonstate that $\mu_b = \mu$. We will thus make no further reference to~$\mu_b$.

\noindent{\bf Proof of Lemma~\ref{l.bridgeproduct}.} 
{\bf (1).} \add{Fix $j$ and $\ell$. If $j > ( 2\ell )^{1/2}$, the sum of heights of any element $\big( b_1,\dots, b_j \big) \in \bridgeprod_{\ell,j}$ is greater than $1+ 2 +\dots + j > \ell$, which proves that no such bridge-list exists. Assume now that $j \leq ( 2\ell )^{1/2}$. }

The concatenation map $\concat$ that sends $\big(b_1,\dots,b_j\big)$ to $b_1 \circ \dots \circ b_j$  
maps
$\bridgeprod_{\ell,j}$ to $\SAB_\ell$. Treating this operation as a multi-valued map, though it is a degenerate one, we apply Lemma~\ref{l.mvm} with $A = \bridgeprod_{\ell,j}$ and $B = \SAB_\ell$. Clearly, $\emm = 1$, while the fact noted early in Section~\ref{s.bridge} implies that $\vert B \vert \leq \mu^\ell$. 
The map $\concat$ is not injective, because of the ambiguity in how several bridges are concatenated to form a given output bridge. We gauge distance from injectivity by asserting that the quantity $\EMM$ in Lemma~\ref{l.mvm} 
 is at most $\ell^j$. The assertion is valid  because, for given $b \in \SAB_\ell$, the identity of an element in $\concat^{-1}(b)$ may be inferred from the $j$ concatenation points of consecutive bridges along $b$, and each of these points lies in a set of cardinality $\ell$. 

{\bf (2).}  Fix $\ell$. When the concatenation map is extended to act on $\bigcup_{j=1}^\ell \bridgeprod_{\ell,j}$, its range again lies in  $\SAB_\ell$. Applying Lemma~\ref{l.mvm} once more, it is enough to argue that $\EMM \leq  \exp (C \ell^{1/2} )$
for a suitable positive constant~$C$. 
To bound $\EMM$, note that for given $b \in \SAB_\ell$, it suffices to determine the strictly decreasing sequence $\big\{ h_i: i \geq 1 \big\}$ \add{of heights of $b_1,\dots, b_j$ in order to determine $b_1,\dots, b_j$.}
This sequence is a decreasing partition of the height of $b$, which is at most $\ell$.

In view of a much finer result of  Hardy and Ramanujan~\cite{HarRam17asymp}, there exists a positive constant $C$ such that the number of decreasing partitions of a given integer $k \geq 1$
is bounded above by $\exp ( C k^{1/2} )$. 
Taking $k$ to equal the height of $b$, we verify that Lemma~\ref{l.mvm}'s $\EMM$ is at most $\exp ( C \ell^{1/2} )$. \qed

\subsection{The Hammersley-Welsh argument}
Although the proof of Theorem~\ref{thm:HW} that we are reviewing  is valid in any dimension $d \geq 2$, our development of this argument requires that $d=2$. We thus take $d=2$ throughout the article; even if in this section it is purely as a matter of notational convenience.

A {\em half-space self-avoiding walk} (which, in accordance with our convention, we will call simply a half-space walk) of length $n$ is a self-avoiding walk $\gamma$ of length $n$ starting at $O$ with $0 < y(\gamma_k)$ for all $0 < k \leq n$. 
Write $\HSW_n$ for the set of such walks. 


Decomposing any self-avoiding walk of length~$n$ at its final lowest point, we find that
\begin{align}\label{eq:SAW_vs_HSW}
	|\SAW_n| \leq \sum_{k = 0}^n |\HSW_{k+1}| \, |\HSW_{n-k}| \, .
\end{align}
The second term of the product accounts for the segment of the walk after its final lowest point; 
it is the translation of a half-space self-avoiding walk. 
The segment of the walk before its final lowest point, with a downward vertical edge added at the end, 
is the translation and reversal of a half-space self-avoiding walk, and contributes to the first term of the product. 
The $+1$ in the index of the first term is due to the added downward step.  

\begin{lemma}\label{l.hsw}
There exists $C > 0$ such that, for $n \in \N$,
$$
  \vert \HSW_n \vert \leq \exp ( C n^{1/2} ) \,  \mu^n \, .
$$
\end{lemma}
\noindent{\bf Proof of Theorem~\ref{thm:HW}.} This is immediate from the two preceding bounds. \qed

\noindent{\bf Proof of Lemma~\ref{l.hsw}.}
Let $n \in \N$. In light of Lemma~\ref{l.bridgeproduct}(2), it is enough to construct an injective map $\Psi: \HSW_n \to 
 \bigcup_{j=1}^n \bridgeprod_{n,j}$ from length-$n$ half space walks to bridge-lists of sum length $n$.

Let $\gamma\in \HSW_n$. 
	Define the {\em record points} of $\gamma$ recursively as follows. 
	Let $a_0=0$. For $k \geq 1$ odd, let $a_k$ be the last time that $\gamma_{[a_{k-1},n]}$ reaches its highest $y$-coordinate. 
		For $k \geq 2$ even, let  $a_k$ be the last time that $\gamma_{[a_{k-1},n]}$ reaches its lowest $y$-coordinate. 
			Formally, 
	\begin{align}
	a_k := \begin{cases}  \,
	\max\{ t > a_k : \, y(\gamma_t) = y_{\max}(\gamma_{[a_{k-1},n]})\} \quad &\text{ if $k$ is odd},\\
\, 	\max\{ t > a_k : \, y(\gamma_t) = y_{\min}(\gamma_{[a_{k-1},n]})\} \quad &\text{ if $k$ is even}.
	\end{cases}
	\end{align}
	Stop this procedure on the first occasion that $a_k = n$.

	We call the walks $\gamma_{[a_k,a_{k+1}]}$, $0 \leq k < r$, the {\em branches} of $\gamma$. 
In this way, we associate to~$\gamma$ 
its {\em branch decomposition} 
\begin{equation}\label{e.branchdecomp}
 ( \gamma_{[0,a_1]} , \gamma_{[a_1,a_2]} , \dots , \gamma_{[a_{r-1},a_{r}]} ) \, ,
\end{equation}
with $r = r(\gamma)$ the number of branches of $\gamma$. The branch heights form a strictly decreasing sequence.
    
When $k$ is even, $\gamma_{[a_k,a_{k+1}]}$ is a bridge after a suitable translation.
When this index is odd, it is the  vertical reflection of $\gamma_{[a_k,a_{k+1}]}$ that is the translate of a bridge. 
More formally, write $\tau$ for the reflection with respect to the horizontal axis. 
	The {\em bridge decomposition} of $\gamma$ is formed by applying $\tau$ to every second component of  $\gamma$'s branch decomposition and suitably translating each resulting component so that it is a bridge: that is, the new decomposition is 
\begin{equation}\label{e.bridgedecomp}	
 ( \gamma_{[0,a_1]} ,  \tau(\gamma_{[a_1,a_2]} - \gamma_{a_2}) , \gamma_{[a_2,a_3]} - \gamma_{a_2} , \dots ) \, .
 \end{equation}
	We define $\Psi(\gamma)$ to be the bridge decomposition of $\gamma$. Since the transformation that specifies this map is invertible, $\Psi$ is injective, and Lemma~\ref{l.hsw} is proved. \qed

\subsection{A rough guide to how polygon insertion will modify the Hammersley-Welsh argument}\label{sec:howpoly}
The injective map $\Psi$ constructed in Lemma~\ref{l.hsw}  sends a half-space walk to its bridge decomposition.
Since this map  is so vital to the proof of Theorem~\ref{thm:HW},  we call $\Psi$ the {\em Hammersley-Welsh} map.
It will be by perturbing the construction of $\Psi$ that we prove Theorems~\ref{thm:HW1} and~\ref{thm:HW2}, and thus we emphasise its role. Note, however, that the proof of Theorem~\ref{thm:HW} amounts to considering the composition $\concat \circ \Psi : \HSW_n \to \SAB_n$, where $\concat$ is the concatenation map from Lemma~\ref{l.bridgeproduct}'s proof, and noting that the number of arrows incoming to any point in the range of this map is bounded above by $\exp ( C n^{1/2} )$. (We have called Hammersley-Welsh's proof  an unfolding argument because  $\concat \circ \Psi$ unfolds walks into bridges.) 
 
Here, we outline the improvement of the Hammersley-Welsh bound. Fix $\eps > 0$ to be the desired reduction of the $1/2$ exponent of \eqref{eq:HWimp1}; we will see below that the improvement is limited, as $\eps$ needs to be smaller than some threshold.

A simple argument will permit us to restrict the domain of $\Psi$ to a special class of half-space walks.
This set will be denoted (for a reason that we explain shortly) by $\SHSW_n$. In essence, a member $\gamma \in \SHSW_n$ is an element of $\HSW_n$ that   
\begin{enumerate}
\item has {\em many branches}: the branch, and thus also the bridge, decomposition of $\gamma$ has at least an order of~$n^{1/2 - \eps}$ terms;
\item and is {\em horizontally confined}: 
$\gamma$ is contained in a vertical strip of width of order  $n^{1/2 + \eps}$. 
\end{enumerate}
Our goal will be achieved if we prove that $|\SHSW_n|\mu^{-n} \leq \exp(C n^{1/2 - \eps})$ for some constant~$C$. 
Indeed, the bound $|\HSW_n \setminus \SHSW_n|\mu^{-n}\leq \exp(C n^{1/2 - \eps})$ is straightforward: see Lemmas~\ref{lem:HSW_CHSW} 
and~\ref{lem:few_branches}. These two bounds will then imply Theorem~\ref{thm:HW1}.

In fact, and as we will now explain heuristically, 
we are able to obtain a bound on $| \SHSW_n \vert \mu^{-n}$ that is stronger than needed, namely
\begin{align}\label{eq:Alan_new}
 \vert \SHSW_n \vert \cdot \mu^{-n} \leq \exp \big\{ - C n^{1/2 - \eps} \log n \big\} \, .
\end{align}
To begin explaining how we will reach this inference, fix $r \in \N$ to be a given value whose order is at least $n^{1/2-\eps}$. 
Set $\SHSW_n(r) = \{\gamma \in \SHSW_n: \text{ $\gamma$ has $r$ branches} \}$.  
We will explain how to derive the counterpart to \eqref{eq:Alan_new} on whose left-hand side $\SHSW_n$ is replaced by $\SHSW_n(r)$. 
The actual bound \eqref{eq:Alan_new} may then be obtained by summing over the concerned values $r$, 
provided that the positive constant $C$ is reduced suitably.

%

Consider the restriction of $\Psi$ to $\SHSW_n(r)$; its image is a subset of $\bridgeprod_{n,r}$.
How may we seek to improve~$\Psi$ for our purpose, when this map is already injective? The answer lies in modifying $\Psi$ into a multi-valued map
with a large number of outgoing arrows from each point of $\SHSW_n(r)$, without significantly compromising the map's injectivity.
\remove{between these same sets so that there are a helpful number of outgoing arrows from any domain point, without significantly compromising this map's injectivity. }

To survey options for introducing outgoing arrows from a given domain point, denote such a point by $\gamma \in \SHSW_n(r)$.
\remove{The sequence of branch heights of $\gamma$ is strictly decreasing, and has   $r$ terms.} 
Due to the high number of branches of $\gamma$, 
it is in essence true  that at least one-half of them have height exceeding  $n^{1/2 - \eps}$.
These will be called the {\em tall branches}. Since the vertical intervals occupied by consecutive branches form  a decreasing sequence under containment, we see that $\gamma$'s tall branches all cross a given vertical interval of length~$n^{1/2 - \eps}$.
Moreover, as $\gamma$ is horizontally confined, every tall branch crosses between its upper and its lower side 
some given translate~$R$ of the rectangle  $[-n^{1/2 + \eps}, n^{1/2 + \eps}] \times [-n^{1/2 - \eps}, n^{1/2 - \eps}]$, 
all the while remaining in~$R$.

Thus, $\gamma$ is necessarily a rather dense object: 
the notation $\raisebox{-0.04em}{\includegraphics{ww.pdf}}$ is intended to evoke $\gamma$'s repeated up-down movement at close quarters. 

	In particular, the tall branches of $\gamma$ often pass close to each other inside~$R$. 
	Say that a point of $\gamma$ in $R$ is a {\em  near self-touch} 
	if it is the rightmost point of a tall branch at a given height among the $y$-coordinates assumed by $R$, 
	and if the horizontal interval of length $n^{2\eps}$ to its right intersects another tall branch of $\gamma$. 
	Each of the at least $r/2$ tall branches of $\gamma$
	contains one rightmost point for each of the $n^{1/2-\eps}$ heights of $R$. 
	Due to the limited volume of $R$, most of these rightmost points may be shown to be near self-touches. 
	In conclusion, $\gamma$ has of order $r\, n^{1/2-\eps}$ near self-touches.
	
%

We will alter the definition of $\Psi$ in a way that seeks to exploit the abundance of near self-touches among elements in its domain. 
%
Our alteration of~$\Psi$, which we will call~$\Phi$, entails the insertion of certain self-avoiding polygons around near self-touches. 
The polygons to be inserted will be drawn from a certain class  $\WSAP^u_m$ of {\em wide} polygons; here,  $u$ and $m$ are parameters. The set $\WSAP^u_m$  is a collection of self-avoiding polygons of length~$m$ each of which has, roughly expressed, width at least $u$.
A choice of $u$ will be made of order $n^{\delta}$, where $\delta > 0$ is a parameter chosen so that the resulting choice  equals or slightly exceeds the horizontal near self-touch distance $n^{2\eps}$; 
while~$m$ will be a given value that is at most $n^{2\delta}$. The value $\delta$ should be at least $2\eps$ and, in this overview, we in fact set $\delta = 2\eps$. 

In a first attempt to modify the definition of $\Psi$, suppose that we attempt to insert a single polygon onto a half-space walk in the domain of~$\Psi$, and then take the bridge decomposition of this modified walk.
That is, suppose that we define $\Phi$ to be a multi-valued map whose domain is the product space $\SHSW_n(r) \times \WSAP^u_m$.

	Fix a domain point $(\gamma, P)\in \SHSW_n(r) \times \WSAP^u_m$ 
    and consider the decomposition $b_1,\dots, b_r$ of $\gamma$ into branches. In this paragraph, we will abusively treat each $b_i$ as if it were bridge, even though alternate vertical reflection is in fact needed to ensure this. 
   The set $\Psi(\gamma,P)$ of bridge-lists associated with $(\gamma,P)$ will be obtained by inserting $P$ onto a given tall branch of~$\gamma$.
   Let $b_j$ denote this branch; we choose the index $j$ so that the cardinality of the collection of near self-touches of $\gamma$ along $b_j$ has the typical order of $n^{1/2-\eps}$. 
    Let $\gamma_t$ denote a generic element in this collection of near self-touches. 
    Define $\tilde{b}_j$ to be the bridge obtained by gluing the polygon $P$ to the right of $b_j$, at the level of the point $\gamma_t$. 
    The procedure used to glue a polygon to a bridge is called the {\em Madras join} and will be reviewed in the next section. 
    By definition of the Madras join, the result $\tilde{b}_j$ has length  $|b_j| +m + 16$ and, excepting a few marginal examples, is a bridge with the same height as $b_j$. 
    Consider now the list $(b_1,\dots, \tilde b_j, \dots, b_r)$ which is given by $\Psi(\gamma)$ except that $b_j$ is replaced by $\tilde b_j$; this is a bridge-list of total length $n + m+16$.
    Then $\Phi(\gamma,P)$ is the set of all bridge-lists obtained in this manner, as $\gamma_t$ ranges over points of near self-touch for $\gamma$ that lie in $b_j$. Thus, $\Phi$ is a multi-valued map from $\SHSW_n(r) \times \WSAP^u_m$ to $ \bridgeprod_{n+m+16,r}$.    

\begin{figure}
	\begin{center}
	   	\includegraphics[width = 0.7\textwidth]{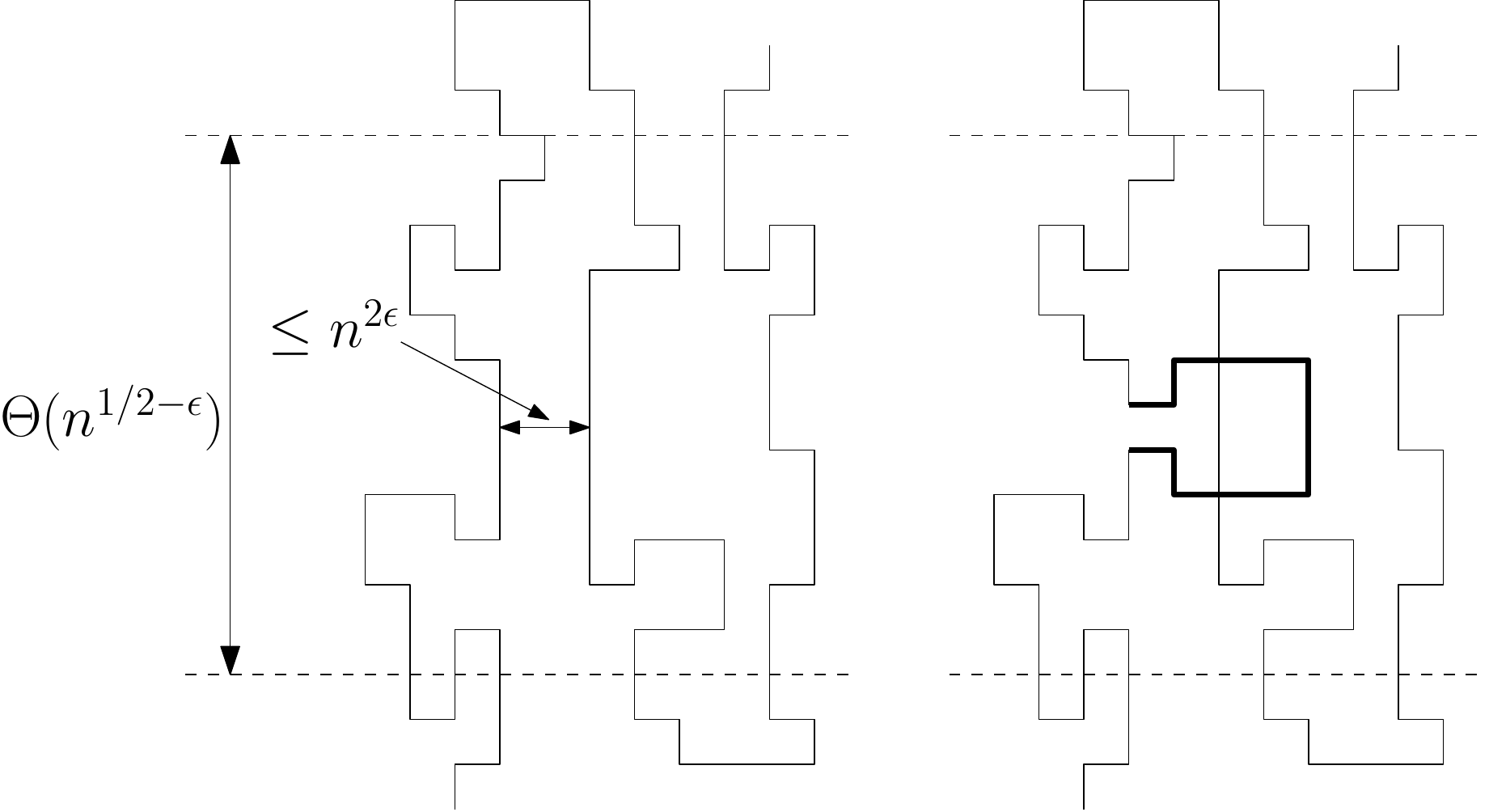}
		\caption{By means of a local deformation on its left, a square-shaped polygon $P$ is inserted into a branch walk $\gamma$
		 at a height where this branch nearly touches the branch to its right. 
		 The implanted copy of $P$ intersects the branch on the right. 
		 The edges introduced into the walk are shown in bold in the right sketch. }
	\label{fig:polygoninsertion}
	\end{center}
\end{figure}

We now evaluate how successful this modification~$\Phi$ of the Hammersley-Welsh map~$\Psi$ is in finding an improved upper bound on $\vert \SHSW_n(r) \vert$.  
Lemma~\ref{l.mvm} and Lemma~\ref{l.bridgeproduct}(1) imply that 
\begin{align}\label{eq:EMM}
	\vert \SHSW_n(r) \vert \cdot \vert \WSAP^u_m \vert 
	\leq \frac{\EMM}{\emm} \vert \bridgeprod_{n+m,r} \vert
	\leq \frac{\EMM}{\emm} \cdot  (n+m)^{r} \mu^{n+m},
\end{align} 
where $\EMM$ is  the maximum number of arrows incoming to a given range point; and $\emm$, the minimum number of arrows outgoing from any domain point, is, in light of the preceding paragraph, at least $n^{1/2 - \eps}$ up to a constant factor which we here neglect.

For $\Phi$ to produce a meaningful result, the constant $\EMM$ associated with it needs to be bounded from above. We next explain how this is done.

Let $\bar{b} = \big(b_1,\dots,b_r \big) \in \bridgeprod_{n+m +16,r}$ be a point in the image set of $\Phi$. 
Consider the walk resulting from the concatenation of the elements $b_1,\dots, b_r$ after vertical reflection of every other bridge in this sequence. 
Vitally, this result is not self-avoiding: the choice of $P$ and of the point $\gamma_t$ of near self-touch 
have been made so that the branch to which $P$ was added intersects its right-adjacent counterpart. 
We may thus infer the branch to which $P$ was added (though in the present case, we know it: it is $b_j$); and, in view of the fixed length of~$P$, we further in essence infer the location of the joining of $P$ to that branch. Only a finite number~$C$ of pre-images of $\bar{b}$ are then possible, where the value $C$ accounts for possible microscopic variations in the surgery used to join $P$ to the branch of $\gamma$. Thus, \eqref{eq:EMM} becomes
\begin{align}\label{eq:EMM2}
	\vert \SHSW_n(r) \vert \cdot \mu^{-n} 
	\, \leq \, C  n^{-(1/2-\eps)} \cdot  (n+m)^{r} \cdot  \big(\vert \WSAP^u_m \vert\cdot \mu^{-m}\big)^{-1} \, .
\end{align} 
This should be compared to the naive bound involving the original map $\Psi$, which yields 
\begin{align}\label{eq:EMM0}
	\vert \SHSW_n(r) \vert \cdot \mu^{-n} 
	\leq n^{r} \, .
\end{align} 
The difference between the right-hand factors of the last two inequalities represents the change in the outcome due to the insertion of a single polygon. 
The order $n^{-(1/2-\eps)} \ll 1$ benefit 
due to polygon placement entropy is possibly compromised by the polygon scarcity cost $(\vert \WSAP^u_m \vert \mu^{-m})^{-1}$ (it will transpire that the loss due to $n^{r}$ being replaced by $(n+m)^r$ is insignificant).
To understand the role of the term $(\vert \WSAP^u_m \vert \mu^{-m})^{-1}$, note that, since $\WSAP^u_m$ is a collection of polygons of length~$m$, we have $\vert \WSAP^u_m \vert \leq \vert \SAP_m \vert \leq \mu^m$: see~\cite[(3.2.5)]{MadSla93}. Insofar as this inequality on $\vert \WSAP^u_m \vert$ fails to be sharp, the resulting polygon scarcity constitutes an opposing force that undoes some of the benefit achieved by modifying~$\Psi$ to~$\Phi$. Thus, it is an upcoming challenge to demonstrate a lower bound on the number of wide polygons. If we posit that, for some constant $\alpha > 0$, there exist $u$ of order $n^{2\eps}$ and $m \leq u^2$ for which $\vert \WSAP^u_m \vert \mu^{-m} \geq u^{-\alpha}$,
then \eqref{eq:EMM2} would become 
\begin{align}\label{eq:EMM3}
	\vert \SHSW_n(r) \vert \cdot \mu^{-n} 
	\leq C n^{-(1/2 - \eps) + 2\alpha \eps} \cdot  (n+m)^{r} \, ,
\end{align} 
where the influence of the modified argument is transmitted again through the first term on the right-hand side, namely \add{$C{n^{\eps(1+2\alpha) - 1/2}}$}. 
Given the information that suitable wide polygons are plentiful in the sense of the parameter~$\alpha$, the modification is seen to return a benefit when $\eps (2\alpha +1)<1/2$.
Even when this condition is met, the improvement of $n^{\eps(1+2\alpha) - 1/2}$ of \eqref{eq:EMM3} compared to \eqref{eq:EMM0} is only very modest: after all, $r$ has order at least $n^{1/2 - \eps}$. 


To obtain a greater improvement, we have little choice but to iterate our procedure: instead of adding one polygon, we will add as many as needed to reduce the right-hand side of \eqref{eq:EMM2} below $\exp(n^{1/2 - \eps})$.
Suppose that we add $K = \kappa r$ polygons independently and uniformly into admissible near-touch slots; note that $\kappa > 0$ is a parameter indicating the order of the number of polygons inserted per tall branch. Adding a polygon to a near-touch slot may block nearby near self-touches from receiving polygons. 
However, only at most order $n^{2\eps}$ other near self-touches are thus affected. 
Indeed, since there is at least an order of $n^{1-2\eps}$ near self-touches in total, the assumption that $K = o(n^{1-4\eps})$, which in practice certainly permits $\kappa$ to be chosen to be of unit order, entails that
there are of order $n^{K(1-2\eps)}$ sets of near self-touches for the $K$ polygons to be added. 

In this light, we reconsider $\Phi$ to be a multi-valued map from $\SHSW_n(r) \times (\WSAP^u_m)^K$ to $ \bridgeprod_{n+K(m+16),r}$
obtained by adding $K$ polygons drawn from $\WSAP^u_m$ to the tall branches of $\gamma \in \SHSW_n(r)$. 
The reasoning offered for the case of the insertion of a single polygon now shows that
$\emm \geq n^{K(1-2\eps)}$ and $\EMM \leq C^K$. Assuming again that  $\vert \WSAP^u_m \vert \mu^{-m} \geq u^{-\alpha}$, we find that 
\begin{align*}
	\vert \SHSW_n(r) \vert \cdot \mu^{-n} 
	&\leq \big[C n^{-(1/2 - \eps) + 2\alpha \eps} \big]^{K}\cdot (n+Km)^{r} \\
	&\leq \exp\bigg[\Big(-\kappa\big( 1/2 - \eps(1+2\alpha)  \big) +1 \Big)r\log n + r \log(2C^\kappa) \bigg] \, ,
\end{align*} 
where the latter bound depends on the practically harmless replacement of $(n+Km)^r$ by $n^r$.
The choice of $\eps >0$ so that $\eps(1+2\alpha) <1/2$, and the further choice of $\kappa$ to be a unit-order quantity satisfying 
$\kappa > \big( 1/2 - \eps(1+2\alpha) \big)^{-1}$,
ensures that the factor multiplying $r\log n$ in the above may be treated as a negative constant. 
That $r$ is of order $n^{1/2 - \eps}$ yields the inference 	
\begin{align*}
 \vert \SHSW_n(r) \vert \cdot \mu^{-n} \leq \exp \big\{ - C n^{1/2 - \eps} \log n \big\} \, .
\end{align*}
As already mentioned, \eqref{eq:Alan_new} follows from the above by summing over $r$. 


These heuristic considerations offer a  road map for the rigorous argument for Theorem~\ref{thm:HW1} which will be presented in Sections~\ref{sec:proofsone} and~\ref{sec:proofstwo} (with the variation needed for Theorem~\ref{thm:HW2} explained in Appendix~\ref{sec:madras_join_hex}). 
The proof hinges on the abundance of polygons, which is to say, on a power-law lower bound on $\vert \WSAP^u_m \vert \mu^{-m}$. 
Such bounds will be studied in Section~\ref{sec:SAP}; they are proved to hold subsequentially for the square lattice and for all $u$ (and well-chosen values of $m$) for the hexagonal lattice.

\section{Self-avoiding polygons: counting and joining}\label{sec:SAP}

In this section, we present two tools which are needed to rigorously implement the plan that has just been sketched. In the preceding discussion, we 
invoked an assertion of polygon abundance, namely $|\WSAP_m^u| \geq u^{-\alpha}$.  In Section~\ref{sec:WSAP}, we will present a precise definition of the wide polygon set $\WSAP_m^u$. Propositions~\ref{prop:many_poly_subseq} and~\ref{prop:many_poly_hex} are our rigorous assertions of polygon plenitude. 
We also review in Section~\ref{sec:MJ} Madras' joining technique which allows polygons to be attached to the branches of $\gamma$.  

\subsection{Wide polygons}\label{sec:WSAP}
%

For a polygon $p \in \SAP_m$, the {\em width} of $p$ is 
\newcommand{\width}{\mathsf{w}}
\newcommand{\lwidth}{\mathsf{Lw}}
\begin{align*}
	\width(p) := x_{\max}(P) - x_{\min}(P) \, , 
\end{align*}
where $P$ is any representative of the equivalence class $p$.
The {\em line-width} of $p$ is 
\begin{align*}
	\lwidth(p) := \max \{ \,  x_{\max}(P \cap \ell_y ) - x_{\min}(P \cap\ell_y):\, y \in \bbZ\text{ for which } P \cap \ell_y\ne \emptyset \,  \} \, ,
\end{align*}
where $\ell_y$ is the horizontal line $\bbZ\times\{y\}$. Note that $\width(p) \geq \lwidth(p)$;   the former  quantity may much exceed  the latter.

Let $m \in 2\N$ and $u \in \bbN$. We define the {\em wide polygon} set
$$
\WSAP_m^{u} \, = \, \{ \, p \in \SAP_m: 
 \lwidth(p) \geq u \, , \, \height(p) \leq 16 u \, \} \, ,
$$
where the height notation~$\height$ from  Section~\ref{s.notation} has been extended to polygons.

For practical purposes, we will choose one particular rooting for wide polygons.
Indeed, the walk $P$ corresponding to any $p \in \WSAP_m^u$ is such that 
\begin{itemize}  
    \item 
    the line-width of $p$ is realized by $P$ at height $0$, and this is the lowest level where it is so realized: 
    \begin{align*}
    	\lwidth(p) 
		&= x_{\max}(P \cap \ell_0) - x_{\min}(P \cap \ell_0)\\
		&> \max\{ x_{\max}(P \cap \ell_y) - x_{\min}(P \cap \ell_y):\, y < 0 \}\quad \text{ and}
    \end{align*}
    \item the origin is the leftmost point of $P$ at height $0$, i.e.~$0 = x_{\min}\big( P \cap \ell_0 \big)$.
\end{itemize}
By convention, the walk that  represents $p$ begins and ends at the origin, and traverses the edges that comprise $P$ in a counterclockwise manner, so that      
 the 
 interior of $P$ is on the left during the walk.
 The correspondence between polygon and walk offered by this convention is a bijection, because it specifies a representative of the equivalence class in question, which amounts to distinguishing a vertex on the polygon; and it then selects one of the two possible orientations of the resulting walk. 

The abundance of wide polygons is important for our application. 
We verify this abundance  by two different means, one for $\Z^2$ and the other for $\bbH$. The latter inference is stronger. 
It is in deriving this inference that the parafermionic observable available for $\bbH$ is used during the proof of Theorem~\ref{thm:HW2}.

\begin{proposition}[Subsequential abundance of polygons for $\Z^2$]\label{prop:many_poly_subseq}
	For any $\alpha > 29$, there exists an infinite set $A \subset \bbN$ such that, 
	for all $u \in A$ there exists $m \in \llbracket 4u+4,4u^2+4 \rrbracket$ for which 
	$$
		\big|\WSAP_m^{u}\big|\geq u^{-\alpha}\mu^m \, .
	$$
\end{proposition}


Write $\WSAP_m^{u}(\bbH)$ for the set of self-avoiding polygons on $\bbH$, with the same properties as those defining wide polygons on $\bbZ^2$. 

\begin{proposition}[Abundance of polygons for $\bbH$]\label{prop:many_poly_hex}
	There exist positive constants $C$ and $c_0$ such that, for all $u \in \bbN$, we may find $m\le c_0 u^2$ for which
	$$
		|\WSAP_m^{u}(\bbH)|\geq C\, u^{-10}\mu(\bbH)^m \, .
      $$
\end{proposition}

Proposition~\ref{prop:many_poly_hex} is a principal element of the variation of the proof of Theorem~\ref{thm:HW1} that is needed to obtain Theorem~\ref{thm:HW2}. Its proof appears after the derivation of Theorem~\ref{thm:HW1} is completed,
 in  Section~\ref{sec:observable}.


The rest of Section~\ref{sec:WSAP} is dedicated to proving Proposition~\ref{prop:many_poly_subseq}.
We start with a lemma that bounds the number of wide polygons  in terms of the number of bridges. 

\begin{lemma}\label{lem:WSAP_SAB}
    For any $m \in \N$, there exists $u \in \N \cap [m^{1/2},m]$ such that
 $$
     \big\vert \WSAP^u_{4m + 4} \big\vert \, \geq \, 
      (\log 2)^2  3^{-3}  2^{-14} \,   \big( \log u \big)^{-2} u^{-21} \big\vert \SAB_m \big\vert^4 \, .
    $$  
    \end{lemma}
\noindent{\bf Proof.}
	Let $m \in \N$. Specifying the diameter $\mathsf{diam}(p)$ of a polygon $p \in \SAP_{2m+2}$ to be $\max\{\height(p),\width(p)\}$, it is easily seen that 
	$m^{1/2} \leq \mathsf{diam}(p) \leq m$. 
	Let $u \in[m^{1/2},m]$ be a number of the form $2^{-j}m$ such that the set 
	$$
	\{ \, p \in \SAP_{2m+2}: \, u/2 \leq  \mathsf{diam}(p) \leq u \, \}
	$$ 
	has maximal cardinality among such sets. Note that, for now, we permit the parameter $u$ to be real, although it is a natural number in the lemma's statement; we will resolve the discrepancy at the end of the proof. 
	
		At least half of the polygons in the just displayed set are at least as wide as they are high; 
	(indeed, any polygon that does not satisfy this condition is the right-angled rotation of one that does). 
Writing $\overline{\WSAP}^u_{2m+2}$ for the set of polygons $p \in \SAP_{2m+2}$ with $\width(p) \geq u/2$ and  $\height(p) \leq u$, we thus see that \begin{equation}\label{eq:diam}
	\big|\overline \WSAP^u_{2m+2} \big|  \geq \frac{1}{2\log_2 m}|\SAP_{2m+2}| \, .
	\end{equation}
	When  $p \in \overline \WSAP^u_{2m+2}$ is depicted as a walk, we may specify the west-south vertex of~$p$, 
	$\mathsf{WS}(p)$, (and its east-north vertex $\mathsf{EN}(p)$) 
	to be the lowest among the leftmost (and the highest among the rightmost) elements of~$\Z^2$ visited by the walk~$p$. 
	Clearly each $p\in \overline \WSAP^u_{2m+2}$ has a unique representation such that $\mathsf{WS}(p) = 0$.
	Then $x\big(\mathsf{EN}(p)\big)= \width(p)$. 
	
\begin{figure}
    	\begin{center}
    	   	\includegraphics[width = 1\textwidth]{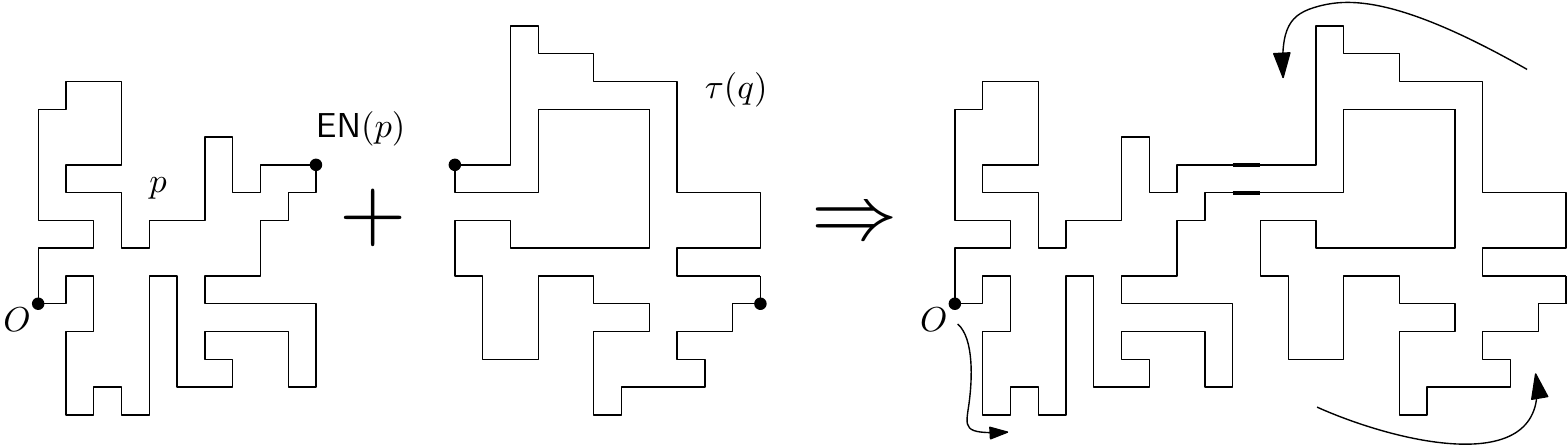}
    		\caption{Two elements $p,q$ in $\overline\WSAP^u_{2m+2}$ with $y\big(\mathsf{EN}(p)\big) = y\big(\mathsf{EN}(q)\big)$ are joined to form $J\big( p,\tau(q) \big) \in \WSAP^u_{4m+4}$. The dot at $O$ and the curved arrows show the root and orientation of the joined polygon that is dictated by our convention.}
    	\label{fig:WSAP}
    	\end{center}
    \end{figure}

	Let $p,q \in \overline \WSAP^u_{2m+2}$ be represented by walks such that $\mathsf{WS}(p)=\mathsf{WS}(q)= 0$
	and be such that $y\big(\mathsf{EN}(p) \big) = y\big(\mathsf{EN}(q)\big)$.
We will join the polygon $p$ to a reflection $\tau(q)$ of $q$
to form a new polygon $J\big( p,\tau(q) \big)$ about which we will claim that $J(p,q) \in \WSAP^u_{4m + 4}$. The joining operation~$J$ does the same job as the Madras join which we will review in Section~\ref{sec:MJ}, but the special geometry arising from the assumption that $y\big(\mathsf{EN}(p) \big) = y\big(\mathsf{EN}(q)\big)$ permits the use of a simpler technique in the present case. 
Figure~\ref{fig:WSAP} illustrates the operation $J$.

The polygon $q$ may be reflected in the vertical axis, via the reflection~$\tau$, and then horizontally translated to such a position that its vertex in correspondence with~$\mathsf{EN}(q)$ is  to be found one unit directly to the right of $\mathsf{EN}(p)$. The union of the edges that comprise $p$ and this translation of $\tau(q)$ intersects the four edges of the plaquette whose north-west corner is $\mathsf{EN}(p)$ along the west and east sides of the plaquette. By replacing these two edges in the union by the plaquette's north and south sides, a polygon of length $2(2m+2)$ results. 
This polygon is $J\big(p,\tau(q) \big)$. 

We now show that $J\big(p,\tau(q) \big)$ belongs to~$\WSAP^u_{4m +4}$.  
The height of this polygon is at most $\height(p) + \height(q) \leq 2u$.  
Moreover,  $\mathsf{WS}(J\big(p,\tau(q) \big))  =\mathsf{WS}(p) = 0$
and $J\big(p,\tau(q) \big)$  realizes its line-width $\lwidth(J\big(p,\tau(q) \big))=  \width(p)+ \width(q) +1 > u$ at height zero.
Finally, $O$ is the leftmost point of $J\big(p,\tau(q) \big)$ at height $0$; hence $J\big(p,\tau(q) \big)$ is rooted according to our convention for wide polygons. 
The claim that $J\big(p,\tau(q) \big) \in \WSAP^u_{4m + 4}$ has been confirmed.

Note that  $J\big(p,\tau(q) \big)$ determines $(p,q) \in \WSAP^u_{2m+2}$. Indeed,
there is a unique vertical line whose coordinate has the form $k + 1/2$ for some $k \in \Z$ that cuts $J\big(p,\tau(q) \big)$ through a consecutive pair of horizontal edges, and for which the replacement in  $J\big(p,\tau(q) \big)$ of the horizontal with the vertical sides of the associated plaquette  results in two polygons, of which the one whose vertex set contains the origin has $2m + 2$ edges (and, as it happens, so does the other). The plaquette in question was used in the surgical formation of  $J\big(p,\tau(q) \big)$, and the operation just described undoes this surgery so that $p$ and $\tau(q)$ are recovered.	
This injectivity on the part of $J$
implies that
$$
 \big\vert \WSAP^u_{4m + 4} \big\vert
		\geq \sum_{j = - \lceil u \rceil}^{\lfloor u \rfloor} \big\vert \big\{ p\in \overline \WSAP^u_{2m+2}:\, y(\mathsf{ne}(p)) = j \big\} \big\vert^2 \, ,
$$
where it is since $u$ may not be an integer that use is made of rounding.
This right-hand side is seen to be  at least $(2u + 1)^{-1} \big\vert \overline\WSAP^u_{2m+2} \big\vert^2$ by means of  the Cauchy-Schwarz inequality. Applying~(\ref{eq:diam}) and Lemma~\ref{l.sapsab} with $n =m$, we see that
\begin{align*}
 \big\vert \WSAP^u_{4m + 4} \big\vert 
 \geq  (2u + 1)^{-1} \bigg(\frac{\big\vert \SAB_m \big\vert^2}{8(2m+1)m(m+1)^3 \log_2 m} \bigg)^2 \, .
\end{align*}
Since $u^2 \geq m  \geq 1$, a short algebraic manipulation leads to 
\begin{equation}\label{e.ubound}
 \big\vert \WSAP^u_{4m + 4} \big\vert \, \geq \, 
 (\log 2)^2  3^{-3}  2^{-14}   \big( \log u \big)^{-2}      u^{-21} \big\vert \SAB_m \big\vert^4 \, . 
\end{equation}
This is precisely the bound asserted by Lemma~\ref{lem:WSAP_SAB}. To complete the lemma's proof, it remains merely to argue that our parameter $u \in [m^{1/2},m]$ may in fact be chosen to satisfy $u \in \N$. For this purpose, we consider $\lceil u \rceil \in [m^{1/2},m]$, and note that, by the definition of the wide polygon set,  
$\WSAP^u_{4m + 4} \subseteq \WSAP^{\lceil u \rceil}_{4m + 4}$.
Thus,~(\ref{e.ubound}) holds when the replacement $u \to \lceil u \rceil$ is made.  \qed

\noindent{\bf Proof of Proposition~\ref{prop:many_poly_subseq}.}
As~\cite[Corollary~$3.1.8$]{MadSla13} reviews, a variation of the Hammersley-Welsh argument for Theorem~\ref{thm:HW} establishes that the partition function for bridges diverges at its radius of convergence~$\mu^{-1}$:
\begin{align*}
	\sum_{m = 1}^\infty \mu^{-m}|\SAB_m| = \infty \, .
\end{align*}
Let $\eps > 0$. There are thus infinitely many values of $m \in \N$ such that $|\SAB_m| \geq m^{-1-\eps}\mu^{m}$.
Inserting this bound into Lemma~\ref{lem:WSAP_SAB},
we obtain infinitely many pairs $(m,u)$ with $u \in \N \cap [m^{1/2},m]$ such that 
\begin{align*}
	\big|\WSAP_{4m+4}^u\big| 
	\,	\geq \, c_0 \, u^{-21}(\log u)^{-2}  m^{-4-4\eps} \, \mu^{4m}
	\,	\geq \,  c_0 \, u^{-29 -8\eps} (\log u)^{-2} \, \mu^{4m} \, ,
\end{align*}
where $c_0 = (\log 2)^2 3^{-3} 2^{-14}$.
Since $\eps > 0$ is arbitrary, the proposition follows by relabelling~$m$. \qed

\subsection{Madras joining}\label{sec:MJ}

In~\cite{Mad95}, Madras introduced a technique by which a pair of polygons may be joined to form a new polygon. We will be using his technique in order to attach polygons to the side of bridges or their vertical reflections.
We begin by briefly reviewing his procedure. Then we state and prove Lemma~\ref{l.madrascondition}. This result gathers together all properties of the Madras join which we will use. 
This section's hexagonal lattice counterpart is postponed to Appendix~\ref{sec:madras_join_hex}.

The eight parts of Lemma~\ref{l.madrascondition} are all needed in this article. However, the reader may lose little  by 
perusing the lemma briefly and omitting its proof; the core of the argument may even be followed by merely
interpreting the Madras join as a procedure that inserts a polygon at a given height on the right of a walk while effecting a bounded change of length near a distinguished plaquette, called the junction plaquette, about which the join is made. 

The reader who wishes to follow details precisely is encouraged, in reading the lemma and its proof, to consult either~\cite{Mad95} or the article~\cite{Ham17} which has also made use of Madras' technique. Section~$4.1$ of the latter article contains a detailed review of the technique,  including a depiction in Figure~$1$ of the various cases of local deformation used in the operation.

\subsubsection{An overview of Madras' joining technique}
%
%
%
For any given polygon $P$, 
the integer interval $\llbracket y_{\min}(P) , y_{\max}(P)  \rrbracket$ comprises the vertical coordinates of vertices in~$\Z^2$ visited by the polygon. 
Suppose given a pair of polygons $P$ and $Q$ such that 
\begin{equation}\label{e.joiningcondition}
\llbracket y_{\min}(P)-1 , y_{\max}(P) +1  \rrbracket \cap \llbracket y_{\min}(Q) -1 , y_{\max}(Q) + 1  \rrbracket \neq \emptyset \, .
\end{equation}
 Let $n,m \in 2\N$ denote the respective lengths of $P$ and $Q$.
The {\em Madras join} of $P$ and $Q$ is a polygon, which we will denote by $M(P,Q)$, formed by attaching $Q$ after horizontal displacement to the right side of $P$.
The procedure for joining begins by shifting $Q$ far to the right of~$P$, and bringing it back, step by step, until there is a vertex in $P$ and a vertex in the translate of $Q$ 
that are either equal or whose displacement is a vector of vertical orientation and of Euclidean norm at most two. We then define the vertex~$Y$ to be the vertex~$z$ of maximal $y$-coordinate among those for which $\big\{ z - e_2,z,z+e_2 \big\}$ intersects $P$ and the translate of $Q$. A local surgery is performed on the two polygons in the neighbourhood of~$Y$. The local geometry of $P$ determines the form of the surgery on $P$ according to a division into several cases; and similarly for the translate of $Q$. The two local surgeries are performed
in such a way that the modified 
translate of $Q$ may be  translated a certain few steps to right into such a position that there is a plaquette whose intersection with the locally modified copy of $P$ consists of its left vertical edge, and whose intersection with the locally modified translate of $Q$ consists of its right vertical edge.  

With $Q$ indeed relocated in this way,  these two vertical edges are replaced in the collection of edges in $P$ and the new translate of $Q$ by the two horizontal edges in the plaquette.  The Madras join polygon $M(P,Q)$ is defined to be the resulting polygon.  The plaquette involved in the final step of this procedure will be called the junction plaquette associated with the joining of  $P$ and $Q$. 
 
 The local surgeries performed are such that the length of $M(P,Q)$ equals $n + m + 16$, so there is always a net gain of sixteen in the number of edges involved as a result of the joining operation. The collection of edges in the symmetric difference of the edge collections of $M(P,Q)$ and of the union of $P$ with the surgically determined translate of $Q$ may be viewed as an error set arising in surgery.
 This error set contains a bounded number of edges (there are at most twenty of them), all these in the locale of the junction plaquette associated with~$P$ and $Q$.

\subsubsection{Our application of Madras joining, with bridges and polygons}

In the present paper, we will not be joining pairs of polygons, but rather pairs $(\gamma,P)$ where $\gamma$ is a bridge, or the vertical reflection of one, and $P$ is a polygon. In a natural extension of notation, the walk resulting from the join, in which a horizontal translate of $P$ is in essence incorporated into the range of $\gamma$, will be denoted by $M(\gamma,P)$. 

In our applications, this outcome $M(\gamma,P)$ will be a bridge, or the reflection of one: indeed, the first part of the next result ensures this property. Note that the condition~(\ref{e.madrascondition}) on vertical alignment of $\gamma$ and $P$
is slightly stronger than that~(\ref{e.joiningcondition}) which ensures that Madras joining for a pair of polygons may occur.

\begin{lemma}\label{l.madrascondition}
Let $\gamma$ be a bridge and let $P$ be a polygon. Suppose that
\begin{equation}\label{e.madrascondition}
 \llbracket y_{\min}(P) - 3 ,  y_{\max}(P) + 3 \rrbracket \subseteq  \llbracket y_{\min}(\gamma)  ,   y_{\max}(\gamma) \rrbracket \, . 
\end{equation} 
\begin{enumerate}
\item
 The Madras join $M(\gamma,P)$ is well-defined and is a bridge whose endpoints are shared with $\gamma$.
 \item Write $n \in \N$ for the length of  $\gamma$ and $m \in 2\N$ for the length of~$P$.  
 The endpoints of the right vertical edge in the junction plaquette of $M(\gamma,P)$ belong to  $M(\gamma,P)$. 
 The graph distance between them in $M(\gamma,P)$ (henceforth the ``chemical distance'') equals $m + 7$. 
 \item Denote the junction plaquette  of $M(\gamma,P)$  by $J$, 
  and by $M(\gamma,P) \Delta J$ the modification formed from $M(\gamma,P)$ by the removal of the horizontal edges in $J$ and the insertion of the vertical ones. 
Then  $M(\gamma,P) \Delta J$  comprises a walk $\gamma'$ and a polygon $Q$. No right translation of $Q$, including $Q$ itself, intersects $\gamma'$. 
 \item Suppose that $P'$ is a polygon that also satisfies~(\ref{e.madrascondition}) and for which
  $$
  \llbracket y_{\min}(P) - 2 ,  y_{\max}(P) + 2 \rrbracket \cap \llbracket y_{\min}(P') - 2 ,  y_{\max}(P') + 2 \rrbracket = \emptyset \, .
  $$
 Then  $M \big( M(\gamma,P), P'  \big) = M \big( M(\gamma,P'), P  \big)$. 
 \item Let $P + (i,0)$ denote the horizontal translate of $P$ which is locally modified to form $M(\gamma,P)$.
 Let $\southeast$ (or $\northeast$) denote 
the rightmost among the lowest  (or the highest)
vertices in $P + (i,0)$. The {\em right-side} of $P + (i,0)$
is the collection of vertices in $P + (i,0)$
encountered on a counterclockwise journey between $\southeast$ and $\northeast$, including the journey's endpoints. Every element of the right-side of $P + (i,0)$ lies in~$M(\gamma,P)$. 
\item All but at most two edges in $\gamma$ lie in $M(\gamma,P)$. 
\item The walk $M(\gamma,P)$ is the concatenation of three subwalks. The first is an initial subwalk of $\gamma$; the second has length either $m + 17$ or $m + 18$; and the third is a final subwalk of $\gamma$.
 \item There are at most four bridge-polygon pairs $(\gamma_0,P_0)$ whose elements' respective lengths are shared by $(\gamma,P)$ 
	and for which $M (\gamma_0 , P_0 ) = M(\gamma, P)$, with the junction plaquette for $(\gamma_0,P_0)$ equalling that for~$(\gamma,P)$.
\end{enumerate}
 If $\gamma$ is instead supposed to be the vertical reflection of a bridge, then the same statements hold, except that  $M(\gamma,P)$ is a vertically reflected bridge, as is $\gamma_0$ in the eighth part.
\end{lemma}
\noindent{\bf Proof: (1).} 
The $y$-coordinate of the vertex $Y$ in Madras' construction belongs to  $\llbracket y_{\min}(P) - 1,  y_{\max}(P) + 1 \rrbracket$.
Any endpoint of an edge that is added or removed during surgery has $y$-coordinate equal to $y(Y)-1$, $y(Y)$ or $y(Y)+1$. 
Thus the $y$-coordinates of affected endpoints are seen to lie in the interval $\llbracket y_{\min}(P) - 2, y_{\max}(P) +2 \rrbracket$. By hypothesis, such coordinates lie in $\llbracket y_{\min}(\gamma) +1, y_{\max}(\gamma) -1 \rrbracket$. Since modifications are not made at $y$-coordinate values that are extremal for $\gamma$, the status of $\gamma$ as a bridge or a vertical reflection of one is unaffected by the change $\gamma \to M(\gamma,P)$.   

\noindent{\bf (2).} Because~\cite[Figure~$1$]{Ham17} presents the left, rather than the right, of the two local deformations made during surgery, it is helpful, in regard to the first assertion,  to argue  the claim in which $\gamma$ is a polygon of length~$n$ and $P$ is a walk (the choice of notation is incongruous but it is adopted only temporarily);  and where it is instead claimed that the chemical distance along $M(\gamma,P)$ between the two endpoints of the left vertical edge in the junction plaquette equals $n+7$.
 In each of the seven cases depicted 
in~\cite[Figure~$1$]{Ham17}, consider the two endpoints of the union of the dotted line segments. The vertices so denoted are elements of $\gamma$ whose longer chemical distance around $\gamma$ equals $n - r$ where $r \in \{1,2\}$
is the number of dotted line segments in the case in question. In each case, the sketch in the third column has two consecutive vertical solid edges in the rightmost coordinate. It is the higher of these two edges that is  
the left vertical edge in the junction plaquette. The combined length of the journeys that do not use this edge between the endpoints of this edge and the endpoints of the union of the dotted line segments is seen by inspection to equal $7 + r$. Thus the post-surgical chemical distance around $\gamma$ between these two vertices of the junction plaquette is~$n + 7$.

We want to establish a similar assertion concerning the journey around $P$. The Madras join operation involves the half-circle rotation of the right element $P$ about~$Y$; the identical surgery on the resulting polygon; and then a further half-circle rotation about~$Y$ of the post-surgical polygon.  Since the form of surgery is the same, the preceding analysis is applicable, and we obtain Lemma~\ref{l.madrascondition}(2).

\noindent{\bf (3).} The polygonal component $Q$ of $M(\gamma,P) \Delta J$ lies in the union $C_1 \cup C_2$, where $C_1$ is the part of the strip $\R \times [y(Y) - 1,y(Y)+1]$ on or to the right of the right vertical edge in the junction plaquette $J$, and $C_2$ is the translate of $P$ used in the Madras join. As such, no translation of $Q$ to the right may intersect $\gamma$. The walk component $\gamma'$  of $M(\gamma,P) \Delta J$ lies in the union of $\gamma$ and certain edges in $\R \times [y(Y) - 1,y(Y)+1]$  on or to the left of the left vertical edge in the junction plaquette. Thus, nor may a translation of $Q$ to the right intersect $\gamma'$. 

\noindent{\bf (4).}
Any deformation to $\gamma$ under the surgery $\gamma \to M(\gamma,P)$ occurs in the strip $\R \times  [ y_{\min}(P) - 2 ,  y_{\max}(P) + 2 ]$.
Thus the joining of $P$ and $P'$ to $\gamma$ is a commutative operation.

\noindent{\bf (5).} Similarly to the derivation of Lemma~\ref{l.madrascondition}(2), \cite[Figure~$1$]{Ham17} makes it convenient to argue for a symmetric statement. Suppose instead that the polygon $P$ is on the left, and another polygon is being joined to its right. Introduce the vertices $\mathsf{SW}(P)$ and $\mathsf{NW}(P)$ in accordance with Lemma~\ref{l.madrascondition}(5)'s notation, and specify the left-side of $P$ as the collection of vertices in $P$ encountered on a clockwise journey between this pair of vertices. It is enough to demonstrate that the left-side of $P$ enters the joined output without modification by surgery, because  Lemma~\ref{l.madrascondition}(5) is an assertion symmetrical to this. Consult~\cite[Figure~$1$]{Ham17}. A moment's thought shows that no black edge depicted in the second column may abut a vertex in $P$'s left side. Since the edges removed in surgery, which are the dotted edges in the third column, are a subset of such edges, the sought assertion is demonstrated.

\noindent{\bf (6).} The edges removed from $\gamma$ in the formation of $M(\gamma,P)$ are the dotted edges in the third column of~\cite[Figure~$1$]{Ham17}. There are at most two such edges in each case.

\noindent{\bf (7).} The walk $M(\gamma,P)$ is formed by the removal of either one edge, or two consecutive edges from $\gamma$, and the insertion of a new walk that connects the two endpoints that arise from this removal. The inserted piece has length either $m + 17$ or $m+18$ because $M(\gamma,P)$ has length $n+m+16$.

\noindent{\bf (8).} It is perhaps helpful to begin with an example of how the presurgical data $(\gamma',P')$ is not uniquely specified given the output $M(\gamma,P)$ and the location of the junction plaquette. Suppose that $\gamma$ in a certain locale makes a sequence of down moves until a vertex $z$; then makes one left move; and then continues with down moves. As~\cite[Figure~$1$]{Ham17} indicates, right-attachment of a polygon $P$ with $Y =z$ in Case~$Ib$, or of a polygon $P'$ with $Y = z - e_2$ in Case~$IIb$, leads to outcomes in the third column of the figure that are indistinguishable up to affine shift. In this example, $M(\gamma,P) = M(\gamma,P')$ and the junction plaquette is the same in the two cases. 

More generally, inspection of cases shows that there are never more than two inputs on the left compatible with any given output data. 
The same considerations apply on the right, so that Lemma~\ref{l.madrascondition}(8) is obtained. 
\qed

\subsubsection{Notation for the Madras join in applications}\label{sec:notapp}
Let $\gamma \in \SAW_n$ and $p \in \SAP_m$.
Recall from Section~\ref{sec:WSAP} that the polygon $p$ is depicted as containing the origin as the leftmost point whose vertical coordinate attains $p$'s line-width. 
In our application, we will wire $p$ into $\gamma$
after translating $p$ vertically; after all, it is the entropic benefit among those translations that we seek to exploit. We now set up the notation to be used for the Madras join in this application.

Fix $j \in \bbZ$ with 
$$ 
y_{\min}(\gamma) + \height(p) + 3  \leq j \leq y_{\max}(\gamma) - \height(p) - 3 \, .
$$
We write $\MJ_j(\gamma,p) \in \SAW_{n+m + 16}$
for the Madras join $M(\gamma, p + je_2)$
of $\gamma$ with the copy of $p$ vertically translated by $j$ units. The condition on~$j$ ensures that the hypothesis~(\ref{e.madrascondition}) is met.

\section{Hammersley-Welsh with polygon insertion: reducing~to~a~key~estimate}\label{sec:proofsone}

We begin to implement  the plan laid out in Section~\ref{sec:howpoly} to prove Theorem~\ref{thm:HW1}. 
The present section states but does not prove Proposition~\ref{prop:many_branches}, an upper bound on the number of horizontally confined half-space walks with many branches. This estimate  is at the heart of the proposed plan. The section contains the proof of Theorem~\ref{thm:HW1} that invokes this proposition. In the next Section~\ref{sec:proofstwo}, we give the proof of this key estimate.

We start the rigorous implementation by recalling from the overview the cast of characters $(\eps,\delta,\alpha,u,m)$, the roles of this quintet's members and the conditions that it is necessary to impose on them.
\begin{itemize}
    \item 
    The sought improvement
     in the exponent of the Hammersley-Welsh bound is~$\eps$, so that $n^{1/2 - \eps}$ will replace $n^{1/2}$;
    \item the line-width of inserted polygons will be at least  
    $u$, a quantity of order~$n^\delta$; 
    \item a lower bound on the number of wide polygons, suitable for insertion, is expressed for $\Z^2$ in terms of the exponent
    $\alpha$ in Proposition~\ref{prop:many_poly_subseq};
    \item  and all inserted polygons have a common length, namely~$m$.
\end{itemize} 
Of these five positive parameters, the first three are exponents on which we impose the conditions that
\begin{align}
	&2\eps \leq \delta < 1/2 - \eps \qquad \text{ and } \label{eq:exponent_condition1}\\
	&\eps + \alpha \delta <1/2 \, . \label{eq:exponent_condition2}
\end{align}

The polygon line-width $n^\delta$ is chosen, in the lower bound of~(\ref{eq:exponent_condition1}), to be at least the horizontal near self-touch distance $n^{2\eps}$ between consecutive crossings. There is no use in strict inequality, and so we choose~$\delta = 2\epsilon$, as we did  in the plan in Section~\ref{sec:howpoly}.  
The upper bound of~\eqref{eq:exponent_condition1} is a non-degeneracy condition, which is not restrictive. 
The vital condition~(\ref{eq:exponent_condition2}) coincides with the upper bound on~$\eps$, noted after~(\ref{eq:EMM3}): it ensures that a given polygon insertion is helpful, with the entropic gain in the resulting outgoing arrows being larger than the cost of purchasing the polygon.

By imposing~(\ref{eq:exponent_condition1}) and~(\ref{eq:exponent_condition2}), we are developing rather directly the suggestion prompted by the heuristic discussion of Section~\ref{sec:howpoly}. In fact, such a choice for $(\eps,\delta,\alpha)$ as these conditions impose will permit the derivation of an upper bound on  the number $\vert \HSW_n \vert$ of length~$n$ half-space walks. In order to use this information to obtain Theorem~\ref{thm:HW1}, namely a similar bound on $\vert \SAW_n \vert$, we will invoke~(\ref{eq:SAW_vs_HSW}).
To be able to do so, we will need control on $\vert \HSW_n \vert$ not merely for infinitely many~$n$, but for infinitely many long ranges of consecutive $n$. To achieve such a bound, we will need to strengthen the conditions~(\ref{eq:exponent_condition1}) and~(\ref{eq:exponent_condition2}) on the triple $(\eps,\delta,\alpha)$. 

Reflecting these considerations, the next proposition has two parts. The first part 
articulates the notion that polygon abundance, which limits the  cost of polygon purchase to which we alluded a few moments ago, secures an improvement in the Hammersley-Welsh bound for half-space walks. 
The second part shows that such an improvement is also secured for all self-avoiding walks. As such, it is in this second part that  a strengthening of the conditions~(\ref{eq:exponent_condition1}) and~(\ref{eq:exponent_condition2}) is hypothesised.

We encourage the reader to focus on the elaboration of the heuristic guide and thus on a choice of parameters satisfying~(\ref{eq:exponent_condition1}) and~(\ref{eq:exponent_condition2}). In this regard, note that, when $n \in \N$ satisfies~(\ref{eq:polygonal_existence}), we formally specify the final two members~$(u,m)$ of the parameter quintet given the first three~$(\eps,\delta,\alpha)$. 
\begin{proposition}\label{prop:HW3}
\textnormal{(1)} Fix a triple $(\eps,\delta,\alpha)$ of positive exponents satisfying \eqref{eq:exponent_condition1} and \eqref{eq:exponent_condition2}. 
	There exists a constant $C = C(\eps) > 0$ such that, whenever $n \in \N$ satisfies the condition
	\begin{align}\label{eq:polygonal_existence}
		\text{there exist $u \in [n^{\delta}, 2n^{\delta}]$ and $m \leq 4^{-1} u^2$, for which $|\WSAP_m^{u}| \geq n^{-\alpha \delta} \mu^m$} \, ,
	\end{align}
	we have that
	\begin{align}\label{eq:HW3}
		|\HSW_n| \leq  C \exp ( 21 n^{1/2 - \eps}\log n ) \cdot \mu^n \, .
	\end{align}
	
	\noindent\textnormal{(2)}
	Fix a triplet $(\eps,\delta,\alpha)$ of positive exponents satisfying the stronger conditions 
	\begin{align}\label{eq:exponent_condition_strong}
		4\eps \leq \delta < 1/2 - \eps\quad \text{ and }\quad \eps + 2\alpha \delta < 1/2.
	\end{align} 
	Then there exists a constant $C = C(\eps) > 0$ such that, 
	whenever $n \in \N$ satisfies~(\ref{eq:polygonal_existence}), any $\ell \in  \llbracket n^{3/2},n^2\rrbracket$ satisfies
	\begin{align*}
			|\SAW_\ell| \leq  \exp ( C \ell^{1/2 - \eps}\log  \ell ) \cdot \mu^\ell \,.
	\end{align*}
	
\end{proposition}

Proposition~\ref{prop:HW3}(1) only has value if $\alpha$, $\delta$ and $\eps$ are such that there are infinitely many values $n$ satisfying~\eqref{eq:polygonal_existence}.
Since $\WSAP_m^u \subset \SAP_m$, and it is expected that $\vert \SAP_m \vert = m^{-5/2 + o(1)} \mu^m$ for $m$ even (see~\cite{Nie82}), 
the best upper bound on $\vert \HSW_n \vert$ that the present implementation of our method may achieve is with $\eps$ smaller than, but arbitrarily close to $1/12$. 
Indeed, $1/12$ is the value of $\eps$ in~(\ref{eq:exponent_condition2}) consistent with $\alpha =5/2$. 
That said, aspects of the method other than the $\alpha$-value could be varied: see Section~\ref{s.method}.


Our main result for~$\Z^2$ follows readily from 
Proposition~\ref{prop:HW3}(2)  and the demonstration of polygon abundance in Section~\ref{sec:SAP}.

\noindent{\bf Proof of Theorem~\ref{thm:HW1}.}
	Fix $\eps_0 < \eps  < \frac{1}{466}$. We will prove the theorem for $\eps_0$.  
	With a view to choosing the triple $(\eps,\delta,\alpha)$ to satisfy the hypotheses of Proposition~\ref{prop:HW3}(2), 
	we set $\delta = 4\eps$. 	
A value of $\alpha$ may then be selected so that 
the  hypothesis $\alpha > 29$ of Proposition~\ref{prop:many_poly_subseq} holds, while the bound $\eps(1 + 8\alpha) < 1/2$ of \eqref{eq:exponent_condition_strong} is also met.

 Proposition~\ref{prop:many_poly_subseq} ensures that there are infinitely many values of $(u,m)$ satisfying 
	$\big|\WSAP_m^{u}\big|\geq u^{-\alpha}\mu^m$.
	For any such value of $u$, choose $n$ such that $u/2\leq n^\delta\leq u$. 	
	Then Proposition~\ref{prop:HW3}(2) implies that 
	$|\SAW_\ell| \, \leq \, \exp ( C  \ell^{1/2 - \eps}\log \ell ) \cdot \mu^\ell$ whenever $\ell \in \llbracket n^{3/2},n^2 \rrbracket$.
	This in turn implies that $|\SAW_\ell| \leq\exp ( \ell^{1/2 - \eps_0} ) \cdot \mu^\ell$ for such $\ell$ when $n$ is supposed to be high enough. Taking $j = n^{3/2}$, we obtain the second, stronger, assertion of Theorem~\ref{thm:HW1}.  \qed

We now prepare to give the proof of Proposition~\ref{prop:HW3}, which is to say, to rigorously implement the plan from Section~\ref{sec:howpoly}.
To that end, \begin{center}
	{\bf fix a triple $(\eps,\delta,\alpha)$ of positive exponents satisfying \eqref{eq:exponent_condition1} and \eqref{eq:exponent_condition2}.}
\end{center}
 
The plan began by roughly explaining why the walks in the domain of the Hammersley-Welsh map~$\Psi$ might be chosen to be a certain subset of~$\HSW_n$. This subset was called $\SHSW_n$. As we turn to precise specification, we mention that we will be working with a more general object, $\SHSW_k^n$.

Here, we are fixing parameters $k,n \in \N$ with $k \leq n$, and specify $\SHSW_k^n$ to be a subset of the space $\HSW_k$ of length-$k$ half-space walks. Elements of $\SHSW_k^n$ will adhere nonetheless to the informal description offered at the start of Section~\ref{sec:howpoly} of members of $\SHSW_n$: they are horizontally confined walks with many branches. 

By horizontally confined walks, we mean the elements of a set $\CHSW_k^n$ defined to be the subset of $\HSW_k$ whose elements are contained in the vertical strip 
$[-n^{1/2 + \eps}, n^{1/2 + \eps}] \times \Z$.
 
The set $\FBHSW^n_k$ of length-$k$ half-space walks with few branches  consists of those members of $\HSW_k$
with fewer than  $7 n^{1/2 - \eps}$ elements
in their branch decomposition.
We then define $\SHSW_k^n$ to be $\CHSW_k^n \setminus \FBHSW_k^n$. 

The implementation of our plan comprises  three tasks: the first short, the second substantial, and the third again short.
 The first task
 ascertains that, for the purpose of bounding $\vert \HSW_n \vert$, attention may indeed by focussed on the sets $\SHSW_k^n$; the second provides an upper bound on the cardinality of these sets, with which the proof of  Proposition~\ref{prop:HW3}(1) may easily be presented; and the third proves Proposition~\ref{prop:HW3}(2).

The next two lemmas carry out the first task. 
The following Proposition~\ref{prop:many_branches} states the outcome of the second. After these statements, we prove Proposition~\ref{prop:HW3}(1)  by invoking the three results.
Proofs of the two lemmas then follow. The second task, namely the proof of 
Proposition~\ref{prop:many_branches}, appears in Section~\ref{sec:proofstwo}; naturally, it is there that 
that the central ideas in our plan are enacted. The present Section~\ref{sec:proofsone} ends with the third task, namely
the proof of Proposition~\ref{prop:HW3}(2). 


\begin{lemma}[Restricting to horizontally confined walks]\label{lem:HSW_CHSW} 
For $n \geq 4$, we have that
	\begin{align*}
    	|\HSW_n|\cdot  \mu^{-n}
    	\, \leq \, 3 \mu \, \exp ( 2 n^{1/2 - \eps} \log n ) \Big(\sum_{k = 0}^{n} |\CHSW_{k}^n|\cdot \mu^{-k}\Big)^2.
	\end{align*}
\end{lemma}

\begin{lemma}[Dispatching walks with few branches]\label{lem:few_branches}
For $n$ high enough and	for $k \leq n$, 
	\begin{align}
		\big\vert \FBHSW^n_k \big\vert \leq \exp ( 8 n^{\frac12 - \eps}\log n ) \, \mu^k \, .
	\end{align}
\end{lemma}

\begin{proposition}[The key bound]\label{prop:many_branches}
For some positive constant $C$, for  $n \in \N$ satisfying~\eqref{eq:polygonal_existence} 
and any $k \leq n$,
	\begin{align}
		\big\vert \SHSW_k^n \big\vert \leq C \mu^k \, . \label{eq:many_branches}
	\end{align}
\end{proposition}
\noindent{\bf Proof of Proposition~\ref{prop:HW3}(1).}  
Since $\CHSW_k^n \subseteq \FBHSW_k^n\, \cup \,\SHSW_k^n$,    Lemma~\ref{lem:few_branches} and Proposition~\ref{prop:many_branches} imply that, when $n \in \N$  satisfying~\eqref{eq:polygonal_existence} is high enough, and $k \leq n$, 
    \begin{align*}
    		|\CHSW_{k}^{n}| \leq |\FBHSW^n_k| +|\SHSW_k^n|
    		\leq \exp ( 9 n^{\frac12 - \eps}\log n ) \, \mu^k \, .
    \end{align*}
    Using Lemma~\ref{lem:HSW_CHSW}, we deduce that, for such $n$, 
    \begin{align*}
    	|\HSW_n| \,  \mu^{-n}
     \,	\leq \, \exp( 21 n^{\frac12 - \eps}\log n) \, .
    \end{align*}
Multiplying this right-hand side by  a suitable constant $C(\eps) > 0$ permits us to discard the condition that $n \in \N$ satisfying~(\ref{eq:polygonal_existence}) is supposed to be high enough. We have obtained Proposition~\ref{prop:HW3}(1). \qed

\begin{figure}
	\begin{center}
	\includegraphics[width = 1\textwidth]{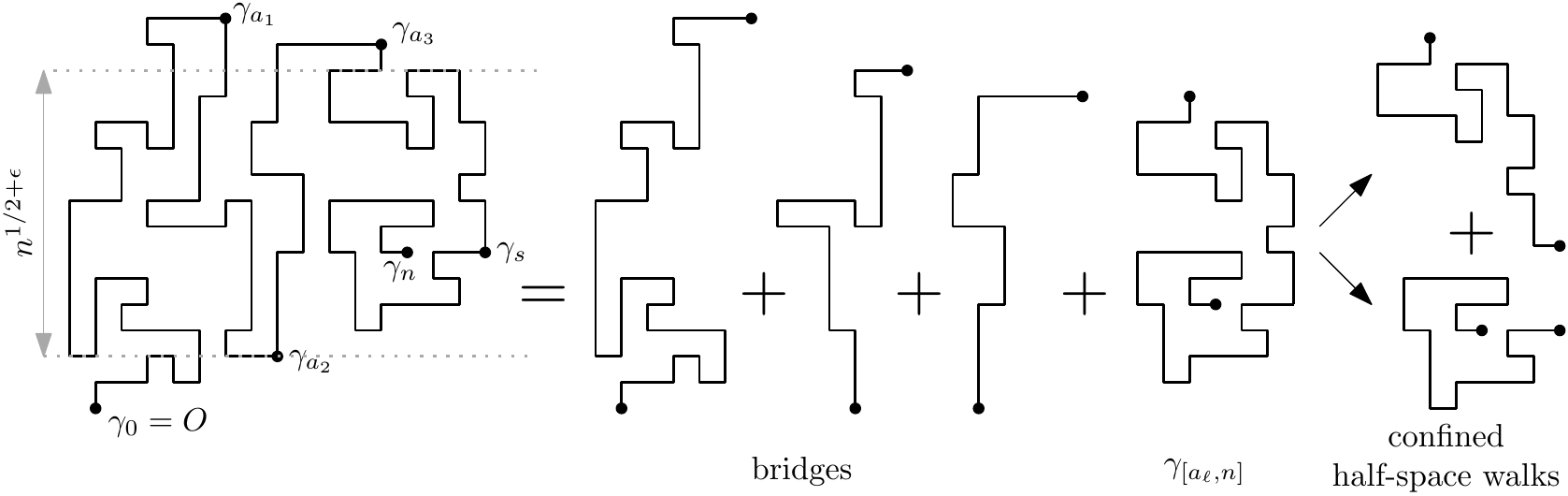}
	\caption{A half-space walk is decomposed into three bridges of height larger than $n^{1/2 + \eps}$ 
	and a final piece of height smaller than  $n^{1/2 + \eps}$. 
	The latter is decomposed into two (rotated) confined half-space walks.}
	\label{fig:CHSW}
	\end{center}
	\end{figure}

\noindent{\bf Proof of Lemma~\ref{lem:HSW_CHSW}.}
	Let $\gamma \in \HSW_n$. Let $\ell$ denote the number of elements in the branch decomposition of~$\gamma$ whose heights exceed~$n^{1/2 + \eps}$, and note that $\ell \leq n^{1/2 - \eps}$.
	We may consider the injective map $\Xi$ defined on $\HSW_n$ and specified by  
    \begin{align*}
  		\Xi: \gamma \, \mapsto  \, \big( \gamma_{[0,a_1]} , \dots , \gamma_{[a_{\ell-1},a_{\ell}]} , \gamma_{[a_{\ell},n]} \big) \, ,
    \end{align*}
    where $\gamma_{[0,a_1]},\dots, \gamma_{[a_{\ell-1},a_{\ell}]}$ are the first $\ell$ branches of $\gamma$. 
	See Figure~\ref{fig:CHSW} for an illustration of the decomposition of~$\gamma$ under~$\Xi$. 
    
 The final component, $\gamma_{[a_{\ell},n]}$, is a half-space walk, contained in a horizontal strip of height $n^{1/2 + \eps}$.
	Let $s$ denote the highest index of a point of maximal $x$-coordinate of $\gamma_{[a_{\ell},n]}$.
After right-angled rotation and translation, $\gamma_{[s, n]}$ becomes an element of $\CHSW_{n - s}^{n}$, while  a suitable such operation on $\gamma_{[a_{\ell}, s]}$ may be followed by the prefixing of a vertically oriented edge  so that  an element of $\CHSW_{s-a_\ell+1}^{n}$ results.
The range of $\Xi$ is thus seen to be a subset of a set in bijection with 
$$
\bigcup_{\ell,a, s} \, \bridgeprod_{a ,\ell} \times \CHSW_{s - a+1}^{n} \times \CHSW_{n-s}^{n} \, ,
$$ 
where the union is over $0 \leq \ell \leq n^{1/2 - \eps}$ and $0 \leq a \leq s \leq n$. By Lemma~\ref{l.bridgeproduct}(1), the displayed set has cardinality at most 
\begin{eqnarray*}
 	& &  \sum_{\ell = 0}^{\lfloor n^{1/2 - \eps} \rfloor} \sum_{s= 0}^n \sum_{a =0}^s \mu^a a^\ell  \cdot  \big\vert \CHSW_{s - a+1}^{n}\big\vert \cdot \big\vert \CHSW_{n-s}^{n} \big\vert \\
 	 & \leq & \big( n^{1/2 - \eps} + 1 \big) n^{n^{1/2 - \eps}}  \mu^{n+1} \\
   	& & \quad \times \sum_{s= 0}^n \sum_{a =0}^s \mu^{-(s-a+1)}    \big\vert \CHSW_{s - a+1}^{n}\big\vert \cdot \mu^{-(n-s)} \big\vert \CHSW_{n-s}^{n} \big\vert \\
    & \leq &  3 n^{2n^{1/2 - \eps}}  \mu^{n+1} \Big( \sum_{k =0}^n \mu^{-k}    \big\vert \CHSW_k^{n}\big\vert \Big)^2
   \, .
\end{eqnarray*}
The final inequality is due to  $\big\vert \CHSW_{k+1}^{n} \big\vert \leq 3 \big\vert  \CHSW_{k}^{n} \big\vert$, a bound seen by the removal of the last edge from each element of $\CHSW_{k+1}^{n}$; to $x + 1 \leq e^x$ with $x = n^{1/2 - \eps}$; and to $n \geq 4$ in the guise $e^x \leq n^x$. Since $\Xi$ is injective, we obtain Lemma~\ref{lem:HSW_CHSW}.
\qed
	
\noindent{\bf  Proof of Lemma~\ref{lem:few_branches}.}
The Hammersley-Welsh map $\Psi$ from Lemma~\ref{l.hsw}
injectively maps $\FBHSW_k^n$ to 
$\bigcup_{j=1}^{\lfloor 7 n^{1/2 - \eps} \rfloor} \bridgeprod_{k,j}$, the set of bridge-lists with at most $7n^{1/2 - \eps}$ terms, which terms have combined length~$k$.
By Lemma~\ref{l.bridgeproduct}(1), the latter set has cardinality at most $\mu^k \sum_{j=1}^{\lfloor 7 n^{1/2 - \eps} \rfloor} k^j \leq \mu^k \cdot 7 n^{1/2 - \eps} \cdot n^{7 n^{1/2 - \eps}}$. By supposing $n \in \N$ to be high enough, we obtain the bound in Lemma~\ref{lem:few_branches}. \qed

Next we prove Proposition~\ref{prop:HW3}(2).
For $\eta,\beta > 0$, $n \in \N$ is said to satisfy Assertion $A(\eta,\beta)$ when 
\begin{align*}
  {A(\eta,\beta)} :\qquad &\text{there exist $u \in [n^{\eta}, 2 n^{\eta}]$ and $m \in \N$ with $m \leq 4^{-1}u^2$}  \\
&\text{for which $|\WSAP_m^{u}| \geq n^{-\beta \eta} \mu^m$.}
\end{align*}
Note that $n \in \N$ satisfies Assertion $A(\delta,\alpha)$ precisely when $n$ satisfies~(\ref{eq:polygonal_existence}).  

\begin{lemma}[Stability of polygon abundance under changes in length]\label{lem:stable}
    Suppose that $n \in \N$ satisfies Assertion~$A(\delta,\alpha)$.
    Then any $\ell \in \llbracket n,n^2\rrbracket$ satisfies Assertion~$A(\delta',2\alpha)$ for some $\delta' = \delta'(\ell) \in [\delta/2,\delta]$.
\end{lemma}

\noindent{\bf Proof.}
Pick $\ell \in \llbracket n,n^2 \rrbracket$.  Assertion~$A(\delta,\alpha)$ for $n$ entails the existence of $u$ and  $m \leq 4^{-1}u^2$ such that $u \in[n^\delta,2 n^\delta]$ and $|\WSAP_m^{u}| \geq n^{-\alpha \delta} \mu^m$.
The first condition can be rewritten as $u\in [\ell^{\delta'},2\ell^{\delta'}]$ with $\delta'=\delta\frac{\log n}{\log \ell}\in[\delta/2,\delta]$. Since $\ell\ge n$ and $2\alpha \delta'\ge\alpha\delta$, the second implies $|\WSAP_m^{u}| \geq \ell^{-\alpha 2\delta'} \mu^m$. Overall, $\ell$ satisfies Assertion $A(\delta',2\alpha)$.
 \qed

\noindent{\bf Proof of Proposition~\ref{prop:HW3}(2).} Assume the fixed values $\varepsilon, \delta,\alpha$ satisfy the stronger condition \eqref{eq:exponent_condition_strong} rather than merely \eqref{eq:exponent_condition1} and \eqref{eq:exponent_condition2}. Fix $n \in \bbN$ that satisfies \eqref{eq:polygonal_existence} and consider $n^{3/2}\le \ell<n^2$; the case where $\ell=n^2$ follows from  $|\SAW_{n^2}|\le 4|\SAW_{n^2-1}|$. 

Applying \eqref{eq:SAW_vs_HSW} and then the Cauchy-Schwarz inequality gives 
\begin{align}\label{eq:hw320}
 &|\SAW_{\ell}|
  \leq \sum_{k=0}^{\ell}	|\HSW_{k+1}| \, |\HSW_{\ell-k}| \le \sum_{k=0}^{\ell+1} |\HSW_{k}|^2\,.
\end{align}
The classical Hammersley-Welsh bound Theorem~\ref{thm:HW} implies that, for all $k\le n$,
\begin{equation}\label{eq:hw322}
	|\HSW_k| \, 
    \,	\leq \, \exp( C  k^{1/2} ) \mu^k ,
\end{equation}
for some constant $C$ independent of $k$ and $n$. Furthermore, by  Lemma~\ref{lem:stable}, any $k \in \llbracket n,n^2 \rrbracket$ satisfies Assertion~$A(\delta',2\alpha)$ for some $\delta' \in [\delta/2,\delta]$.
The strengthened hypothesis \eqref{eq:exponent_condition_strong} on $(\eps,\delta,\alpha)$    
ensures that the triple $(\eps,\delta',2\alpha)$ satisfies the conditions~(\ref{eq:exponent_condition1}) and~(\ref{eq:exponent_condition2}).
Thus, Proposition~\ref{prop:HW3}(1) gives
\begin{equation}\label{eq:hw321}
	|\HSW_k|
    \,	\leq \, \exp( C k^{1/2 - \eps}\log k ) \mu^k ,
\end{equation}
for all $k \in \llbracket n,n^2 \rrbracket$, where $C = C(\eps) > 0$ is some constant (which we can assume to be larger than the constant in \eqref{eq:hw322}). 

Applying \eqref{eq:hw322} for $k\le n$ and \eqref{eq:hw321} for $n < k \leq \ell +1$ to \eqref{eq:hw320}, we obtain 
\begin{align*}
 &|\SAW_{\ell}|\le (\ell+2)\exp(2C\ell^{1/2-\varepsilon}\log \ell)\mu^\ell\end{align*}
since $n^{1/2} \leq  \ell^{1/3} \leq  \ell^{1/2 - \eps}$, where $\eps<1/6$ is imposed by \eqref{eq:exponent_condition_strong}. The desired bound is thus obtained  for  $\ell \in \llbracket n^{3/2} , n^2 - 1 \rrbracket$  by a suitable adjustment of  the constant~$C$. \qed

\section{Proving the key bound on confined walks with~many~branches}\label{sec:proofstwo}

Here we prove Proposition~\ref{prop:many_branches},  implementing the essential aspects of the polygon insertion plan.
For the whole section, the triple $(\eps,\delta,\alpha)$ satisfying \eqref{eq:exponent_condition1} and \eqref{eq:exponent_condition2} is fixed. 

Recall from the conclusion of the heuristically presented plan that we will endeavour to alter a horizontally confined walk $\gamma$ that has  a high number~$r = r(\gamma)$ of branches by the surgical insertion of $\kappa r$
polygons. Here, the parameter $\kappa > 0$ will be of unit order, chosen so that $\kappa r \in \N$.
(This membership of $\N$ will be straightforward to arrange: recall that $r \geq 7n^{1/2 - \eps}$ when $\gamma \in \SHSW_k^n$.) 
The value of  $\kappa > 0$ 
will be fixed when  needed; the consideration that determines its value has been explained at the end of the guide in Section~\ref{sec:howpoly}, though the precise condition will be slightly modified, for a rather inconsequential reason.


\subsection{Locations for polygon insertion}\label{sec:polyinsert}
Let $\gamma \in \SHSW_k^n$. 
In this subsection, we specify a counterpart $\Rect$ of the rectangle~$R$ through which, as we argued in Section~\ref{sec:howpoly}, so many vertical crossings would be made by such walks as~$\gamma$; we specify a set $\joinloc(\gamma)$ of indices for join locations at the side of these vertical crossings at which a polygon may be inserted so as to overlap the adjacent crossing; and we record the set  $\AJLL(\gamma)$ of length-$\kappa r$ lists of such sites. Section~\ref{sec:polyinsert} ends with a lower bound, Lemma~\ref{lem:many_join_loc}, on the prospective forward arrow number, $\big\vert \AJLL(\gamma) \big\vert$, of the multi-valued map~$\Phi$ that we aim to construct.

Let $(n,u,m)$ satisfy~\eqref{eq:polygonal_existence}. 
We denote by 
$\big( \gamma_{[0,a_1]} , \gamma_{[a_1,a_2]} , \dots ,  \gamma_{[a_{r-1},a_r]} \big)$
the branch decomposition of $\gamma$.  
Let $\gamma_{[a_{\tall - 1},a_\tall]}$ be the last branch in the decomposition of $\gamma$ whose height exceeds $n^{1/2 - \eps}$. 
Since the branch heights form a strictly decreasing sequence,  we have $\tall \geq r - n^{1/2 - \eps}$;
and, since $r \geq 7 n^{1/2 - \eps}$, we see that $\tall$ has the same order as $r$.

Write $j_{\max} = y_{\max}(\gamma_{[a_{\tall - 1},a_\tall]})$ and $j_{\min} = y_{\min}(\gamma_{[a_{\tall - 1},a_\tall]})$, so that $j_{\max} - j_{\min} > n^{1/2 - \eps}$.
The branches $\gamma_{[0,a_1]}$,$\,\dots$, $\gamma_{[a_{\tall-1},a_\tall]}$ 
traverse vertically the rectangle 
$$
\Rect := [-n^{1/2 + \eps}, n^{1/2+\eps}] \times [j_{\min},j_{\max}] \, ,
$$ 
which permits us to order them from left to right.
More precisely, a given such branch may traverse this rectangle several times, and these traversals may be ordered from left to right. The leftmost traversals of the $\tall$ branches are themselves ordered from left to right, and it is in this order that we record the branches $\phi^{(1)},\dots, \phi^{(\tall)}$. 
It is these objects that we call the {\em tall branches} of $\gamma$.

For $l \in \llbracket 1, t \rrbracket$, denote by $z_{l,j}$ the rightmost among the vertices of $\phi^{(l)}$ whose $y$-coordinate equals $j$.
Now let $(\ell,j) \in \llbracket 1,\tall-1 \rrbracket \times \bbZ$.
The pair $(\ell,j)$ is called a		
		 {\em viable join index} if
\begin{itemize}
	\item $j_{\min} + 16 u + 3 \leq j \leq j_{\max}  - 16 u - 3$, 
    \item and $x(z_{\ell + 1,j})  - x(z_{\ell,j}) \leq u$.
\end{itemize} 
Let $\joinloc(\gamma)$ be the set of viable join indices of $\gamma$. The next lemma indicates why we call such pairs~$(\ell,j)$ viable: when a wide polygon is inserted on the right of $\phi^{(\ell)}$ at height~$j$, it will overlap the right-adjacent tall branch~$\phi^{(\ell + 1)}$ of~$\gamma$. Recall from Subsection~\ref{sec:notapp} that the subscript in the notation $\MJ_j$ indicates the vertical displacement to which the right polygon is subjected before Madras joining is undertaken.

\begin{lemma}\label{lem:intersections}
Suppose that $n \geq 9^{1/(1-\eps)}$.
    Let $(\ell,j)$ be a viable join location for $\gamma$. 
    Then, for any $p \in  \WSAP_m^{u}$, 
    there exists a vertex of $\bbZ^2$ that belongs to both $\MJ_j \big( \phi^{(\ell)}, p \big)$ and $\phi^{(\ell+1)}$, 
    but not to $\phi^{(\ell)}$.
\end{lemma}
\noindent{\bf Proof.}
By assumption on $(p,j)$, 
	the polygon $P := p + (0,j)$ is contained in the region $\Strip := \R \times [j_{\min} + 3, j_{\max} - 3]$. The pair $\big(\phi^{(\ell)},P\big)$ is thus seen to satisfy~(\ref{e.madrascondition}), so Lemma~\ref{l.madrascondition}(1) permits the construction of  $\MJ_j \big( \phi^{(\ell)}, p \big) = M\big(\phi^{(\ell)},P\big)$. 
	
In the Madras join of $\gamma$ and $P$,  let $(i,0)$ 
denote the vector of horizontal translation which relocates $P$ before surgery takes place. Any vertex of  $P + (i,0)$ at height $j$ lies to the right of any vertex in $\phi^{(\ell)}$ at the same height. Thus, the $x$-coordinate of any such vertex in $P + (i,0)$ exceeds $x(z_{\ell,j})$. 
In particular, $i > x ( z_{\ell,j})$. 
%
%
%
%
The polygon $p$ is wide, and thus $P + (i,0)$ 
attains its line-width of at least~$u$ at height~$j$.
Denoting by $w$ the rightmost point at height~$j$ in this polygon, we thus see that~$w$ lies on the semi-infinite line 
$z_{\ell,j} + [u, \infty) \times \{0\}$. 
Since $z_{\ell + 1,j}$ lies at a distance of at most $u$ directly to the right of $z_{\ell,j}$, the point~$w$ lies directly to the right of~$z_{\ell+1,j}$.
    
    The branch $\phi^{(\ell+1)}$ separates $\Strip$ into several regions, of which two are infinite. 
    Call these two the left and right regions of $\Strip$. 
    By definition, the half line $z_{\ell+1,j} + (0, \infty)\times \{0\}$ is contained in the right region of $\Strip$. The point $w$ lies on the right-side of $P + (i,0)$ in the sense of Lemma~\ref{l.madrascondition}(5); this result thus implies that  $\MJ_j \big( \phi^{(\ell)}, p \big)$ contains~$w$, which is a point in the right region of $\Strip$. 
    However, $\phi^{(\ell)} \cap \Strip$ is  contained in the left region of $\Strip$. Note that  $\MJ_j \big( \phi^{(\ell)}, p \big)$ has non-empty intersection with $\phi^{(\ell)} \cap \Strip$, since, by Lemma~\ref{l.madrascondition}(6), $\MJ_j \big( \phi^{(\ell)}, p \big)$ contains all but at most two edges in $\phi^{(\ell)}$, while $\phi^{(\ell)}$ contains at least $j_{\max} - j_{\min} - 6 \geq n^{1 - \eps} - 6 \geq 3$ edges in view of the assumption that $n \geq 9^{1/(1-\eps)}$.
    We see then that     $\MJ_j \big( \phi^{(\ell)}, p \big)$ intersects $\phi^{(\ell+1)}$ inside $\Strip$. Any vertex in that intersection satisfies the assertion made by Lemma~\ref{lem:intersections}, because such a vertex cannot lie in $\phi^{(\ell)}$ in view of the disjointness of $\phi^{(\ell)}$ and~$\phi^{(\ell + 1)}$. \qed

A {\em join location list} for $\gamma$ is a set $\ajll \subset \joinloc$ of cardinality $\kappa r$ with the property that 
if $(\ell_1, j) \in \ajll$ and $(\ell_2,j')\in \ajll$ are two viable join indices  whose first components satisfy $\vert \ell_1 - \ell_2 \vert \leq 1$, then 
$| j - j' | \geq  32 u + 5$. (The term `list' is used because we prefer to think of inserting polygons one at a time; note, however, that a join location list is an unordered set. Note also that the term `join location list' is a misnomer: the elements of $\ajll$ merely {\em index} physical locations in $\Z^2$.)
Let $\AJLL(\gamma)$ denote the set of join location lists for $\gamma$. 

Set $\SHSW_k^n(r) = \big\{ \gamma \in \SHSW_k^n: \gamma \, \, \textrm{has $r$ branches} \big\}$.

\begin{lemma}\label{lem:many_join_loc}
	If $n$ is large enough, then for every $r$ and every $\gamma \in \SHSW_k^n(r)$,
		\begin{align}\label{eq:many_join_loc}
	|\joinloc(\gamma) | \geq \tfrac14\, r\, n^{1/2- \eps} \, .
	\end{align}
    By further increasing $n$ if need be, we have that 
    \begin{align}\label{eq:many_join_loc2}
    \vert \AJLL(\gamma) \vert \, \geq \, \exp  \, \big( \, \kappa  (1/2-\eps) \, r \log n \, - \, \kappa  \log(8\kappa) \, r \, \big) \, .
    \end{align}
\end{lemma}
\noindent{\bf Proof.}
	We start by proving~\eqref{eq:many_join_loc}.	Recall that $\tall$ denotes the number of $\gamma$'s branches whose height exceeds $n^{1/2-\eps}$. There are  $(\tall-1) (j_{\max} - j_{\min} -32u -5)$ pairs $(\ell,j)$ 
	with $1 \leq \ell  < \tall$ and $j \in \Z \cap \big[ j_{\min} +16 u + 3 , j_{\max} - 16u - 3 \big]$.
If such a pair fails to be a viable join index, it is because the second condition in the definition of such an index is not met. 
	
	Fix such a pair $(\ell,j)$, so that indeed the second condition fails. 
	Let $L$ denote the elements of $\Z^2$ that are encountered strictly and directly to the right of $z_{\ell,j}$, up to and including $z_{\ell,j} + (u,0)$. No vertex of the form $z_{\ell',j}$ lies in~$L$.
	For this reason, we say that $(\ell,j)$ {\em blocks} the $u$ elements of~$L$. 
	Observe that if $(\ell,j), (\ell',j') \in [1,\tall-1] \times \big[ j_{\min} + 16u + 3, j_{\max} - 16u - 3 \big]$ 
	are two distinct pairs that are not viable join indices, 
	then the sets of points that they block are disjoint. 
	Moreover, all blocked points lie in $\big[-n^{1/2 + \eps}, n^{1/2 + \eps} + u \big] \times \big[ j_{\min}+16u + 3 ,j_{\max} - 16u - 3 \big]$. 
	
	It follows that there are at most $u^{-1} \big( 2n^{1/2 + \eps} + u + 1 \big)\big( j_{\max} - j_{\min} -32u - 5 \big)$ 
	 pairs that are not viable join indices. 
	Thus, we may bound from below the number of viable join indices:
	\begin{eqnarray*}
	\big\vert	\joinloc(\gamma) \big\vert 
		& \geq & \big( \tall-1 - 2 u^{-1} n^{1/2 + \eps} - 1 - u^{-1} \big) \big( j_{\max} - j_{\min} - 32u - 5 \big) \\
		& \geq & \big( r - 3n^{1/2 - \eps} - 3 \big) \big( j_{\max} - j_{\min} - 32u - 5 \big) \, , 
	\end{eqnarray*}
	where the latter inequality is due to $u \geq n^{2\eps} \geq 1$ and $r - \tall \leq n^{1/2 - \eps}$. Since $r \geq 7 n^{1/2 - \eps}$, the first term in the latter product is bounded below by $r/2$ when $n$ is high enough. Since $u \leq 2 n^\delta$, while $\delta < 1/2 - \eps$ by~(\ref{eq:exponent_condition1}), the second term in this product is bounded below by $2^{-1} n^{1/2 - \eps}$, provided that $n$ is supposed high enough. 
	This proves~\eqref{eq:many_join_loc}.	
	
	Suppose now that $n$ is indeed large enough for~\eqref{eq:many_join_loc} to hold. 
	In order to count join location lists $\ajll \subset \joinloc(\gamma)$, we record such an element~$\ajll$, one item at a time.  In view of Lemma~\ref{l.madrascondition}(4),
	each choice of a viable join index disallows the prospect of later inclusion in the growing list  for at most $3(32u + 5)$ other elements of $\joinloc(\gamma)$, where the factor of three is contributed because the relevant disjointness condition concerns pairs of indices whose first components differ by at most one from each other. 
	Thus, we find that  
	\begin{eqnarray*}
		 \vert \AJLL(\gamma) \vert 
		& \geq &  \frac{1}{(\kappa r)!} \prod_{i = 0}^{\kappa r - 1}   \big( |\joinloc(\gamma)| - 3(32  u+5)i \big) \\
		& \geq & \frac{1}{(\kappa r)!} \big( 8^{-1} r\, n^{1/2- \eps} \big)^{\kappa r} 
		  \geq  \frac1{2e\sqrt{\kappa r}}  e^{\kappa r} \big( 8^{-1} \kappa^{-1} n^{1/2- \eps} \big)^{\kappa r} \\
		& \geq & (2e\sqrt{\kappa r})^{-1} \exp  \big((1/2- \eps) \kappa r \log n  - Cr \big) \, ,
\end{eqnarray*}
	where $C = \kappa \big( \log (8\kappa) - 1 \big)$. 
	The second inequality is due to~\eqref{eq:many_join_loc}, and, since $u \leq 2 n^\delta$, holds when  
	$n$ is high enough that 
$3(64 n^\delta + 5) \kappa \leq 8^{-1} n^{1/2 - \eps}$. 	
	The third makes use of $h! \leq 2e \, h^{h+1/2} e^{-h}$ with $h = \kappa r$. Replacing $C$ by $\kappa \log(8\kappa)$ and invoking again that $n$ is high enough, we obtain Lemma~\ref{lem:many_join_loc}.
\qed
	
\subsection{Key properties of the multi-valued map $\Phi$}	
	
Henceforth, we let $r \in \N$ denote any {\em given} value satisfying  $r \geq 7n^{1/2 - \eps}$. In order to bound above the number of elements of $\SHSW_k^n$ with $r$ branches, we construct a multi-valued map 
$$
\Phi:    \SHSW_k^n(r) \times \big( \WSAP_{m}^u \big)^{\kappa r} \to \bridgeprod_{k + \kappa r (m + 16),r} \, ,
$$
where the bridge-list set that is the range was specified in Section~\ref{sec:bridge}.

 An element in the domain of $\Phi$ will be recorded in the form $( \gamma \, ; \, p_1,\dots, p_{\kappa r} )$. The arrows under $\Phi$ outgoing from this element will be indexed by the set 
$\AJLL(\gamma)$ of length-$\kappa r$ join location lists.

Next we state a key property that $\Phi$ will be constructed to satisfy. We then use this property to close out the proof of Proposition~\ref{prop:many_branches}. In two further subsections, $\Phi$ is constructed, and the key property is proved. 
\begin{lemma}[The multi-valued map $\Phi$ is close to injective]\label{l.phikeyprop}	
	When $m \geq 17$,
	\begin{itemize}
	\item  	the minimum number $\emm$ of arrows of $\Phi$ outgoing from any element of its domain is at least 
	$\exp \big( \kappa r\big[ (1/2-\eps)\log n - \log(8\kappa)\big] \big)$; 
\item and the maximum number $\EMM$ of arrows incoming to any element in its range is at most $L^{\kappa r}$, where $L=12$.
	\end{itemize} 
\end{lemma}
	
\noindent{\bf Proof  of Proposition~\ref{prop:many_branches}.}
It suffices to prove the statement for $n$ large enough, since smaller values may be incorporated by adjusting the value of the constant $C$. 
Fix $n$ satisfying \eqref{eq:polygonal_existence}.
By the multi-valued map Lemma~\ref{l.mvm},
~\eqref{eq:polygonal_existence} and  Lemma~\ref{l.phikeyprop},
$$
	\big\vert \SHSW_k^n(r) \big\vert  \, \leq \, 
		\Big(\frac{L}{n^{-\delta \alpha } \mu^{m}} \exp \big(- (\tfrac12-\eps)\log n + \log(8\kappa) \big)\Big)^{\kappa r}\,
		\vert \bridgeprod_{k + \kappa r (m + 16),r} \vert \, .
$$
(The hypothesis in Lemma~\ref{l.phikeyprop} that $m \geq 17$ is satisfied because $m$, being the length of a polygon whose line-width is  at least $u$, must satisfy $m \geq u$; while $u \geq n^{2\eps}$; and $n \geq 17^{1/(2\eps)}$, since we are permitted to suppose that $n$ is sufficiently high.)
Using Lemma~\ref{l.bridgeproduct}(1) to bound the right-hand side and the crude bound $k +\kappa r (m+16) \leq n^{1 + 4\eps}$, we find that
\begin{equation}\label{eq:bb}
	\big\vert \SHSW_k^n(r)   \big\vert \, \mu^{-k} \,
	 \leq  \,
	\exp \Big( - r \big[ \kappa  (\tfrac12-  \eps - \delta \alpha) - (1+4\eps) \big] \log n  + r C' \Big) \, ,
\end{equation}
where $C' =   \kappa \log (8\kappa)   + \kappa \log L + 16 \kappa  \log \mu = \kappa \log \big( 96 \kappa \mu^{16} \big)$, since $L=12$ for $\bbZ^2$. 
By~\eqref{eq:exponent_condition2}, the coefficient of $\kappa$ on the right-hand side of~\eqref{eq:bb} is strictly positive. 
Set $\kappa = \lceil \frac{1 + 4\eps}{1/2- \eps - \delta \alpha}+1\rceil$; the rounding up is a device that ensures that $\kappa r \in \N$ in view of $r \in \N$.
Then \eqref{eq:bb} reads 
\begin{align*}
\big\vert \SHSW_k^n(r)   \big\vert \, \mu^{-k} \,
	\leq \, \exp (- r \log n + r C') \leq  n^{-1} \, ,
\end{align*}
the latter by assuming that $n$ is large enough. 
Then
\begin{align*}
	\big\vert \SHSW_k^n \big\vert \, \mu^{-k} \, = \,
	\sum_{r = \lceil 7n^{1/2-\eps} \rceil}^n \big| \SHSW_k^n(r) \big| \, \mu^{-k}
	\leq 1 \, , 
\end{align*}
and Proposition~\ref{prop:many_branches} is proved. \qed

\subsection{Construction of $\Phi$}
An element in the domain of $\Phi$ may be recorded in the form  $( \gamma \, ; \, p_1,\dots, p_{\kappa r} )$,
and an arrow under $\Phi$ outgoing from this element may be recorded as $(\gamma \, ; \, p_1,\dots, p_{\kappa r} \, ;\, \ajll )$, where $\ajll \in \AJLL(\gamma)$.
In a notational device that is intended to draw attention to the proposed surgical locations of the polygons, we  instead denote the respective polygons 
$p_1,\dots, p_{\kappa r}$ in the form 
$\big( p_{(\ell,j)}: (\ell,j) \in \ajll \big)$, where the elements of $\ajll$ are ordered by the value of $\ell$, and by the value of $j$ when the value of $\ell$ is shared between elements.

	Recall from Section~\ref{sec:polyinsert} that  
	$\gamma$'s tall branches $\phi^{(1)},\dots, \phi^{(\tall)}$ are ordered according to the left-right order of the leftmost crossings  of $\Rect$ that they make.

\begin{figure}
    	\begin{center}
 	   	\includegraphics[width = 0.87\textwidth]{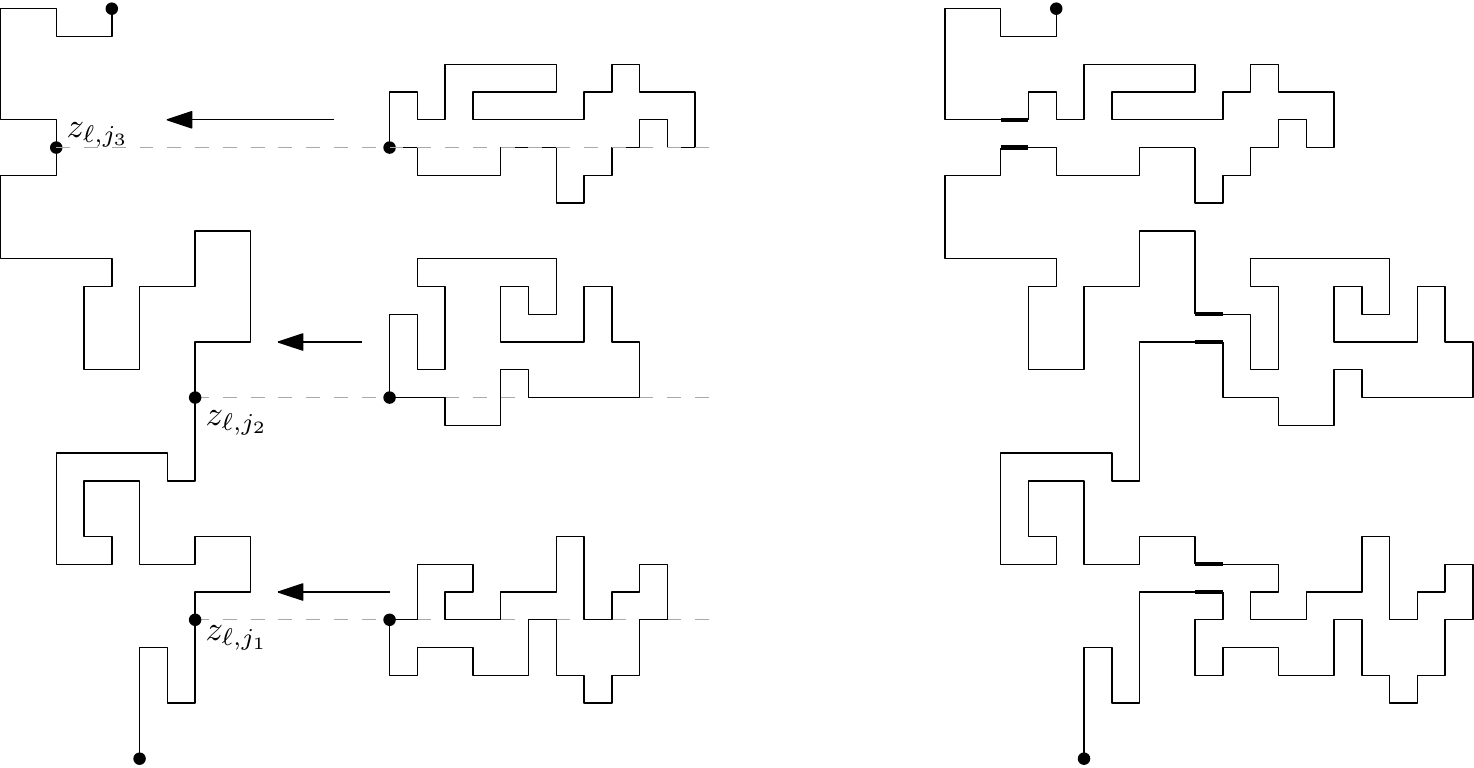}
    		\caption{A toy depiction of the Madras join of a branch $\phi^{(\ell)}$ and three polygons $p_{(\ell,j_1)}$, $p_{(\ell,j_2)}$ and $p_{(\ell,j_3)}$ at respective heights $j_1$, $j_2$ and $j_3$. The edges introduced in surgery have been represented merely by two thick edges in each case.}
    	\label{fig:MJ_multi}
    	\end{center}
	\end{figure}

	For $\ell \in \{1,\dots,\tall-1\}$, define $\tilde \phi^{(\ell)}$ to be the Madras join of $\phi^{(\ell)}$ with all polygons of the form 
	$p_{\ell,j}$ with $(\ell,j) \in \ajll$ at the corresponding heights $j$. 
	Formally, if the pairs $(\ell,j) \in \ajll$ are those with values $j_1,\cdots, j_J$, where $J = J(\ell)$,
	then 
	\begin{align*}
	\tilde \phi^{(\ell)} 
	= \MJ_{j_J}\Big(\cdots \MJ_{j_2}\big( \MJ_{j_1}(\phi^{(\ell)},p_{(\ell,j_1)}), p_{(\ell,j_2)}\big) \cdots \Big) ,p_{(\ell,j_J)}\Big).
	\end{align*}
		Although this specification of the iterated Madras joins is perhaps notationally cumbersome,
the joins are commutative in view of the condition imposed on $\ajll$ by its membership of $\AJLL$.
	Indeed, since distinct $j$-values for given $\ell$ differ by at least $32u + 5$, while the intervals of  $y$-coordinates adopted by the translations used in surgery of the  polygons  $\{ p_{(\ell,j_i)}: 1 \leq i \leq J \}$ each have length at most $16u$, Lemma~\ref{l.madrascondition}(3) ensures this commutativity.
	Moreover, all of these $j$-values are contained in $[j_{\min} + 3, j_{\max} - 3]$, 
so that~(\ref{e.madrascondition}) holds for each proposed join. Hence, Lemma~\ref{l.madrascondition}(1), or its counterpart concerning vertical reflection, implies that  
	$\tilde \phi^{(\ell)}$ inherits from  $\phi^{(\ell)}$ the status of  bridge or vertical reflection of one. 
	The length of $\tilde \phi^{(\ell)}$ is the sum of that of $\phi^{(\ell)}$ and  the quantity $J (m +16)$, 
	since each insertion contributes $m + 16$ edges,
	with $m$ due to the polygon's length and a net gain of sixteen arising from surgery. 
Figure~\ref{fig:MJ_multi} illustrates the several surgical attachments that a given tall branch may endure. 
	
	The arrow under $\Phi$ that we index by $(\gamma \, ; \,  p_1, \dots, p_{\kappa r} \, ; \, \ajll )$ is incoming to a certain range point in $\bridgeprod_{k +  (m + 16)\kappa r,r}$. We now specify this point precisely, so that $\Phi$ is indeed constructed. The half-space walk $\gamma$ has a branch decomposition~(\ref{e.branchdecomp}) in which appear each of the tall branches ~$\phi^{(\ell)}$, $1 \leq \ell \leq \tall$, (though they may not do so in increasing order). 
Replace each of these terms by 	its post-surgical counterpart~$\tilde \phi^{(\ell)}$. Then apply the operation, mapping~(\ref{e.branchdecomp}) to~(\ref{e.bridgedecomp}), that replaces any given branch decomposition by its bridge counterpart. 
The outcome reported by this operation is an element of $\bridgeprod_{k +  (m + 16)\kappa r,r}$. The arrow under $\Phi$ that we consider points to this element. We have constructed~$\Phi$.

\subsection{The near injectivity of $\Phi$: Proof of Lemma~\ref{l.phikeyprop}}
	That $$\emm \geq \exp \big( \kappa r\big[ (1/2-\eps)\log n - \log(8\kappa)\big] \big)$$ for $\Phi$ follows directly from its construction and \eqref{eq:many_join_loc2}. We turn to the upper bound on $\EMM$.

	Let $( b_1,\dots, b_r )$ be an element of $\bridgeprod_{k + \kappa r (m + 16),r}$ to which an arrow under $\Phi$ is incoming. 
	We wish to show that the domain point $( \gamma \, ; \,  p_1,\dots, p_{\kappa r} )$  from which the arrow is outgoing 
	may be determined from $(b_1,\dots, b_r )$ provided that $\kappa r$ choices, each from among at most twelve
	options, are made successively.
	We will establish the stronger assertion in which the data $\ajll$ is also accurately determined by these successful choices. 
	
	
	Let $\chi$ be the concatenation of suitable translations of $b_1,\tau(b_2),b_3,\dots$.
Then $\chi$ is a post-surgical walk,
whose branches are the translations of the just listed terms. Incorporated in these branches are local deformations of the polygons  $p_1,\dots, p_{\kappa r}$ which have been attached by the Madras join operation. 	
Violations of self-avoidance on the part of $\chi$ are the residues of surgery and they will permit a very accurate inference to be made by an observer of $\chi$ as to the form of the $\kappa r$ surgeries by which this walk has been formed.
	
Recall from the start of Section~\ref{sec:polyinsert}	
that the tall branches $\phi^{(1)},\dots,\phi^{(\tall)}$
of $\gamma$ are those among $\gamma$'s branches whose heights exceed $n^{1/2 - \eps}$. It is only onto  these branches that polygons are grafted in the formation of~$\chi$, and any insertion of a polygon onto a branch leaves the height of the branch unchanged. Thus it is that a record of the branches of $\chi$ whose heights exceed $n^{1/2 - \eps}$ is given by the first $\tall$ terms in the sequence
$b_1,\tau(b_2),b_3,\dots$. 

We have recalled  
that the tall branches  $\phi^{(1)},\dots,\phi^{(\tall)}$
of $\gamma$ are recorded in the left-to-right order of their leftmost crossings of $\Rect$. The insertion of polygons into these branches does not produce new crossings, and nor does it affect the intersections of crossings with the horizontal sides of $\Rect$, at heights $j_{\min}$ and $j_{\max}$. The leftmost crossing made by a given branch of either $\gamma$ or $\chi$ occupies the leftmost coordinates at height $j_{\min}$ (and also at height $j_{\max}$)
among all crossings made by this branch. 
For this reason, inspection of the first~$\tall$ terms of the sequence 
$b_1,\tau(b_2),b_3,\dots$ permits a permutation to be recorded in the form 
$\psi^{(1)}, \dots, \psi^{(\tall)}$ so that the terms correspond to the list $\phi^{(1)},\dots,\phi^{(\tall)}$ of the tall branches of~$\gamma$.

The branch decompositions of $\gamma$ and $\chi$
coincide beyond the first $\tall$ terms. This is also true for the $\tall\textsuperscript{th}$ term, because no polygon insertion is attempted in this, the rightmost, case.
Regarding the earlier terms,
let $\ell \in \{ 1,\dots,t-1\}$.
The branch $\psi^{(\ell)}$ is the Madras join $\phi^{(\ell)}$ 
	with a certain number of polygons  which we label $p_{(\ell,j_1)},\dots,p_{(\ell,j_J)}$, where $j_1 < \dots < j_J$ are the height parameters~$j$ in the Madras join operations.

	We aim to retrieve each $\phi^{(\ell)}$ and the polygons $p_{(\ell,j_1)},\dots,p_{(\ell,j_J)}$ by interpreting intersections between $\psi^{(\ell)}$ and $\psi^{(\ell + 1)}$ as signposts that mark where these polygons have been inserted into~$\phi^{(\ell)}$.

First we set $\ell = 1$. Several polygon insertions modify $\phi^{(1)}$ into $\psi^{(1)}$. Any such insertion alters the present walk by removing or adding edges whose endpoints have heights that occupy an integer interval. 
These intervals are disjoint for distinct polygon insertions, and so they may be ordered, from bottom to top. The first polygon insertion that we seek to undo is the lowest, in which $p_{(1,j_1)}$ is inserted into $\phi^{(1)}$.  	
	Lemma~\ref{lem:intersections} indicates that a vertex in~$\phi^{(2)}$ exists that lies in the outcome $\MJ_{j_1}(\phi^{(1)},p_{(1,j_1)})$ of this insertion but does not lie in the input walk~$\phi^{(1)}$.
For definiteness, we select the leftmost of the lowest among such  vertices, and call  it $v$: this vertex  is supposed to offer a clue as to the location of the junction plaquette used in this surgery. 
	Naturally, the observer of $\Phi$'s range point $\big( b_1,\cdots , b_r \big)$ should be able to find~$v$. This observer computes $\psi^{(1)}$ and $\psi^{(2)}$, and identifies $v$ as the lowest point of intersection between these walks, choosing the leftmost among such points if need be. 
	
The vertex $v$ is correctly identified by this means because, 
	in view of  the disjointness conditions placed on join location lists in Section~\ref{sec:polyinsert}, 
	the form of $\psi^{(1)}$ within the horizontal strip whose vertical coordinate interval 
	is that in which perturbations occur in the surgery  $\phi^{(1)} \to \MJ_{j_1}\big( \phi^{(1)},p_{(1,j_1)} \big)$   
	coincides with the form of $\MJ_{j_1}(\phi^{(1)},p_{(1,j_1)})$ within this strip; 
	and likewise the form of $\psi^{(2)}$ within this strip coincides with the form of $\phi^{(2)}$ therein 
	(indeed, the properties of join location lists prevent any polygon from being attached to $\phi^{(2)}$ within this strip). 
	That is, the observer is able to infer the location of $v$ because the essential attributes 
	of the first polygon insertion are not disrupted by later insertions. 
	 	 
	 The next definition and lemma are needed alongside  Lemma~\ref{l.madrascondition}(2) and (3) 
	 to locate the junction plaquette given the value of $v$. 
		

	
\begin{definition}
Let $\ell \in \N$ and $\rho \in \SAW_\ell$.
An index pair $(j,k)$, $0 \leq j < k \leq \ell$,
is called a right-detachable $\rho$-adjacency of gap $k -j$ if 
\begin{itemize}
\item $\big\{ \rho_j,\rho_k \big\}$ is the endpoint pair of a vertically oriented unit edge;
\item the horizontal edges in the plaquette $P$ whose right border is this edge belong to $\rho$ and the vertical edges do not;
\item the modification $\rho \, \Delta \, P$ formed from $\rho$ by the removal of the horizontal edges in $P$ and the insertion of the vertical ones is a disjoint union of a walk $\rho'$ and a polygon $Q$;
\item and any translation of $Q$ directly to the right is disjoint from $\rho'$.
\end{itemize}
\end{definition}

\begin{lemma}\label{l.adjacency}
Let $\ell \in \N$ and let $\rho \in \SAW_\ell$ attain its minimal height at $\rho_0$ and its maximal height at $\rho_\ell$.
 For $h \in \N$, two right-detachable $\rho$-adjacencies $(j_1,k_1)$ and $(j_2,k_2)$ of gap~$h$
satisfy $\vert k_2 - k_1 \vert \geq h$.
\end{lemma}

Lemma~\ref{l.adjacency} depends on planarity, so this tool has a certain importance. However, we defer its slightly irksome proof until the derivation of Lemma~\ref{l.phikeyprop} is finished.

Set $\rho = \MJ_{j_1}\big( \phi^{(1)},p_{(1,j_1)} \big)$.
To find the junction plaquette used in the surgery  $\phi^{(1)} \to \rho$,  let $i \in \N$ satisfy $v = \rho_i$. Since $v$ does not belong to $\phi^{(1)}$,
Lemma~\ref{l.madrascondition}(7) implies that $v$ lies in $\rho$'s surgically implanted middle section, which has length at most $m + 18$.
 The upper-right vertex of the junction plaquette also lies in this section, and thus this vertex's 
 index in $M(\gamma,P)$
 lies in the interval $I := [i - m - 18,i+m + 18]$.
By Lemma~\ref{l.madrascondition}(2) and~(3), the indices in $M(\gamma,P)$ of the  lower-right and upper-right corners of the junction plaquette are the elements of a right-detachable $\rho$-adjacency of gap $m+7$. By Lemma~\ref{l.adjacency}, the  distance between the higher indices in a pair of such adjacencies is at least $m+7$. Since the condition $m \geq 17$ ensures that $2m + 37 < 3(m + 7)$, the interval $I$ may contain the higher index of at most three such adjacencies.

	
	In this way, the junction plaquette in
	$\chi$ associated with the lowest polygon insertion
	among those that construct  $\psi^{(1)}$ may be detected by correctly choosing an element from a set of size at most three. (With further analysis of the Madras join, the junction plaquette may perhaps be detected uniquely, but we will not attempt this detection.) This choice made, the walk and polygon used in this surgery may be determined via Lemma~\ref{l.madrascondition}(8).
	There are at most four forms for the presurgical walk-polygon pair given the surgical outcome alongside the identity of the junction plaquette. 
	
	Thus, the data $\big(  p_{(1,j_1)} , j_1 \big)$ and the form of the deformation to $\phi^{(1)}$ as a result of the lowest insertion onto this branch  may be recovered by correctly choosing an element in a set of size at most twelve.
	
The counterpart recovery should be accomplished for the second lowest polygon insertion for the leftmost tall branch; and then the third; and then all higher  such; and then we may increase the value of $\ell$ so that it equals two; and recover all polygons attached to this branch, from the lowest to the highest; and likewise for all tall branches; and so we remove all inserted polygons. Each of the $\kappa r$ removals introduces a factor of at most twelve to our estimate of incoming arrow number to the given range point of $\Phi$. At the end, we have identified $\big( \phi^{(1)}, \dots , \phi^{(\tall-1)} \big)$, the form of each polygon, and the order of the polygons. Since the higher indexed $\phi^{(i)}$ are already determined, $\gamma$ is recovered as the concatenation $\phi^{(1)}\circ\dots\circ  \phi^{(r)}$.
That is, $(\gamma \, ; \, p_1,\dots, p_{\kappa r} )$ is known;  and so is $\ajll$, though this is incidental.
This completes the proof of Lemma~\ref{l.phikeyprop}. \qed

\begin{figure}
	\begin{center}
	   	\includegraphics[scale = 0.65]{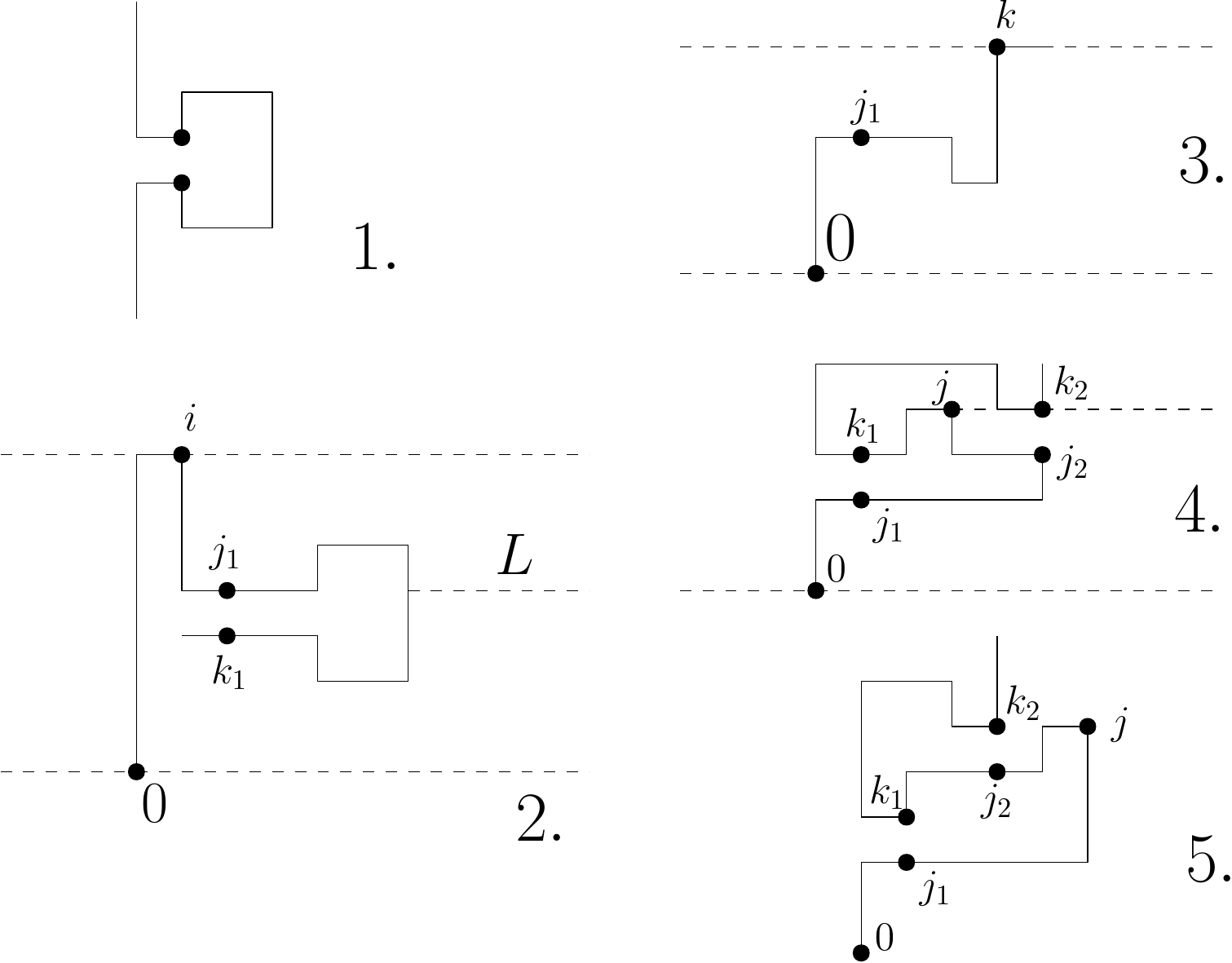}
		\caption{Claims and cases in the proof of Lemma~\ref{l.adjacency}. Each sketch depicts a walk~$\rho$. Each label refers to the chemical distance~$i$ along $\rho$ of the vertex $\rho_i$ indicated by the nearby dot.}
	\label{fig:adjacency}
	\end{center}
\end{figure}

\noindent{\bf Proof of Lemma~\ref{l.adjacency}.}
Let $P$ denote the plaquette whose right border has endpoints $\rho_{j_1}$ and $\rho_{k_1}$, and let $Q$ denote the polygonal component of $\rho \, \Delta \, P$. Note that the vertices of $Q$ are those of $\rho_{[j_1,k_1]}$. 
The rightmost vertex of maximal height in~$Q$, $\northeast(Q)$, will be denoted by $\rho_j$, where note that the index $j$ lies in $\llbracket j_1,k_1 \rrbracket$.  Let $\rho'$ denote the walk component of $\rho \, \Delta \, P$. Figure~\ref{fig:adjacency}(1) depicts the vertex pair,  such as $\{ \rho_{j_1},\rho_{k_1}\}$, corresponding to a right-detachable $\rho$-adjacency as a pair of dots along the walk~$\rho$.

First we show two claims.

\noindent{\bf Claim~1.}  The vertex $\rho_{k_1}$ lies one unit upwards from $\rho_{j_1}$. 

\noindent{\bf Proof.} In the opposing case, $\rho_{k_1}$ is located one unit downwards from $\rho_{j_1}$. Let $\rho_i$ denote a vertex of maximal height on $\rho_{[0,k_1]}$:  see Figure~\ref{fig:adjacency}(2). Let $L$ denote a semi-infinite horizontal line segment whose left endpoint is a vertex  of maximal $x$-coordinate on $\rho_{[j_1,k_1]}$; and note that this vertex lies in $Q$. By our assumption on $\rho_{k_1}$ and planarity, $\rho_{k_1}$ lies in a component of the upper half-plane after the removal of the union of $\rho_{[0,k_1]} \setminus \{ \rho_{k_1} \}$ and $L$ which contains no point of height greater than or equal to that of $\rho_i$. Since $\rho$ terminates at such a height, the path along $\rho$ from $\rho_{k_1}$ must cross~$L$. Since the left endpoint of $L$ lies in $Q$, we see that some right translation of $Q$ intersects $\rho_{[k_1+1,\ell]}$ and thus also intersects $\rho'$. However, this is contrary to the hypothesis that $(j_1,k_1)$ is a right-detachable $\rho$-adjacency. \qed

\noindent{\bf Claim~2.}
The journey along $\rho_{[j_1,k_1]}$ from $\rho_{j_1}$ to $\rho_{k_1}$
first encounters a vertex of maximum height at $\rho_j = \northeast(Q)$. 

\noindent{\bf Proof.} After this journey first reaches  this height, at $\rho_k$ say, a rightward turn on the part of the walk would, as  Figure~\ref{fig:adjacency}(3) depicts,  leave the walk in a different  component of the strip $\R \times [y(\rho_0),y(\rho_k)]$
than the point that is one unit above~$\rho_{j_1}$. By Claim~$1$, the latter location is the endpoint of the journey, $\rho_{k_1}$. 
Thus such a rightward turn is impossible, and Claim~$2$ is proved. \qed

Recall that $(j_2,k_2)$ is also a right-detachable $\rho$-adjacency, and that we seek to show that $k_2 - k_1$ is at least the shared value $k_1 -j_1 = k_2 - j_2$.  
It is thus enough to argue that the intervals $[j_1,k_1]$ and $[j_2,k_2]$ are disjoint, except possibly for some coincidence between endpoints.
To establish this, suppose that the contrary holds; and, without loss of generality, that $j_2 \geq j_1$. Note then that $j_2 < k_1$.

We will analyse three cases: $j_2 < j$; $j_2 =j$; and $j > j_2$. 


\noindent{\em Case~$1$.} Suppose that $j_2 < j$, and see  Figure~\ref{fig:adjacency}(4). Claim~$2$ implies that, in this case, when
the journey along $\rho$ from $\rho_{j_1}$ to $\rho_j$
passes through $\rho_{j_2}$, it does so at a strictly lower height than is achieved at its ending point. 
That is,  $y(\rho_{j_2}) < y(\rho_{j})$.    
Recall that $\rho_{k_2}$ is supposed to be adjacent to $\rho_{j_2}$.  
The path $\rho_{[k_1, k_2]}$ is not permitted to visit heights less than that of $\rho_0$. 
In order to reach its endpoint $\rho_{k_2}$, this path must,  by planarity, 
intersect, at a point other than $\rho_{k_1}$,  the union of $\rho_{[0,j]} \cup \rho_{[j_2,k_1]} = \rho_{[0,k_1]}$
and the horizontal half-line whose left endpoint is $\northeast(Q) = \rho_j$. 
Since self-avoidance prevents such crossing of  $\rho_{[0,k_1]}$, we see that such crossing occurs on the semi-infinite line. 
Thus, some translation of $Q$ to the right intersects $\rho_{[k_1 + 1, k_2]}$. 
But the latter is a subpath of $\rho'$, the walk component of $\rho \, \Delta \, P$, so this inference is contrary to assumption.

\noindent{\em Case~$2$.} Suppose that $j_2 = j$. The vertex $\rho_{j_2 - 1}$ lies directly to the left of $\rho_{j_2}$, because $j_2$ is an element in a right-detachable $\rho$-adjacency. However, $\rho_{j-1}$ lies directly downwards from $\rho_j$, by Claim~$2$. Case~$2$ is thus impossible. 

\noindent{\em Case~$3$.} Suppose instead that $j < j_2 <  k_1$, and see  Figure~\ref{fig:adjacency}(5).  
Write $\overline{P}$ for the plaquette whose right border has endpoints $\rho_{j_2}$ and $\rho_{k_2}$, and write $\bar{\rho}'$ and $\overline{Q}$
for the walk and polygon of which $\rho \, \Delta \, \overline{P}$ is comprised.

By Claim~$2$, the walk $\rho$ turns left after reaching $\rho_j$. It does not reach a higher $y$-coordinate until at least index $k_1$, because $\rho_{[j,k_1]}$ describes part of the polygon~$Q$. 
Thus, $\rho_{[0,j]}$ disconnects every point in $\rho_{[j+1,k_1]}$ from points arbitrarily far to the right within the horizontal strip of boundary heights $y(\rho_0)$ and $y(\rho_j)$.
In Case~$3$, $\rho_{[j_2,k_1]}$ is a subwalk of $\rho_{[j+1,k_1]}$. Since  $\rho_{[j_2,k_1]}$ is a subwalk of $\overline{Q}$, we see that, if $\overline{Q}$ is moved out to the right, it will encounter $\rho_{[0,j]}$ at some point.   The walk $\bar{\rho}'$ contains  $\rho_{[0,j_2 - 1]}$ and thus also  $\rho_{[0,j]}$. Thus, the rightward movement of $\overline{Q}$ will encounter $\bar{\rho}'$. This is contrary to the assumption that $(j_2,k_2)$ is a right-detachable $\rho$-adjacency. \qed

\subsection{Prospects for polygon insertion}\label{s.method}

As we remarked after Proposition~\ref{prop:HW3}, 
and since the parameter $\alpha$ in (\ref{eq:polygonal_existence}) is predicted to be at least $5/2$,
it seems unlikely that the present method may improve Theorem~\ref{thm:HW1} beyond values of $\eps$ 
that are smaller than, but arbitrarily close to $1/12$. 
Indeed, our method of polygon insertion has been presented in an effort to communicate how an improvement of some $\eps > 0$ may be achieved, and we have not been concerned with the explicit value of $\eps$ so obtained. Here we make two comments about how further progress could be made by varying the method's implementation. 

\noindent{\em Obtaining suitably wide polygons more efficiently.}
Elements of $\WSAP_{4m + 4}^u$ were produced in the proof of Lemma~\ref{lem:WSAP_SAB} by joining, under a suitable circumstance, pairs of elements of $\overline\WSAP^u_{2m + 2}$.
Since the latter elements were formed of two bridges in Kesten's proof of Lemma~\ref{l.sapsab}, four bridges are needed to build an element of $\WSAP_{4m + 4}^u$. 
However, we might replace the use of the set  $\WSAP^u_{m}$ in the proof of Theorem~\ref{thm:HW1}
by a set such as $\{ p\in \overline \WSAP^u_{m}:\, y(\northeast(p)) = j\}$, for any given value of $j$. The new polygons, each formed of merely two bridges, may not have sufficient line-width, but they may have significant width when oriented suitably. A potential improvement in Proposition~\ref{prop:HW3}(1) and Theorem~\ref{thm:HW1} may arise from this approach, due to an improvement on the polygon abundance lower bound $\alpha > 29$ in Proposition~\ref{prop:many_poly_subseq}.

\noindent{\em Making inserted polygons of variable length.}
Every polygon that is inserted in the proof of the key estimate, Proposition~\ref{prop:many_branches}, is drawn from $\WSAP_m^u$. As such, each of these polygons has given length $m$. 
One could attempt polygon insertion by permitting this length to vary, say between $m$ and $2m$. There would seem to be a gain in outgoing arrow entropy in the construction of the multi-valued map $\Phi$
which leads to an effective drop  of one in the value of the parameter $\alpha$. In principle, then, the limit of this variant of the method takes $\alpha$ equal to $3/2$, rather than $5/2$. The challenge for this variant is that incoming arrow number may increase because, in the counterpart to the proof of Lemma~\ref{l.phikeyprop}, the value of polygon length on the interval $[m,2m]$ must be surmised as each polygon insertion is undone. 
The post-surgical walk crosses the lower side of the junction plaquette used in surgery, journeys over the inserted polygon, and travels back across this plaquette's upper side. As such the plaquette's location can be roughly surmised, unless the inserted polygon has many locations at which chemically distant vertices are neighbours. The rarity of such locations may be gauged by~\cite[Proposition $4.5$]{Ham17}.

\section{Abundance of polygons for the hexagonal lattice}\label{sec:observable}

Two remaining sections treat self-avoiding walks on the hexagonal lattice.  
In the present section, we prove the polygon abundance estimate for this lattice, Proposition~\ref{prop:many_poly_hex}, which in its exponent and its validity for all lengths improves on Proposition~\ref{prop:many_poly_subseq}. Theorem~\ref{thm:a} will be proved in Section~\ref{s.thm:a}.
Both results depend vitally on a very specific integrability enjoyed by self-avoiding walk on $\bbH$, expressed in terms of the {\em parafermionic observable}. Use of the observable entails some changes in model definition that we present in a first subsection. In a second, we offer a quick introduction to the observable. In a third, we prove  Proposition~\ref{prop:many_poly_hex}. 

\subsection{Definitions for the hexagonal lattice}\label{s.domaindef}

Recall that~$\bbH$ denotes the hexagonal lattice dual to the triangular lattice~$\mathbb T = \bbZ+{\rm e}^{{\rm i}\pi/3}\bbZ$.
The vertices of~$\bbH$ are in the centres of faces of~$\mathbb T$ and the edges of~$\bbH$ are perpendicular to those of~$\mathbb T$. 
In particular, $O = (0,0)$ is the centre of a face of $\bbH$. 
See Figure~\ref{fig:triangle} for an illustration. 

A union~$\Omega$ of faces in~$\mathbb H$ is called a {\em domain} if there exists a self-avoiding polygon~$P$ 
on ~$\mathbb T$ such that~$\Omega$ is the set of faces fully contained in the finite connected component of~$\mathbb R^2\setminus P$. 
Let~$\partial\Omega$ be the intersection of~$P$ with the edges of~$\mathbb H$. 
Note that the elements of~$\partial\Omega$ are the midpoints of those edges of~$\mathbb H$ 
that have exactly one endpoint in the interior of~$P$. 
A triangular domain is depicted in the upcoming Figure~\ref{fig:triangle}. 

An internal edge in $\Omega$ is an edge in $\bbH$
that is contained in $\Omega$. The midpoint of an edge in $\Omega$ refers either to the midpoint of an internal edge or to an element of $\partial \Omega$.
A vertex in $\Omega$ is a vertex in $\bbH$ that lies in $\Omega$.

We now specify notation for self-avoiding walks and polygons on $\bbH$ and on graphs useful for their analysis. The definitions are at variance with our earlier usage.

Let $\Omega$ be a domain. Two midpoints of edges in $\Omega$ are said to be adjacent if the edges share an endpoint.
A self-avoiding walk on $\Omega$
is a sequence of adjacent edge midpoints in $\Omega$
in which there are no repetitions among the list of vertices in $\Omega$ that are the shared endpoints of consecutive members of the sequence. This formulation permits the possibility that the first and final elements in a self-avoiding walk are equal; in this case, the walk is called a self-avoiding polygon. As such, polygons are rooted, and have orientation. When the root and the orientation are forgotten, we recover a counterpart to the definition of polygon used for~$\Z^2$; here, we will call this object a {\em polygon trace}.

The length of a self-avoiding walk is one less than the number of sequence elements; namely, the number of vertices of $\Omega$ visited by the walk. A length zero walk thus refers to one that begins and ends at the same edge midpoint without passing through any domain vertex.

Let $\SAW_n(\Omega)$ denote the set of self-avoiding walks of length $n$ on $\Omega$; the usage is at variance with that of Section~\ref{s.def}, because the starting point of the walks is not fixed. 
When the subscript $n$ is absent, the set of such walks of arbitrary length is designated. 
The notation $\SAP$ refers to polygons in place of walks. Omission of $\Omega$ is understood to imply $\bbH$.


Henceforth $x$ denotes $\mu(\bbH)^{-1}$. Let $z_0$ and $z$ be midpoints of edges of $\Omega$. 
The {\em partition function} of walks from $z_0$ to $z$, also known as the critical two point function, is defined by
\begin{equation}\label{e.gcount}
	G_\Omega(z) = G_\Omega(z_0,z)=\sum_{\substack{\gamma\in \SAW(\Omega)\\ \gamma:z_0\rightarrow z}}x^{|\gamma|} \, , 
\end{equation}
where $\gamma: z_0 \rightarrow z$ indicates that the starting and ending points of $\gamma$ are $z_0$ and $z$, 
and where  $G$'s first argument $z_0$ may be omitted when this value is understood.
When $z=z_0$, a normalized count of polygons is being made.

\subsection{A quick introduction to the parafermionic observable}

That $\mu=\mu(\bbH) = \sqrt{2+\sqrt2}$
is a consequence~\cite{DumSmi12} of analysis of the observable. Applications also include~\cite{BeaBouDum12}.



 
Given $\sigma \in \bbR$, a domain $\Omega$ and $z_0 \in \partial \Omega$,
define the parafermionic observable at any midpoint~$z$ of an edge in $\Omega$ to be
\begin{align}\label{eq:F}
F(z)=F_{\Omega}(z_0,z):=\sum_{\substack{\gamma\in \SAW(\Omega)\\ \gamma:z_0\to z} } e^{-i\sigma{\rm wind}(\gamma)}x^{|\gamma|} \, ,
\end{align}
where ${\rm wind}(\gamma)$ is the total rotation of $\gamma$ from~$z_0$ to~$z$.

It is by now classical (see \cite[Lemma 4]{Smi10a}) that $F$ satisfies the following relations when  
$\sigma$ is set equal to $5/8$ in addition to setting $x= \mu(\bbH)^{-1}$.
For the midpoints $p$, $q$, $r$ of the three edges incident to a vertex $v\in \Omega$,
\begin{equation}
\label{eq_CR}
	(p-v)F(p)+(q-v)F(q)+(r-v)F(r)=0 \, ,
\end{equation}
where $p-v$, $q-v$ or $r-v$ are interpreted as complex numbers.

A simple and important observation is that, for $z\in\partial\Omega$ and $x >0$,  the observable can be related to the generating function $G(z)=G_{\Omega}(z_0,z)$ of walks from $z_0$ to $z$ staying in $\Omega$. 
Indeed, since $z$ lies on the boundary of $\Omega$, the winding of all paths going from $z_0$ to $z$ is the same, so that 
\begin{equation}\label{eq_B}
	F(z)=e^{-i\sigma\, {\rm wind}(z)}G(z) \, ,
\end{equation}
where ${\rm wind}(z)$ is a deterministic constant. 
Here it is crucial that $\Omega$ be simply connected, as is the case here by $\Omega$'s definition. 

We set $\sigma = 5/8$ henceforth.

\subsection{An improved lower bound for $\WSAP$: deriving Proposition~\ref{prop:many_poly_hex}}\label{s.polyhex}

It is proved in~\cite[Remark~$2$]{DumSmi12} that there exists $c > 0$ such that, for $u \geq 1$, 
\begin{align}\label{eq:Bh}
	B_u := \sum_{n\geq 1} \big|\{\gamma \in \SAB_n :\, \height(\gamma) = \tfrac{\sqrt3}2 u \}\big|\cdot\mu^{-n} \, \geq \, c \, u^{-1} \, .
\end{align}
This bound does not directly yield a lower bound on  $b_n \mu^{-n}$, but Lemma~\ref{lem:WSAP_SAB} may perhaps be adapted to express 
a lower bound on $|\WSAP_m^{u}|$ in terms of $B_u$, so that some form of Proposition~\ref{prop:many_poly_hex} would follow from~(\ref{eq:Bh}).
Here, we apply (\ref{eq:Bh}) in a different fashion in order to obtain this proposition.
A similar argument has appeared in \cite{GlaMan17}.

Recall from the introduction that, 
for $k \in \N$,  $\Lambda_{k}$ is the set of faces in $\bbH$ whose distance from the face containing the origin is at most $k$, in the sense of distance on the triangular lattice; and note that $\Lambda_{k}$ is a domain of $\bbH$. 
Let $\FSAP^{(k)}$ be the set of {\em round} self-avoiding polygons: those polygons that are included in $\Lambda_{8k}$; begin and end at $(k+1/2,0)$;  surround, but do not intersect, the planar line segment that interpolates $(0,0)$ and $(k,0)$; and which intersect the positive $x$-axis only at $(k+1/2,0)$.  
There is no demand that the polygons' length be given.

Round polygons are readily available:

\begin{proposition}\label{prop:polygon hexagonal}
    There exists a constant $c>0$ such that, for  $k\in \N$, 
    $$ 
    \sum_{\gamma\in \FSAP^{(k)}}x^{|\gamma|} \, \geq \, c \,  k^{-1} B_{8k}^6 \, .
    $$
\end{proposition}

This proposition will be used in the proof of Theorem~\ref{thm:a}. 
We now use it to prove Proposition~\ref{prop:many_poly_hex} and then give its proof.

\noindent{\bf Proof  of Proposition~\ref{prop:many_poly_hex}.}
This result offers a lower bound on the cardinality of $\WSAP^{u}_m$. We mention that, in our present notation, this is a set of polygon traces, rather than polygons.

	Due to Proposition~\ref{prop:polygon hexagonal},
	\begin{align*}
		\sum_{\gamma\in \FSAP^{(u)}}x^{|\gamma|} 
		\, = \, \sum_{m = u}^{|\Lambda_{8u}|} \big|\{\gamma \in \FSAP^{(u)}:\, |\gamma| = m\}\big|\,\mu^{-m}
		\, \geq \,   c \, u^{-1}B_{8u}^6 \, .
	\end{align*}
	The equality is due to each $\gamma \in \FSAP^{(u)}$ having length between $u$ and $|\Lambda_{8u}|$. 
	The number of terms in the sum is bounded by $c_0 u^2$ for some constant~$c_0$. 
	Thus, there exists at least one value for $m$ such that
	\begin{align*}
		\big|\{\gamma \in \FSAP^{(u)}:\, |\gamma| = m\}\big| \, \mu^{-m} \, \geq \, c_1 \, u^{-3} B_{8u}^6 \, .	
	\end{align*}
	for some universal $c_1 > 0$. 
	Finally, notice that, for any $\gamma$ in the set above, the trace of $\gamma$ is an element of $\WSAP_m^{u}$. 
	Indeed, the line-width of $\gamma$ is at least $u$, since $\gamma$ surrounds the segment $[0,u] \times \{0\}$; and
	its height is bounded by $16u$ because it is contained in $\Lambda_{8u}$.

	\begin{figure}
    	\begin{center}
    	   	\includegraphics[width = 0.35\textwidth]{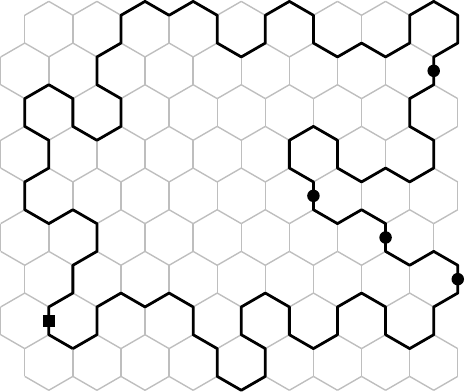}
    		\caption{A polygon trace. 
			Dots mark certain starting points for which the polygon is in $\FSAP^{(u)}$.
The square marks the starting point 
used in the parameterization of elements of $\WSAP^{u}_m$.}
    	\label{fig:rooted}
    	\end{center}
    \end{figure}
 
	The passage $\{\gamma \in \FSAP^{(u)}:\, |\gamma| = m\} \to \WSAP_m^{u}$ from polygon to trace is not injective. 	Any $\gamma \in \FSAP^{(u)}$ is rooted at the rightmost point of its intersection with the horizontal line containing its root: see Figure~\ref{fig:rooted}. 
	Since $\gamma$ visits at most $16u$ different heights, any given trace corresponds to at most $16u$ polygons of $\FSAP^{(u)}$. 
	As such, using \eqref{eq:Bh}, we find that
	\begin{align*}
		\big\vert \WSAP_m^{u}\big\vert \, \mu^{-m} 
		\, \geq \, \frac{c_1}{16u^4}\, B_{8u}(x_c)^6 
		\, \geq \, c_2 \, u^{-10}, 
	\end{align*}
	for some constant $c_2 > 0$. \qed

In the rest of this subsection, Proposition~\ref{prop:polygon hexagonal} is proved.
We start with a lemma. 

\begin{figure}
    \begin{center}
    \includegraphics[scale=0.75]{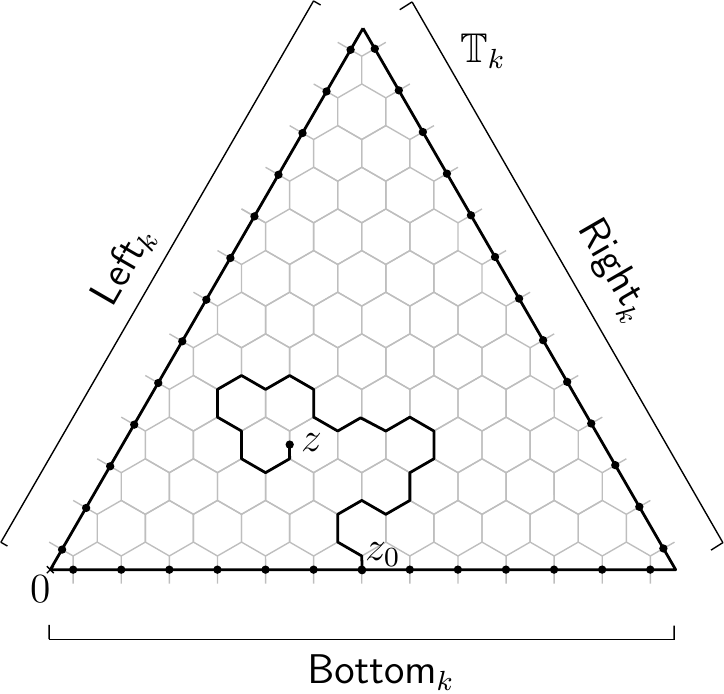}
    \end{center}
    \caption{The domain~$\Omega=\calT_k$: the polygon~$P$ around it is depicted in bold; the edges of~$\calT_k$ are all those with at least one endpoint in~$P$; and the boundary edges are those with exactly one endpoint inside $P$. 
    The centres of the boundary edges, marked with dots, form~$\partial\calT_k$. 
    The three sides of~$\partial\calT_k$ are denoted by~$\mathsf{Bottom}_k$, $\mathsf{Left}_k$, $\mathsf{Right}_k$;
	$z_0$ is the midpoint of the edge in $\calT_k$ emanating from~$a$. 
	The depicted path starts at~$z_0$ and ends at the midpoint~$z$ of an edge inside~$\calT_k$, as in the definition of the parafermionic observable. This path has winding~$2\pi$ at $z$.}
    \label{fig:triangle}
\end{figure}

Fix $k\in \N$. Let the strip $\calS_k$ of height $k$ and the equilateral triangle  $\calT_k$ of side length $2k +1 $  
be defined as the domains whose internal edges are precisely those edges of $\mathbb H$ that are respectively contained in
\begin{align*}
    &\big\{(x,y)\text{ such that }0<y<\tfrac{\sqrt{3}}2 k\big\}\subset\bbR^2\qquad\qquad \text{ and}\\
    &\big\{(x,y)\text{ such that }0< y <\sqrt{3}(k +\tfrac12 -|x - k - \tfrac12|)\big\}\subset\bbR^2 \,;
\end{align*}
see Figure~\ref{fig:triangle}. 
The factor $\sqrt3/2$ above accounts for the height of horizontal layers of $\bbH$, 
that is the vertical difference between the centres of faces of two successive layers.
In an evident notation, let  $\mathsf{Bottom}$ and $\mathsf {Top}_k$ partition $\partial\calS_k$; and let $\mathsf{Bottom}_k$, $\mathsf{Left}_k$ and $\mathsf{Right}_k$ do so for $\partial \calT_k$.
Note that $B_k$ is the partition function of walks contained in $\calS_k$ that start at some fixed point on $\mathsf{Bottom}$ and end on $\mathsf{Top}_k$. 

\begin{lemma}\label{lem:triangle}
    For every $k\in \N$ even,
    $$
    \sum_{\substack{\gamma \in \SAW(\calT_k) \\ \gamma:z_0\to \mathsf{Left}_k } } x^{|\gamma|} \ge x\,B_{2k} \, ,
    $$
    where $z_0= \big(k+1/2,0 \big)$  is in the middle of $\mathsf {Bottom}_k$. 
\end{lemma}

\noindent{\bf Proof.}
Consider the observable $F = F_{\calT_k}$ defined in \eqref{eq:F}. 
Summing~\eqref{eq_CR} over all vertices $v\in \calT_k$, we find that the contributions of each internal edge to the relations around its endpoints cancel each other out; whence
\begin{equation}\label{eq:eq}
e^{\pi i/3}\sum_{z\in\mathsf {Left}_k} F(z)+e^{-\pi i/3}\sum_{z\in \mathsf {Right}_k}F(z)-\sum_{z\in \mathsf {Bottom}_k\setminus\{z_0\}}F(z)=F(z_0) \, .
\end{equation}
Now, $F(z_0)=1$, since only the walk of length zero contributes. Note also that ${\rm wind}(z)$, as defined in~\eqref{eq_B}, is equal to $\pi/3$ on ${\mathsf {Left}}_k$; to $-\pi/3$ on ${\mathsf {Right}}_k$; and to $\pm\pi$ on ${\mathsf {Bottom}}_k$, the choice of sign depending on whether $z$ is on the left or on the right of $z_0$. 
In particular,~\eqref{eq_B} enables us to transform~\eqref{eq:eq} into 
\begin{align}\label{eq:T_k}
	2\cos\big[\tfrac{\pi}{3}(1 - \sigma)\big]\sum_{z\in {\mathsf {Left}}_k} G_{\calT_k}(z)
	+ \cos[\pi(1 - \sigma)\big] \!\!\!\!\!\! \sum_{z\in {\mathsf {Bottom}}_k\setminus\{z_0\}} \!\!\! G_{\calT_k}(z)=1 \, ,
\end{align}
where we used the symmetry of $\calT_k$ with respect to the vertical line running through~$z_0$.

Proceeding for $\calS_{2k}$ as we did for $\calT_{k}$, we obtain
\begin{align}\label{eq:S_k}
	\sum_{z\in {\mathsf{Top}_{2k}}} G_{\calS_{2k}}(z)+\cos[\pi(1 - \sigma)\big]\sum_{z\in {\mathsf {Bottom}}\setminus\{z_0\}} G_{\calS_{2k}}(z)=1.
\end{align}
Here, it may seem problematic that $\calS_{2k}$ is infinite.
However, since $\mu \big( \calS_{2k} \big) < \mu(\bbH) = x^{-1}$, the summation in~(\ref{eq:S_k}) is permitted.
That $\mu \big( \calS_{2k} \big) < \mu(\bbH)$ may be shown using the approach of~\cite[Sec.~6]{HamWhi85},
which implies that $[\mu ( \calS_{2k} )]_k$ is a strictly increasing sequence that converges to $\mu(\bbH)$; 
in~\cite{HamWhi85} the corresponding result is obtained for self-avoiding walks on $\bbZ^d$.  

Finally, observe that, since $\calT_k \subset \calS_{2k}$, for any $z \in {\mathsf{Bottom}}_k$, $G_{\calT_k}(z) \leq G_{\calS_{2k}}(z)$. 
Moreover, the coefficients $\cos[\pi(1 - \sigma)\big]$ and $2\cos\big[\tfrac{\pi}{3}(1 - \sigma)\big] = 1/x$ are both positive. 
Thus, by subtracting~\eqref{eq:S_k} from~\eqref{eq:T_k}, we obtain that, for every $k \in \N$,
\begin{equation}\label{eq:pp}
	\sum_{z\in {\mathsf {Left}}_k} G_{\calT_k}(z) \, \geq \,  \frac{1}{2\cos\big[\tfrac{\pi}{3}(1 - \sigma)\big]}\sum_{z\in \mathsf{Top}_{2k}} G_{\calS_{2k}}(z) \, = \, x\, B_{2k} \, .
\end{equation}
\qed

{\em Remark.} Equation~\eqref{eq:S_k} also applies to strips $\calS_{2k+1}$ of odd height. 
Thus, writing~\eqref{eq:S_k} for $\calS_k$, we find that 
$B_{k} = \sum_{z\in \mathsf{Top}_{k}} G_{\calS_{k}}(z)$ is decreasing in~$k$.

{\em Aside.}
Equation~\eqref{eq:S_k} also shows that
	\begin{align*}
		B_k - B_{k+1} 
		= \cos[\pi(1 - \sigma)\big]\!\!\! \sum_{z\in {\mathsf {Bottom}}\setminus\{z_0\}}\!\!\!\!\!\!\! [G_{\calS_{k+1}}(z) -  G_{\calS_k}(z)]
		\leq \tfrac{\cos[\pi(1 - \sigma)\big]}{x} B_{k+1}^2 \, ,
	\end{align*}
	since any walk contributing to $G_{\calS_{k+1}}(z) -  G_{\calS_k}(z)$ is what may be called an arc from $z_0$ to $z$: a walk that visits $\calS_{k+1} \setminus \calS_k$; and hence such a walk decomposes into two bridges of height $k+1$ (after a small alteration).
	The above was used in \cite{DumSmi12} to prove \eqref{eq:Bh}.

We now present the proof of Proposition~\ref{prop:polygon hexagonal}. A similar strategy was used in \cite{DumGlaPelSpi17} for loop $O(n)$ models. 

\begin{figure}
    \begin{center}
    \includegraphics[scale=1]{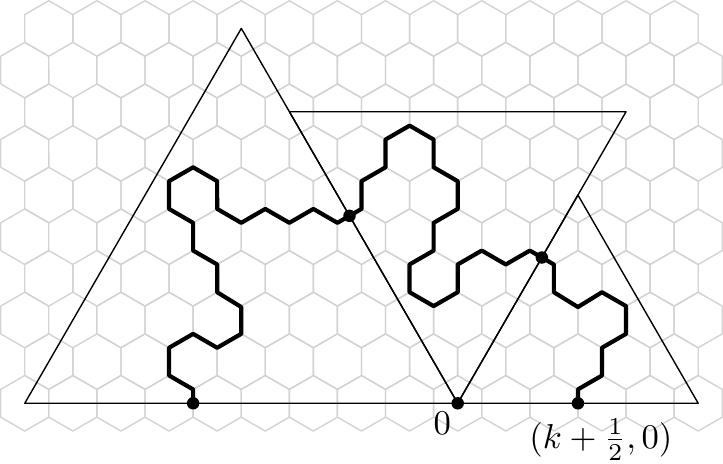}
    \end{center}
    \caption{The concatenation of three walks rotations of which contribute to $G(k,k_1)$, $G(k_1,k_2)$ and $G(k_2,k')$ produces a walk contained in the upper half of $\La_{8k}$, from $(\tfrac12+k,0)$ to $(\tfrac12-k',0)$.}
    \label{fig:triangle_to_poly}
\end{figure}

\noindent{\bf Proof of Proposition~\ref{prop:polygon hexagonal}.}
    For every $k',k'' \geq 0$, define the quantity
    $$
    G(k',k''):=G_{\calT_{k'}}\big(  (\tfrac12+k',0) \, , \, ( \tfrac12+k'' ,0){\rm e}^{{\rm i}\pi/3} \big) \, ,
    $$
    and note that  $G(k',k'') = 0$ if $k'' >2k'$. Let $\La_{8k}^+$ denote the intersection of $\La_{8k}$ and the upper half-plane. 
    Concatenating suitable rotations of three walks that contribute to $G(k,k_1)$, $G(k_1,k_2)$ and $G(k_2,k')$, we obtain 
   a walk in $\La_{8k}^+$ from $(\tfrac12+k,0)$ to $(\tfrac12-k',0)$: 
    see Figure~\ref{fig:triangle_to_poly}.

    Thus,
    \begin{align}\label{eq:arcs}
    &\sum_{k' = 0}^{8k} G_{\La_{8k}^+} \big( z_0, (\tfrac12-k',0)\big) \\
	& \qquad\geq \sum_{k_1 \leq 2k} \sum_{k_2 \leq 2k_1}\sum_{k'\le 2k_2} G(k,k_1)G(k_1,k_2)G(k_2,k')
\, 	,\nonumber
    \end{align}
    whose right-hand side is at least  $x^3\,B_{8k}^3$
    in view of \eqref{eq:pp} and the monotonicity of~$B$ on which we remarked before the proof.  
    
    Consider two walks that contribute to the same summand in the preceding left-hand side. 
    The concatenation of one with the vertical reflection of the other is a polygon of $\FSAP^{(k)}$.
    Thus, and by the Cauchy-Schwarz inequality, 
    \begin{align*}
    	\sum_{\gamma\in \FSAP^{(k)}}x^k
		\, \geq \, \frac{x^6}{8k}\,B_{8k}^6 \, .
	\end{align*}
\qed

\section{The proof of Theorem~\ref{thm:a}}\label{s.thm:a}

The three arenas in which the proof of Theorem~\ref{thm:a}
will unfurl are the lattice $\bbH$; its universal cover $\bbUi$; and its eight-fold cover, which we will call $\ueight$. 
Our description of $\bbUi$ in the introduction was informal. We first precisely specify $\bbUi$ and $\ueight$ and some pertinent subgraphs. 

\subsection{The universal and eight-fold covers of the hexagonal lattice}\label{s.cover}

The positive $x$-axis intersects those edges in $\bbH$ that are vertical, symmetric about this axis and 
whose $x$-coordinates have the form $k + 1/2$ with $k \in \N \cup \{ 0 \}$.
When this semi-infinite line is removed from the plane, each of these edges is cut into two half-edges at its midpoint. The endpoints of the resulting half-edges that lie on the $x$-axis
may be labelled with the terms `lower' and `upper', 
according to whether the argument values of these points, when specified via continuity, are equal to zero or to $2\pi$. The branch-cut hexagonal lattice is the graph so formed: its vertices are the midpoints of edges in $\bbH$, with the exception of those on the positive $x$-axis, for which there are two copies, labelled `upper' and `lower'. 

The universal cover $\bbUi$ of $\bbH$ is formed from a $\Z$-indexed collection of branch-cut hexagonal lattices; any upper vertex on the positive $x$-axis in one lattice is identified with its lower counterpart in the lattice whose index is one greater.

Let $j\in \N$. The $j$-fold cover $\mathbb U^j$ of $\bbH$ 
is formed from a collection of $j$ such lattices. The collection is cyclically ordered, and the cover results when the same  identifications are made using this ordering of the lattice copies. 

The canonical projection from $\bbUi$ to $\bbH$
has been denoted by $\pi_\infty$. The counterpart map from $\ueight$ to $\bbH$ will be called $\pi_8$. 

 At the centre of each face of $\bbH$
is a vertex of the triangular lattice $\mathbb{T}$.
The distance between a pair of these faces is equal to the graph distance in $\mathbb{T}$ between the corresponding vertices. 
Let $k \in \N \cup \{ 0 \}$.
Recall $\La_k$, the set of faces in $\bbH$ whose distance from 
the face containing the origin is at most $k$.
The union of the six planar line segments that interpolate the consecutive rotations of the point $(k+1,0)$ about the origin by the sixth roots of unity is a planar polygon~$P_k$. When $P_k$ is removed from the plane, the lattice $\bbH$ fragments into the bounded domain $\La_k$ and an unbounded piece~$\bbU_k$.
(Reuse of the symbol $\bbU$ is intentional: note that $\pi_\infty^{-1}(\bbU_k)$ equals $\bbUi_k$ from Conjecture~\ref{conj:1}. In fact, the specification of $\bbUi_k$ in the introduction
was rather informal, and we may define it to be  $\pi_\infty^{-1}(\bbU_k)$.)
The two pieces $\La_k$ and $\bbU_k$ share a common set $\partial \La_k$ of boundary edge-midpoints, whose elements lie in the intersection of $P_k$ with the edges in $\bbH$ and are midpoints of such edges.

\subsection{The road map for proving Theorem~\ref{thm:a}}

The expression~(\ref{e.gcount}) is a weighted count of walks.
Upper bounds on such counts are a close surrogate for bounds on the cardinality of the elements of $\SAW_n(\bbH)$ with a given starting point.
Indeed, we begin deriving Theorem~\ref{thm:a} by reducing it to the next result.

\begin{proposition}\label{p.hex}
	There exists $C > 0$ such that, 
	for any $n,k \in \N$ with $n \geq k$ and $\mu(\bbU_k^\infty) = \mu(\bbH)$,
	any $z_0 \in \partial \La_n$ and any midpoint $z$ of an internal edge in $\La_{n-k}$, 
	$$
		G_{\La_n}(z,z_0) \leq  C(k) \, n^{2C(6k+7) + 103} \, ,
	$$
    where $\big\{ C(k): k \in \N \big\}$ is a sequence with $\sup_k k^{-1} \big( \log C(k) \big)^{1/3} < \infty$. 
\end{proposition}

\noindent{\em Remark.} The constant $C > 0$ in the exponent in this result arises from the upcoming Lemma~\ref{lem:a3}. We have not expanded effort on optimizing exponents in our bounds. We mention, however, that $C=32$ is obtained in the proof of Lemma~\ref{lem:a3}, with corresponding explicit values for exponents resulting in Proposition~\ref{p.hex} and in Theorem~\ref{thm:a}, the latter shown in the proof which we now give.	

\noindent{\bf Proof of Theorem~\ref{thm:a}.}
    First, we note that the definition of self-avoiding walk implied in the introduction, in which such walks begin and end at vertices in $\bbH$, is at variance with our actual usage, in which it is the midpoints of edges in $\bbH$ which accommodate these endpoints. It is a trivial matter, however, to interpolate between the two forms of Theorem~\ref{thm:a}.

    Let $\originwalk$ denote the midpoint of the vertical edge in $\bbH$ that borders the right side of the face containing the origin. 
    Since we are assuming Conjecture~\ref{conj:1}, fix $k \geq 1$ for which $\mu(\bbUi_k)=\mu(\bbH)$.

	Let $\gamma$ be a walk in $\bbH$ of length $n$ that starts at $\originwalk$.
   	Write $m = \min \big\{r \in \N: \gamma \subseteq \La_r \big\}$; call $m$ the radius of $\gamma$. 
	Let $\gamma_\ell$ be the first point of $\gamma$ adjacent to $\partial \La_m$, 
	and let $z_0 \in \partial \La_m$ be adjacent to $\gamma_\ell$.
	Then $\gamma_{[0,\ell]}$ may be extended by one step to form a walk that contributes to $G_{\La_m}(\originwalk,z_0)$, while
	$\gamma_{[\ell,n]}$ may be reversed and then so extended to form a walk that contributes to $G_{\La_m}(\gamma_n,z_0)$.
	Suppose for now that the edge midpoint $\gamma_n$ lies in $\La_{m-k}$. 
	We may apply Proposition~\ref{p.hex} with $n =m$ twice to find that
	\begin{equation*}
		\sum_{\substack{\gamma \text{ of radius $m$}\\
		\gamma_{n} \in \La_{m-k} }}x^{|\gamma|}
		\leq \sum_{\substack{z_0 \in \partial \La_m \\z \in \La_{m-k}}} G_{\La_m}(\originwalk,z_0) G_{\La_m}(z_0, z)
		\leq \, c_0 C(k)^2 \, m^{4C(6k+7) + 209}	\,,
	\end{equation*}
	where the first sum is over self-avoiding walks $\gamma$ starting at $\originwalk$ 
	and $c_0$ is a universal constant accounting for  $|\partial\La_{m}|m^{-1}$ and $|\partial\La_{m-k}|m^{-2}$.
	
	It remains to treat the case where $\gamma_n$ lies in $\La_m \setminus \La_{m-k}$. 
	The reversal of $\gamma_{[\ell,n]}$ may be extended from its endpoint $\gamma_\ell$ by $2k$ steps 
	so that its new starting point $z_1$ is on $\partial \La_{m+k}$. Thus,
	\begin{equation*}
		\sum_{\substack{\gamma \text{ of radius $m$}\\
		\gamma_{n} \in\La_m\setminus \La_{m-k} }}x^{|\gamma|}
		\leq \sum_{\substack{z_0\in \partial \La_m\\ z_1\in \partial \La_{m+k}  \\z \in \La_{m}\setminus \La_{m-k} }} 
		G_{\La_m}(\originwalk,z_0) G_{\La_{m+k}}(z_1,z)x^{-2k}.
	\end{equation*}
	As before, by two applications of Proposition~\ref{p.hex} with $n =m$ and $n = m+k$, respectively, 
	the above may be bounded by $c_1(k) C(k)^2 \, m^{4C(6k+7) + 209}$, where $c_1(k)$ is a $k$-dependent constant. 
	
	Restricting the sums in the above displays to walks of length $n$ and summing over all possible values of $m \leq n$, 
	we find that the normalized count of walks of length $n$ that start at $\originwalk$ satisfies 
	\begin{equation*}
		\mu^{-n} \, c_n(\bbH) \, \leq \, C'(k) \, n^{4C(6k+7) + 210}	\, .
	\end{equation*}
\qed

When  $z = \originwalk$ is taken in Proposition~\ref{p.hex},
the proof of the result coincides with the general case's except for the elimination of some mildly distracting notation. There are four principal steps in the proof of the proposition. 

In the next subsection, we state four results, attached to these steps, in the special case when $z = \originwalk$,
and give the proof of Proposition~\ref{p.hex} for $z = \originwalk$.  Four ensuing subsections present the proofs of the four upcoming results.
In a final subsection, the notational changes needed to obtain  Proposition~\ref{p.hex}  in its general form are explained.
Our proof of Theorem~\ref{thm:a} in fact succeeds under a slightly weaker hypothesis than Conjecture~\ref{conj:1}. The present Section~\ref{s.thm:a} ends with a remark in this regard.

\subsection{The four steps to Proposition~\ref{p.hex}}

 Lemma~\ref{lem:a1} is the first step to Proposition~\ref{p.hex} with $z = \originwalk$. 
It is this step which dictates the use of the eight-fold cover. 
  Set $\eightLa_n= \pi_8^{-1}(\Lambda_n)$,  $\eightOmega_n= \pi_8^{-1}(\bbU_n)$, and note that $\partial \eightLa_n =  \pi_8^{-1}(\partial \Lambda_n)$.
 A walk in the eight-fold cover $\ueight$
that begins and ends in 
$\partial \eightLa_n$
is called `inner' if every vertex in $\bbU^8$ that it visits lies in $\eightLa_n$; and `outer' if every such vertex lies in~$\eightOmega_n$. Such definitions may equally be made for the universal cover, with the change $8 \to \infty$.


Walks in $\bbH$ from $\originwalk$
 to $\partial \La_n$, as measured by $G_{\La_n}(\originwalk,z_0)$, are shown to be no more numerous than those walks in the eight-fold cover that begin and end in $\partial \eightLa_n$ and  that are inner.
\begin{lemma}\label{lem:a1}
 Let $n \in \N$. For $z_0\in\partial\Lambda_n$,
 let $\overline{z}_0 \in \partial \eightLambda_n$ satisfy $\pi_8 (\overline{z}_0) = z_0$. 
 Then 
 \begin{equation}\label{e.a1}
    G_{\Lambda_n}(\originwalk,z_0) \, \le \, \frac8{\cos(3\pi/16)}\sum_{z\in \partial\eightLambda_n} G_{\eightLambda_n}\big( \overline{z}_0,z \big) \, .
  \end{equation} 
\end{lemma}

The second principal step provides a lower bound on the number of walks in the universal cover $\bbUi$ between given suitable elements in $\pi_\infty^{-1} \big( \partial \La_n \big)$ that are outer and have  lengths of order roughly $Cn^2$.

 We extend the notation~(\ref{e.gcount}), writing
 $G_{\Omega,\ell}(z_0,z)$ when the sum is indexed by walks whose length is at most $\ell$; we also use it in cases such as $\Omega = \inftyOm_n$ which are not domains, as specified in Section~\ref{s.domaindef}, but for which associated definitions have clear counterparts.
\begin{lemma}\label{lem:a2}
    There exist $c,C>0$ such that, for $n \in \N$,    \begin{align}\label{eq:a3}
    	G_{\inftyOm_n,Cn^2}(z_0,z) \, \ge \, c \, n^{-103} \, ,
    \end{align}
  whenever $z_0,z\in \pi_\infty^{-1}(\partial\Lambda_n)$ are at  distance in $\bbUi$ of at most $48 n$.   
\end{lemma}

Any eight-fold cover inner walk that contributes to the sum on the right-hand side of~(\ref{e.a1}) may be lifted to the universal cover $\bbUi$. An outer walk that shares the endpoints of this lift is  furnished by Lemma~\ref{lem:a2}. 
The concatenation at both endpoints of inner and outer walk is a polygon in the universal cover.
 Thus it is  that, for roughly the price of a walk in $\bbH$ from the origin to distance~$n$, we obtain a polygon in $\bbUi$ of length of order at most $n^2$ that runs through a given vertex in $\partial \inftyLa_n$.
 Our purpose will be served, however, only if we succeed in obtaining a polygon in $\bbUi_k$, rather than one merely in the superset $\bbUi$. To this end, the polygon in~$\bbUi$ will be broken into several polygons so that each avoids $\inftyLa_k$. 
 Our third step thus asserts that, for the given price, we may obtain possibly several such polygons in~$\bbUi_k$.

\begin{corollary}\label{cor:G_ub}
   For $n,k \in \N$ with $n \geq k$, and  $z_0\in\partial\Lambda_n$,
 \begin{align}\label{eq:G_ub}
    	G_{\Lambda_n}\big(\originwalk,z_0\big)
		\, \leq \, C(k) \, n^{103} \Big(\max_{z \in \inftyOm_k} G_{\inftyOm_k,c n^2}(z,z)\Big)^{6k +7} \, ,
    \end{align}
    where $\big\{ C(k): k \in \N \big\}$ is a sequence with $\sup_k k^{-1} \big( \log C(k) \big)^{1/3} < \infty$
    and~$c$ is a positive constant independent of $n$, $k$ and $z_0$.   
\end{corollary}

Our goal is the walk rarity bound  $G_{\La_n}(\originwalk,z_0) \leq n^C$. 
The first three steps have shown that polygons in $\bbUi_k$ of length at most $Cn^2$ are as readily available as are walks contributing to this $G$-count. To reach our goal of a polynomial bound on $G_{\La_n}(\originwalk,z_0)$, our fourth step demonstrates that such polygons are rare in a sense suitable for this purpose.
It is at this moment in the proof of Theorem~\ref{thm:a} that Conjecture~\ref{conj:1} is needed.

\begin{lemma}\label{lem:a3}
  There exist positive constants $C_0$ and $C$  such that, whenever $k \in \N$ satisfies $\mu(\bbUi_k)=\mu(\bbH)$, 
    $n\in \N$ satisfies $n \geq k$, and  $z\in\bbUi_k$,
        $$
        G_{\bbUi_k,n}(z,z)\le C_0 \, n^C \, .
        $$
\end{lemma}

\noindent{\bf Proof of Proposition~\ref{p.hex} with $z = \originwalk$.}
By Corollary~\ref{cor:G_ub} and  Lemma~\ref{lem:a3}, 
$$
	G_{\Lambda_n}\big(\originwalk,z_0\big)  \, \leq \, (c\,C_0)^{6k+7} C(k) \, n^{2C(6k+7) + 103}  \, .
$$
The relabelling of $C(k)$ so that it denotes $(c\,C_0)^{6k+7} C(k)$ does not affect the growth rate of this sequence as the corollary describes it. Thus is the special case of Proposition~\ref{p.hex} obtained. \qed

\subsection{Proof of Lemma~\ref{lem:a1}}
Recall that $n \in \N$; that $z_0 \in \partial \La_n$; that $\overline{z}_0 \in \partial \eightLa_n$ projects to $z_0$ under $\pi_8$; and that $\originwalk$ is the midpoint of a particular edge that borders the face in $\bbH$ containing the origin.  
Eight points in $\eightLa_n$ project to $\originwalk$ under $\pi_8$.
Let $v$ be one among these. The edge in $\bbU^8$ 
of which $v$ is the midpoint will be denoted by~$e$.
A domain to be called $\eightLambda_n(v)$ is formed from $\eightLa_n$ by cutting $e$ in two half-edges. 
These two half-edges may be labelled `plus' and `minus' in an arbitrary fashion.

The domain $\eightLambda_n(v)$ contains, in place of $v$, two vertices $v^+$ and $v^-$, which are endpoints of the new half-edges, and where $v^+$
is the endpoint of the plus half-edge. 
Note then that $\partial \eightLa_n(v) = \partial \eightLa_n \cup\{v^- ,v^+\}$.

In $\eightLa_n(v)$, self-avoiding walks may end at $v^+$, their journey ending via the plus half-edge generated by the cut at $v$; and they may end at $v^-$, if they finish by moving towards this vertex via the minus half-edge; they may not, however, traverse the edge~$e$.

We define the observable in $\eightLa_n(v)$ as we did before, in~(\ref{eq:F}): 
\begin{equation}\label{eq:a1}
	F(z)=F_{\eightLambda(v)}(\overline{z}_0,z) \,
	:= \, \sum_{\substack{\gamma \, \in \, \SAW( \eightLambda_n(v)) \\ \gamma:\overline{z}_0\to z} } e^{-i\sigma\, {\rm wind}(\gamma)}x^{|\gamma|} \, .
\end{equation}
Since $\overline{z}_0 \in \partial \eightLambda_n$
and $v_+$ is the endpoint of a half-edge bordering the singular face of $\bbH$, the difference of value for  ${\rm wind}(\gamma)$
between any two walks $\gamma: \overline{z}_0 \to v_+$ in $\eightLambda_n$ is an integer multiple of $16 \pi$.
Since $\sigma = 5/8$, the $\sigma$-multiple of this difference is divisible by $2\pi$, 
and the phase factor in~(\ref{eq:a1}) of all walks contributing to $F(v_+)$ is a constant $e^{-i \sigma {\rm wind}(v_+)}$. 
The same applies to walks contributing to $F(v_-)$, for which the phase factor is $e^{-i \sigma ({\rm wind}(v_+)-\pi)}$. 
(It is this consideration that determines our use of the eight-fold cover.) 

The local relation~\eqref{eq_CR} is obtained in the same way as it is in the planar case.
Using the equivalent for $\eightLambda_n(v)$ of  \eqref{eq:eq}, we deduce that 
\begin{equation}\label{eq:a2}
	\big\vert F(v_+)-F(v_-) \big\vert 
	\, \leq \!\!\!\!\!\! \sum_{z\in\partial\eightLambda_n(v) \setminus \{v_+,v_-\}} \!\!\!\!\!\!  
		G_{\eightLambda_n(v)}\big( \overline{z}_0,z \big)
	\, \leq \, \sum_{z\in\partial \eightLambda_n}G_{\eightLambda_n}\big( \overline{z}_0,z \big) \, .
\end{equation}
Now, keeping in mind the phase difference $e^{i \sigma  \pi} = -e^{-i3\pi/8}$ between $F(v_+)$ and $F(v_-)$, 
after multiplication by $e^{i \sigma{\rm wind}(v_+) + i3\pi/16}$,  we find that
$$
    \big\vert F(v_+)-F(v_-) \big\vert 
    \, = \, \Big\vert \,  G_{\eightLambda_n(v)}(\overline{z}_0, v_+) e^{i3\pi/16} + G_{\eightLambda_n(v)}(\overline{z}_0, v_-) e^{-i3\pi/16} \, \Big\vert \, . 
    $$
For $a \in \mathbb{C}$, $\vert a \vert \geq  (a + \bar{a})/2$; thus is the last right-hand side at least    
    $$
  \cos(\tfrac{3}{16}\pi)\, \Big( G_{\eightLambda_n(v)}(\overline{z}_0, v_+) + G_{\eightLambda_n(v)}(\overline{z}_0, v_-) \Big)
    = \cos(\tfrac{3}{16}\pi) \, G_{\eightLambda_n}(\overline{z}_0, v) \, .
$$ 
Any self-avoiding walk in $\La_n$ is the projection of at least one in $\eightLambda_n$; whence 
$$
	G_{\Lambda_n}( o , z_0 )
	\leq \sum_{v :\, \pi_8(v) = o}G_{\eightLambda_n}(\overline{z}_0, v) \, .
	$$
Since there are eight such $v$,	
	$$
		G_{\Lambda_n}( o , z_0 )
		\leq\tfrac8{\cos(3\pi/16)}\sum_{z\in\partial\eightLambda_{n}}G_{\eightLambda_{n}}(\overline{z}_0,z),
	$$ 
	by~(\ref{eq:a2}) and the bound that is derived after it. This completes the proof of Lemma~\ref{lem:a1}.
\qed

\subsection{Deriving Lemma~\ref{lem:a2}}
\begin{figure}
	\begin{center}
	   	\includegraphics[width=0.47\textwidth]{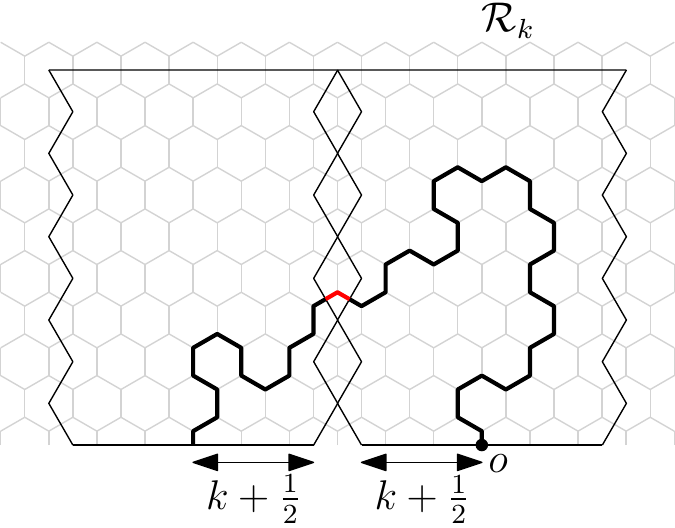}\hspace{0.04\textwidth}
		\includegraphics[width=0.47\textwidth]{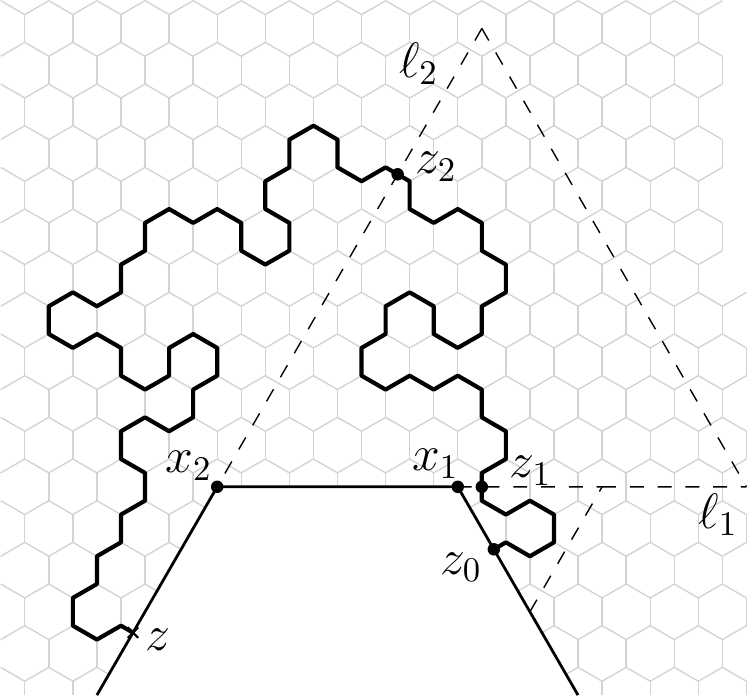}
		\caption{{\em Left:} 
		The domain  $\calR_k$ and a translate of it. Two walks contributing to $G_{\calR_k}(o,z)$ for some $z \in \mathsf{Left}'_k$ 
		(one is reflected and translated) may be used to create an arc of between $o$ and $(-2k-\tfrac12,0)$. The red edges are added in the concatenation 
		{\em Right:} To connect $z_0$ to $z$ we use \eqref{eq:b1} to arrive to some point $z_i$ on the same line as $z$, 
		then use \eqref{eq:b2} to connect $z_i$ to $z$.}
	\label{fig:rect}
	\label{fig:going_around}
	\end{center}
\end{figure}

Two estimates for partition functions in the half-plane hexagonal lattice 
$\mathbb H^+ = (\bbR \times [0,+\infty) ) \cap \bbH$ will be needed in this derivation. 
Recall the triangle~$\calT_{k}$ of 
Figure~\ref{fig:triangle}; 
translate it so that $\originwalk$ is the centre of $\mathsf{Bottom}_k$. 
Due to \eqref{eq:pp} and \eqref{eq:Bh}, there exists a constant $c >0$ such that, for $k \in \N$,
\begin{align}
	\sum_{z\in\mathsf {Left}_k}G_{\calT_k}(\originwalk,z)&\ge  c \, k^{-1} \, . \label{eq:b1}
\end{align}
Our second estimate refers to walks from $\originwalk$ to arbitrary edge midpoints on the positive $x$-axis. 
It states that there exists $c > 0$ such that, for  $k \in \N$, 
\begin{align}
	G_{\La_{4k} \cap \bbH^+}\big( \originwalk,( k + 1/2,0) \big)\geq c \, k^{-7} \, .\label{eq:b2}
\end{align}
We now prove~(\ref{eq:b2}). 
Consider the rectangle $\calR_{k}$ illustrated in Figure~\ref{fig:rect}.  
Its base, $\mathsf{Base}_k$, is centred at $\originwalk$ and has width $2k+1$; the height of $\calR_k$ is $4\sqrt3 \, k$. 
Write $\mathsf{Left}'_k$, $\mathsf{Top}'_{8k}$ and $\mathsf{Right}'_k$ for its left, top and right sides. 
Writing \eqref{eq:T_k} for $\calR_k$ and subtracting \eqref{eq:pp} for $\calS_{8k}$, we find that
\begin{align*}
    \sum_{z\in \mathsf{Left}'_k}\!\!\! 2\cos\big(\tfrac{(3\pm1)\pi}{16}\big)\, G_{\calR_k}(\originwalk,z) 
    \, \geq \, \cos(\tfrac{3\pi}8)\hspace{-0.7cm}\sum_{z \in \mathsf{Bottom} \setminus \mathsf{Base}_k} \hspace{-0.7cm} G_{\calS_{8k}}(\originwalk,z) \,\, 
    +\hspace{-0.5cm}\sum_{z \in \mathsf{Top}_{8k}  \setminus \mathsf{Top}'_{8k} }\hspace{-0.5cm}G_{\calS_{8k}}(\originwalk,z) .
\end{align*}
The inequality comes from the terms omitted from the right-hand side, 
namely those corresponding to arcs ending on $\mathsf{Base}_k$ and bridges ending on $\mathsf{Top}'_{8k}$,
which are contained in $\calS_{8k}$ but not in~$\calR_k$.
These have prefactors $\cos(\tfrac{3\pi}8)$ and $1$, respectively, hence would contribute positively to the right-hand side. 
The angle in the sum on the left-hand side depends on the orientation of $z$; 
it alternates between $\pi/8 = (1-\sigma) \frac{\pi}3$ and $\pi/4 = (1-\sigma) \frac{2\pi}3$. 
Both values lead to prefactors that are strictly positive.
Thus, there exist constants $c_1,c_2 > 0$ such that 
\begin{align*}
	\sum_{z\in \mathsf{Left}'_k} G_{\calR_k}(\originwalk,z) 
	\, \geq 
   	\, c_1 \hspace{-0.6cm} \sum_{z \in \mathsf{Bottom} \setminus \mathsf{Base}_k} \hspace{-0.6cm}G_{\calS_{8k}}(\originwalk,z) 
	\, \geq 2 c_1 x^3 B_{8k}^3
	\, \geq \,  c_2 \, k^{-3} \, .
\end{align*}
The second inequality is due to \eqref{eq:arcs} and the phrase following it;
the last is implied by the lower bound~\eqref{eq:Bh} on $B_{8k}$.

Next, observe that, for any $z \in \mathsf{Left}'_k$, two walks contributing to $G_{\calR_k}(\originwalk,z)$ may be combined, after one of them is reflected horizontally and translated, to form a walk from $\originwalk$ to $(-2k-3/2,0)$; see again Figure~\ref{fig:rect}.
Moreover, this walk is contained in the upper half-plane $\bbH^+$ and in $\La_{8k}$.
Thus, there exists $c_3 > 0$ such that 
\begin{align*}
	G_{\La_{8k} \cap \bbH^+}\big(\originwalk,(-2k-\tfrac32,0)\big) \, 
	\geq \,
	\sum_{z\in \mathsf{Left}'_k} x\,G_{\calR_k}(\originwalk,z)^2
	\, \geq \,
   	c_3 \, k^{-7} \, .
\end{align*}
(The factor $x$ is due to the red piece in Figure~\ref{fig:rect}.) 
A similar construction, with a slightly altered shape for $\calR_k$, 
yields the same lower bound for $G_{\La_{8k} \cap \bbH^+}\big(o,(-2k-\tfrac12,0)\big)$.
By relabelling $2k \to k$, monotonicity of $G$ in the domain, and symmetry, \eqref{eq:b2} is proved. 

\noindent{\bf Proof of Lemma~\ref{lem:a2}.}
Recall that $n \in \N$ and that the distance in $\bbUi$ between $z_0,z \in\pi_\infty^{-1}(\partial\Lambda_n)$ is at most $48n$. 
The elements of $\pi_\infty^{-1} ( \partial\Lambda_n )$ may be labelled by $\Z$, with increasing index corresponding to counterclockwise movement along the copy of $\partial \Lambda_n$
obtained by projection via $\pi_\infty$. The index set $\Z$ is partitioned into countably many intervals, each of length~$n$, according to the segment in the planar polygon~$P_n$ on which the $\pi_\infty$-image of the indexed element lies. Indeed, an unending counterclockwise journey along this polygon may be lifted via $\pi_\infty$ so that the elements of $\pi_\infty^{-1}(\partial\Lambda_n)$ are visited in order; in fact, in order for this to make sense, $\pi_\infty$ should be extended, compatibly with the present definition, to be a projection from the universal cover of $\mathbb{C}$ without the origin to this base space. Any point on this journey whose projection under $\pi_\infty$ is the endpoint of one of the segments comprising the polygon $P_n$ is called a corner.
Corners naturally carry half-integral indices, and each partition interval is bookended by the indices of a pair of corners.

It is without loss of generality that we may suppose that  $z$ is higher than $z_0$ in this ordering of $\pi_\infty^{-1}(\partial\Lambda_n)$.
Let $x_1,x_2,\dots$ denote the corners 
encountered on the journey from $z_0$ to $z$ along $\pi_\infty^{-1}(\partial\Lambda_n)$ dictated by the ordering. 
This sequence may be empty, and may have at most forty-eight terms. 
Write $\ell_i$ for the half-line that emanates from $\pi_\infty (x_i)$ in the direction opposite to $\pi_\infty (x_{i+1}) - \pi_\infty(x_i)$.
See Figure~\ref{fig:going_around} for an illustration.

Let $c > 0$ be such that \eqref{eq:b1} and \eqref{eq:b2} are valid. 
If there is no corner between $z$ and $z_0$, then \eqref{eq:b2} shows that 
\begin{align*}
	G_{\inftyOm_n,3Cn^2}(z_0,z) \, \geq \,  c \, n^{-7} \, ;
\end{align*}
since the number of edge midpoints in $\La_{4n} \cap \bbH$ is at most $3 \cdot 3/2 \cdot (4n)^2 =  72 n^2$, this bound holds provided that $C \geq 24$. 

Suppose now that the sequence  $x_1,x_2,\dots$  has  $i \in \N$ terms. Thus, $\pi_\infty(z)$ lies on the planar line segment that contains  
$\pi_\infty(x_i)$ and $\pi_\infty(x_{i+1})$. Applying \eqref{eq:b1} with $k = n$ after a suitable translation, and using $\vert \mathsf {Left}_n \vert = n$, we find that  there exists a point $z_1$ on $\ell_1$ such that $G_{\calT(z_0)}(z_0,z_1) \geq  c \, n^{-2}$, 
where $\calT(z_0)$ is the equilateral triangle the centre of whose base is $z_0$  and whose left corner is $x_1$. 
(More accurately, the image under $\pi_\infty$ of  $\calT(z_0)$ is an equilateral triangle with corresponding data; we abuse notation in a moment by disregarding a similar nicety.)
The distance between $\pi_\infty(x_2)$ and $\pi_\infty(z_1)$ is at most $3n$. 
Applying \eqref{eq:b1} again, there exists a point $z_2$ on $\ell_2$ such that $G_{\calT(z_1)}(z_1,z_2) \ge c \, (3n)^{-2}$,
where $\calT(z_1)$ is the equilateral triangle whose base is centred at $z_1$ and whose left corner is $x_2$. 
This process may be repeated until a point $z_i$ on $\ell_i$ is reached. 
There are at most forty-eight steps and the size of each triangle is at most three times the preceding one's. 
When walks that respectively contribute to the terms $G_{\calT(z_0)}(z_0,z_1)$, $G_{\calT(z_1)}(z_1,z_2),\dots$ are concatenated, an outer walk between $z_0$ and $z_i$ is obtained. Thus, we find that
\begin{align*}
	G_{\bbUi_n,C n^2}(z_0,z_i) \, \ge \, c'  \,  n^{-96} \, ,
\end{align*}
where $C = 3^{2 \cdot 47} \cdot 9/8$ and $c'  = c \cdot 3^{-2(1 + 2 + \dots + 47)} \geq c \, 3^{-(47)^2}$.
Finally, we apply \eqref{eq:b2} to connect $z_i$ to $z$. By relabelling $C$ and $c$, we obtain Lemma~\ref{lem:a2}. \qed

\subsection{Proof of Corollary~\ref{cor:G_ub}}
Recall that $n \in \N$ and  $z_0 \in \partial\Lambda_n$. In view of Lemma~\ref{lem:a1}, we wish to bound above
$\sum_{z\in \partial\eightLambda_n} G_{\eightLambda_n}\big( \overline{z}_0,z \big)$, where $\overline{z}_0 \in \partial \eightLa_n$ satisfies $\pi_8(\overline{z}_0) = z_0$.

Let $\gamma$ be a walk contributing to this sum; $\gamma$ is an inner walk in $\eightLa_n$, moving from $\partial \eightLa_n$ to itself.
Suppose that $\gamma$ visits $\eightLa_k$. We may decompose $\gamma$'s lifetime between  its first and last visits to $\eightLa_{k+1}$ into outer excursions $\gamma_s$
from $\partial \eightLa_{k+1}$ to itself; and intervening journeys $\psi_s$ within $\eightLa_{k+1}$.
That is, for some $S = S(\gamma) \leq \vert \partial \eightLa_{k+1} \vert$, we may write 
$\gamma = \phi_{{\rm start}} \circ \psi_1 \circ \gamma_1 \circ \psi_1 \circ \dots \circ \gamma_{S-1} \circ \psi_S \circ \phi_{{\rm end}}$, 
where
$\phi_{{\rm start}}: \overline{z}_0 \to \partial \eightLa_{k+1}$ and 
$\phi_{{\rm end}}: \partial \eightLa_{k+1} \to z$ lie in $\eightLa_n \setminus \eightLa_{k+1}$;
for $1 \leq i \leq S$, $\psi_i: \partial \eightLa_{k+1}\to \partial \eightLa_{k+1}$ lies in~$\eightLa_{k+1}$; and, 
for $1 \leq i \leq S-1$, $\gamma_i: \partial \eightLa_{k+1}\to \partial \eightLa_{k+1}$ lies in~$\eightLa_n \setminus \eightLa_{k+1}$. In this description, and later in this proof, `$\phi$ lies in $\Omega$' means `excepting $\phi$'s endpoints, every $\phi_i$ is the midpoint of an internal edge in $\Omega$'.

The endpoints of $\phi_{{\rm start}}$ and $\phi_{{\rm end}}$ on $\partial\Lambda^8_{k+1}$ 
may be joined by  a walk $\gamma'$ that lies in $\eightLa_{k+1} \setminus \eightLambda_{k}$. 
For $1 \leq i \leq S-1$, the endpoints of $\gamma_i$ may be joined by  a walk that lies in 
$\eightLa_{k+1} \setminus \eightLa_k$ so that a polygon $p_i$ 
in $\eightOm_k$ results. The ordering of the endpoints of $\gamma_i$ used to produce the walk may be selected so that $p_i$ does not disconnect the origin from infinity in $\ueight$.  Intersections between the walk $\gamma'$ and the polygons so produced are possible and do not concern us.

The topological property of the produced polygons permits that they be lifted to the universal cover, where they are contained in $\inftyOm_k$; we redesignate $p_i$ to denote the polygon so lifted.
 
When $\gamma$ does not visit $\eightLambda_k$, set 
$\gamma' = \gamma$ and $S = 1$.
Any point in the range of the application $\gamma \to \big( \gamma', p_1, \dots, p_{S-1} \big)$ 
has at most $C_0$ preimages, where $C_0$ is a certain constant that depends on $k$. 
Indeed, $\gamma$ may be inferred from the image point by the reconstruction of the walks $\psi_i$, 
of which there are at most $\vert \partial \eightLa_{k+1} \vert \leq 6 (k+1)$.
Since each $\psi_i$ lies in $\eightLa_{k+1}$, we may choose 
$C_0 = 2^{432(k+1)^3} \geq 2^{\vert \eightLa_{k+1} \vert \cdot 6(k+1)}$; 
here, and later in this proof, $\vert \eightLa_\ell \vert$ abusively denotes the number of internal edges in $\eightLa_\ell$. 
Each polygon $p_i$ has length at most $\big\vert \eightLambda_n \big\vert \leq 72 n^2$.
Moreover, there are at most $12(k+1)$ edges added to each $\gamma_i$ to form $p_i$. 
Finally, the number of polygons satisfies $S-1 \leq 6(k+1)$. 
Thus, by setting $C_1  = C_0\, x^{-72(k+1)^2}$, we find that
\begin{equation}\label{eq:op}
	G_{\eightLa_n} \big(\overline{z}_0,z \big)
	\leq  C_1 \, \Big(\max_{z_1 \in  \partial \inftyLa_{k+1} } G_{\inftyOm_k,72 n^2}(z_1,z_1)\Big)^{6(k+1)} 
	\, G_{\eightLambda_n \setminus\eightLambda_k} \big( \overline{z}_0 , z \big) \, .
\end{equation}

Next, we bound $G_{\eightLambda_n\setminus\eightLambda_k}\big(\overline{z}_0,z \big)$.
Every walk contributing to this sum is an inner walk in $\eightLa_n$, moving from $\partial \eightLa_n$ to itself. The lattice $\eightLa_n$ is a fusion of eight copies of branch cuts of $\La_n$; by applying the map $x \to x^{1/8}$ before fusing these copies, $\eightLa_n$ may be embedded in the plane. This embedding easily demonstrates that the lift $\ga: \overline{z}_0 \to z$ to $\bbUi$ of any contributing walk has  $\vert {\rm wind}(\ga) \vert \leq 16\pi$. (By an abuse of notation, we continue to denote the endpoints of $\ga$ by $\overline{z}_0$ and $z$.)
Recall the constant $C$ provided by Lemma~\ref{lem:a2}. 
Any path $\gamma'$ contributing to $G_{\inftyOm_n,Cn^2}(\overline{z}_0,z)$ may be concatenated at both its endpoints to the lift $\ga$, 
so that a polygon $p:\overline{z_0} \to \overline{z_0}$ in $\inftyOm_k$ results. 
This polygon has length at most  $\vert \eightLa_n \vert + 1 + Cn^2 \leq C_2 n^2$ for some constant $C_2$. 
The application $(\gamma,\gamma') \to p$ is injective, because  $p$ crosses $\pi_\infty^{-1}(\partial\La_n)$ exactly twice, 
at $\overline{z}_0$ and $z$. In conclusion
\begin{align}\label{eq:opp}
	 \sum_{z\in \partial\eightLambda_n} G_{\eightLambda_n \setminus\eightLambda_k}(\overline{z}_0,z)G_{\inftyOm_n,Cn^2}(\overline{z}_0,z) 
	 \,\leq \, G_{\inftyOm_k , C_2 n^2}(\overline{z}_0,\overline{z}_0)\,,
\end{align}
where our notational abuse permits $\overline{z}_0$ and $z$ to refer both to certain points in the eight-fold cover 
and to their counterparts in the universal cover.
In the latter case however, $\overline{z}_0$ and $z$ are at distance at most $48n$ from each other 
due to the limited winding of the lift $\gamma$. 
Thus, Lemma~\ref{lem:a2} implies that 
$G_{\inftyOm_n,Cn^2}(\overline{z}_0,z)  \geq c \, n^{-103}$, 
for all $z$ in the sum above. 
By summing~(\ref{eq:op}) over $z$ and using the above, we find that
\begin{align*}
	\sum_{z\in \partial\eightLambda_n} G_{\eightLambda_n}(\overline{z}_0,z)
	\leq 
	n^{103}c^{-1}C_1 \, \Big( \max_{z \in \inftyOm_k} G_{\inftyOm_k,C_3 n^2}(z,z)\Big)^{6k+7}, 
\end{align*}
where $C_3 = \max\{72, C_2\}$ is a universal constant, 
and both the maximum of \eqref{eq:op} and the right-hand side of \eqref{eq:opp} were incorporated in the maximum above. 
Finally, Lemma~\ref{lem:a1} yields the desired inequality \eqref{eq:G_ub} 
with $C(k) = c^{-1} C_1 \cdot 8 \big(\cos  (3\pi/16) \big)^{-1}$. 
The condition on the growth of $C(k)$ may be verified by recalling that $C_1=x^{-72(k+1)^2}\,2^{432(k+1)^3}$.
\qed

\subsection{Proof of Lemma~\ref{lem:a3}}   
Our presentation will be fairly detailed but inexplicit regarding certain more minor aspects.

It is straightforward that $\mu(\bbUi_l)$, defined in~(\ref{e.uik}), decreases in $l$ to a limit that is at least $\mu(\bbH)$.
We suppose that $k \in \N \cup \{ 0 \}$ satisfies $\mu(\bbUi_k)=\mu(\bbH)$; and we may thus harmlessly suppose that $k$ is even. 

The positive $y$-axis meets edge midpoints in $\bbH$ at $y$-coordinates of the form $2^{-1}\sqrt{3} (1 + 2l)$, where $l \in \N \cup \{ 0 \}$. Note that $\big( 0, 2^{-1}\sqrt{3} (1 + k) \big)$ is the only element in this set of points that lies in $\partial \La_k$.

Set $z_0 = \big( 0, 2^{-1}\sqrt{3} (1 + k) , 0 \big)$.
We embed $\bbUi_k$ in $\R^3$ so that the vertical coordinate of any point $z \in \bbUi_k$ 
is the winding around the singularity of any path $\gamma : z_0 \to z$ in $\bbUi_k$. 

The cover $\bbUi_k$ lacks the translation invariance of $\bbH$, so that  two self-avoiding polygons in $\bbUi_k$ cannot in general be concatenated; thus, the standard method which permits an upper bound on the number of polygons on $\bbH$ is not directly available.
However, $\bbUi_k$ has a rotational $\Z$-invariance,
and thus pairs among a certain class of  polygons can be joined. 
We are about to make a definition in this regard, but first we mention that, in the present proof, we abuse the notation of Section~\ref{sec:observable} and interpret polygons as unrooted objects: a polygon here can be said to be a polygon trace as understood elsewhere in Sections~\ref{sec:observable} and~\ref{s.thm:a}.

A polygon in $\bbUi_k$ is $(i,j)$-{\em good} if it has
a unique lowest edge, and this edge is centred on $2^{-1} \sqrt{3}\big(0,i,0\big)$; and it has a unique highest edge, which is the rotation of the edge centred at 
 $2^{-1} \sqrt{3}\big(0,j,0\big)$ by an angle in $\pi/ 3 \cdot \bbN$. Here, and later, 
highest and lowest refer to the third coordinate in $\bbR^3$; note that $i$ and $j$ must be odd for such a polygon to exist.  
Let $g_n(i,j) = g_n(i,j;k)$ be the number of $(i,j)$-good polygons of length $n$. Then, for every $i,j,n,l$, we have that 
$g_{n+l-2}(i,i)\ge g_n(i,j)g_l(j,i)$. It is a simple consequence of the $\Z$-invariance of $\bbUi_k$ that $g_n(i,j) = g_n(j,i)$.  By Fekete's lemma applied to the sequence $g'_n(i,i)=g_{n+2}(i,i)$, we deduce that, for every odd $i$ and~$j$,
\begin{equation}\label{eq:jg}
	g_{n+2}(i,j)^{1/n}\le \limsup_{N} g_N(i,i)^{1/N} \, .
\end{equation}
Recall the specification of $\mu(\bbUi_k)$ in~\eqref{e.uik}.
Each polygon is a walk and thus $g_N(i,i) \leq  \sup_{v \in \bbUi_k} \big\vert \SAW_N(\bbUi_k,v) \big\vert$. Thus is the right-hand side of~(\ref{eq:jg})
seen to be at most $\mu(\bbUi_k)$. Since $\mu(\bbUi_k) = \mu(\bbH)$, we obtain that, for $n \in \N$, and $i,j \in \N$ odd,
\begin{equation}\label{e.gnij}
	g_{n+2}(i,j)^{1/n} \leq \mu(\bbH) \, . 	
\end{equation}

Now that we have bounded the number of good polygons, we wish to bound  $G_{\bbUi_k,n}(z,z)$ for every
$z \in \bbUi_k$. The idea will be to transform a polygon containing $z$ into a good polygon by concatenating polygons at its ends. 
  
The trace of a walk  $\gamma: z \to z$ that contributes to the sum $G_{\bbUi_k,n}(z,z)$ is a polygon $p$ containing $z$, of length at most $n$. We may index contributions by such polygons, there being two walks $\gamma$ in correspondence to each polygon~$p$. 
Those $p$ that contain a point at graph distance greater than $n/2$ from the centre of the singular hexagonal face are planar polygons, 
as they are too short to wind around the singularity. 
Using supermultiplicativity, it is standard to derive that the number of (unrooted) polygons of length $n$ on $\bbH$ 
is bounded above by $\mu(\bbH)^n$ (see for instance \cite{MadSla93}). 
When accounting for the choice of the root and the varying length, we conclude that the 
contribution to $G_{\bbUi_k,n}(z,z)$ of planar polygons, in particular those who venture outside $\inftyLambda_{n/2}$, is at most $2n^2$.

Thus, we may limit the rest of our study to contributing polygons which are not planar, and thus contained in $\inftyLambda_{n/2}$. 
Fix such a polygon $p$ and assume it is rotated so that its lowest point has height between zero and $\pi/3$.
Write $\ell$ for the line $\big\{ (0,y,0): y \geq 0 \big\}$ and write $\ell'$ for the line at height $\pi/3$
given by a counterclockwise rotation of $\ell$ by $\pi/3$ about the centre of the singular hexagon.
Let $u = (u_1,u_2,\pi/3)$ denote the point on $\ell' \cap p$ whose distance from this centre 
 is the greatest.
Note that $\vert u_1 \vert$ is the planar distance, measured using the first two coordinates, between $u$ and the line $\ell$. 

Suppose for now that $\vert u_1 \vert \geq 3k/2$.	
A polygon in $\bbH$ is called {\em useful} (for the values $n$ and $|u_1|$ above) if 
its width  is at least $n$ and, when $e$ denotes the most northerly among the rightmost (and thus vertical) edges in the polygon, 
there exists a vertex in the polygon whose $x$-coordinate is $\vert u_1 \vert - 3k/2$ less than $e$'s 
and whose $y$-coordinate is at most that of the lower of $e$'s endpoints. 
A convenient class of useful polygons is offered by the set $\FSAP^{(n)}$ from Section~\ref{s.polyhex}: any given element of this set is either useful or its vertical reflection is (recall that $|u_1| \leq n/2$).
		
\begin{figure}
	\begin{center}
		\includegraphics[width = 0.47\textwidth]{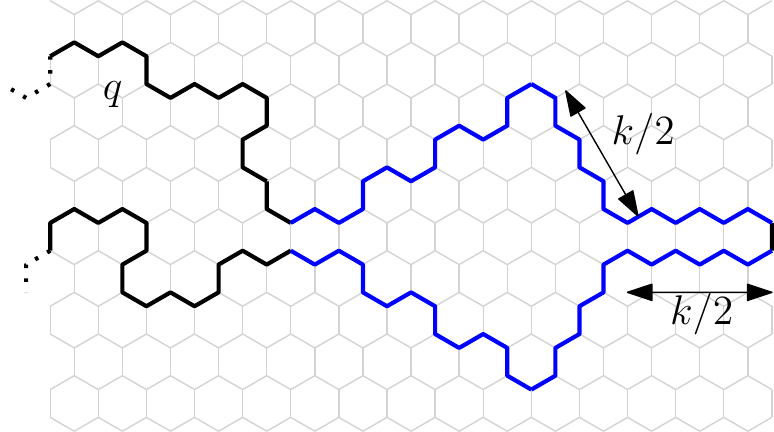}\quad
			\includegraphics[width = 0.47\textwidth]{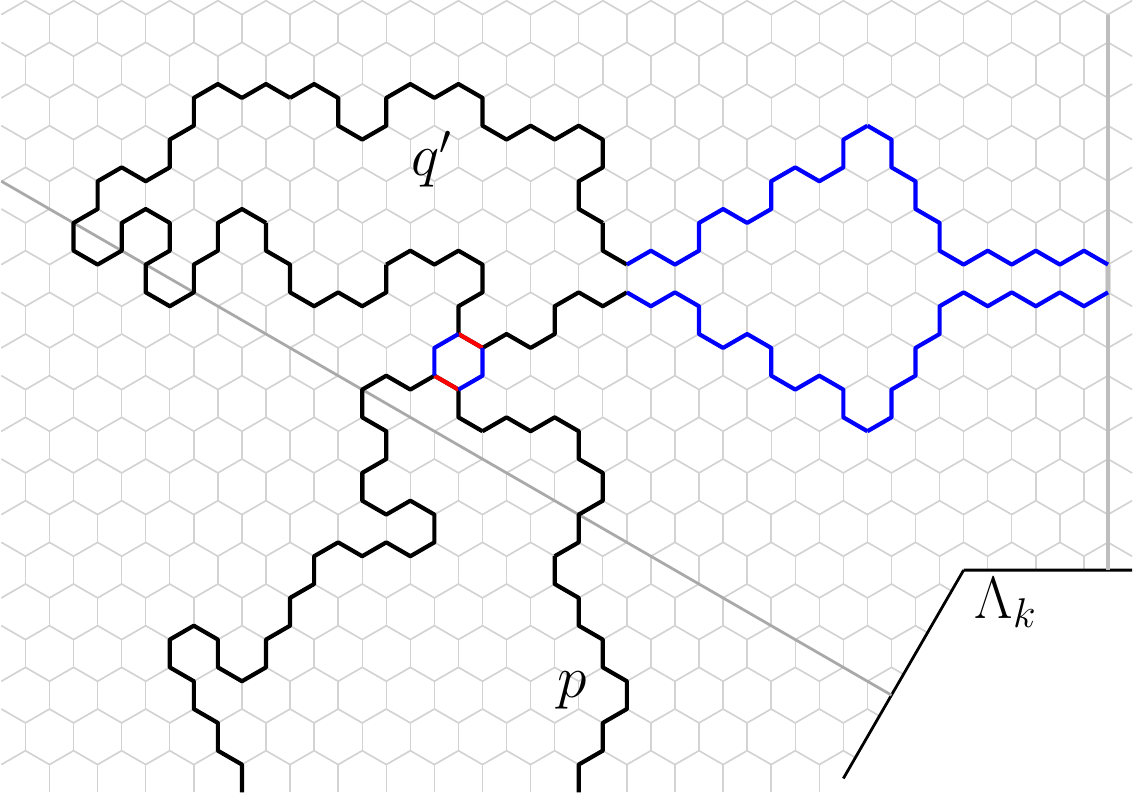}
		\caption{{\em Left:} Construct $q'$ by joining to the rightmost edge of $q$ the blue polygon shaped as an arrow. 
			{\em Right:} The polygon $q'$ is manoeuvred along  $\ell$ into a position whereby it may be joined to $p$ via the indicated hexagon.
			 In this instance, the join occurs by the removal of the two red edges  and the addition of the four blue ones. }
		   \label{fig:SAP_good}
	\end{center}
\end{figure}

Let $q$ denote a useful polygon.
We now join to the rightmost edge $e$ of $q$ the blue arrow-shaped polygon that is depicted in 
Figure~\ref{fig:SAP_good}. The result will be called $q'$. 
Embed $q'$ in $\bbUi_k$ so that its rightmost edge coincides with the lowest vertical edge of $\ell$.
This is indeed possible: the underside of the added blue polygon fits the contour that runs along the boundary of $\La_k$ 
from its intersection with $\ell$ to the intersection with $\ell'$, then along $\ell'$ until the contour recovers its original $y$-coordinate;
the part of $q'$ coinciding with $q$ lies to the left of the singularity.
	
Next $q'$ is translated away from the singularity in the direction of the line $\ell$, in steps of $\sqrt 3$. 
In doing so, there will be a moment when $q'$ will contain the point $u$, and thus intersect $p$. 
Further translate $q'$ until the last time that $q'$ and $p$ are within graph distance two of each other. 

At this stage, a hexagon $h$ may be found that borders both $p$ and the translate of $q'$, 
with contact with any given one of these polygons occurring along consecutive edges in $h$. 
Using the exclusive-or operation on $h$, we join $p$ and the translate of $q'$ and denote by $\overline{p}$ the result.

In fact, the same circumstance will arise when $\vert u_1 \vert < 3k/2$, 
whatever the definition of a useful polygon that is adopted in this case, 
because the underside of the blue polygon added on the right will contact $u$ when $q'$ is translated in the direction $\ell$. 
Again, we denote by $\overline{p}$ the result of the joining of $p$ with the appropriate translate of $q'$.

The polygon $\overline{p}$ has a length that differs from the sum of the constituents' by $-2$, $0$ or $2$.  
This polygon enjoys one of the two defining features of a good polygon, 
in that its least height is attained by a unique edge at height zero.
In order to produce a good polygon, a similar circumstance should be effected at the greatest height. 
This is ensured by a symmetric construction, where $\ell$ and $\ell'$ are now specified so that the sector between these line segments includes the highest point in $\overline{p}$, and with the useful polygon and its enhancement called $\bar{q}$ and~$\bar{q}'$.
	
	  The polygon that results from this double surgery  has length at most $|p|+|q'|+|\bar{q}'|+ 8 \le C' n^2$, where $C'$ may be chosen independently of $n$ and $k$ (provided that $n \geq k$). 
    Given the result of such a procedure, the values of $p$, $q$ and $\bar{q}$ may be retrieved 
    given the join locations  and the starting points of these three polygons.
    Thus, the map sending $\big( \gamma,q,\bar{q} \big)$ to the resulting good polygon is $C_0 n^{9}$-to-one for some constant $C_0$.
    We deduce that 
    $$
   \tfrac{1}{2} \, \overline{G}_{\bbUi_k,n}(z_0,z_0)\cdot\Big( \tfrac{1}{2}  \sum_{q\in \mathsf{FSAP}^{(n)} }x^{\vert q \vert}\Big)^2
    \, \leq \, C_0 \, n^{9}  \sum_{i,j,m\leq  C' n^2} g_m(i,j)\,x^m \, ,
    $$
    where recall that $x = \mu(\bbH)^{-1}$, and
    where the bar over  $G$ on the left-hand side is used to indicate that contributing $p$
    are contained in $\inftyLa_{n+3k/2}$.
    
    Applying~\eqref{eq:Bh} and Proposition~\ref{prop:polygon hexagonal}, we obtain
    $$
    \overline{G}_{\bbUi_k,n}(z_0,z_0)\le C_0 \, n^{23}  \sum_{i,j,m\le  C' n^2}g_m(i,j)\, x^m 
    $$
    after a suitable $k$-independent increase in the value of $C_0$. By~(\ref{e.gnij}), we see that, with a further such increase,  $\overline{G}_{\bbUi_k,n}(z_0,z_0) \leq  C_0 n^{32}$.
    Recalling, as noted earlier, that   $G_{\bbUi_k,n}(z_0,z_0)  - \overline{G}_{\bbUi_k,n}(z_0,z_0)  \leq 2n^2$,
    we obtain Lemma~\ref{lem:a3}. 
\qed

\noindent{\em Remark.} 
The proof of Lemma~\ref{lem:a3} succeeds under the hypothesis that, for some $k \in \N$, 
the number $g_N(i,i;k)$ of $(i,i)$-good polygons of length $N$ on $\bbU_k^\infty$ satisfies $\limsup_{N} g_N(i,i;k)^{1/N} \leq \mu(\bbH)$ for every odd $i$. 
Thus, Theorem~\ref{thm:a} is valid when Conjecture~\ref{conj:1}, which concerns walks, is replaced by the last condition, 
which merely concerns polygons. 
            	 
\subsection{Proof of Proposition~\ref{p.hex} in the general case}
%
%
Fix $n \geq k \geq 1$ such that $\mu(\bbU_k^\infty) = \mu(\bbH)$. 
In addition, let $y$ be the midpoint of an internal edge of $\La_{n-k}$, which may be supposed vertical with no loss of generality. 
Write $\La_n(y)$ for the translate of $\La_n$ by $y$. 
Then the conclusion of Proposition~\ref{p.hex} for $z = - y$ may be rewritten as 
\begin{align}\label{eq:hex_adapt}
	G_{\La_n(y)}(\originwalk,z_0) \leq  C(k) \, n^{2C(6k+7) + 103} \quad \text{ for all $z_0 \in \partial \La_n(y)$}, 
\end{align}
where $C$ and $C(k)$ satisfy the conditions of the proposition. 
Below, we adapt the steps in the proof of Proposition~\ref{p.hex} for $z =\originwalk$ to show~\eqref{eq:hex_adapt}.

Fix $z_0 \in \partial\La_n(y)$. Write $\eightLambda_n(y)$ for the lift of $\La_n(y)$ to $\bbU^8$. 
A straightforward adaptation of the proof of Lemma~\ref{lem:a1} yields the altered conclusion that
\begin{equation}\label{e.a1_adapt}
    G_{\La_n(y)}(\originwalk,z_0) 
    \, \le \, \frac8{\cos(3\pi/16)}\sum_{z\in \partial\eightLambda_n(y)} G_{\eightLambda_n(y)}\big( \overline{z}_0,z \big) \, .
\end{equation} 
Lemma~\ref{lem:a2} requires no alteration. 
Corollary \ref{cor:G_ub} may be adapted to yield 
\begin{align}\label{eq:G_ub_adapt}
	G_{\Lambda_n(y)}\big(\originwalk,z_0\big)
	\, \leq \, C(k) \, n^{103} \Big(\max_{z \in \inftyOm_k} G_{\inftyOm_k,c n^2}(z,z)\Big)^{6k +7} \, ,
\end{align}
Indeed, it suffices to replace $\eightLambda_n$ with $\eightLambda_n(y)$ in the original proof
and decompose walks contributing to the right-hand side of \eqref{e.a1_adapt} 
in terms of their visits to $\eightLambda_k \subset \eightLambda_n(y)$. 
In doing so, and when following the construction of the previous proof, 
we obtain a family of at most $6(k+1)$ polygons and one walk between points $z_0, z\in\partial\eightLambda_n(y)$, 
all contained in $\eightLambda_n(y) \setminus \eightLambda_k$. 
The latter walk is completed into a polygon of $\bbU^\infty_k$ of length at most $Ck^2$
using a walk on $\bbU^\infty \setminus \pi_\infty^{-1}(\La_n(y))$ between $z_0$ and $z$ (or rather between their lifts to $\bbU^\infty$). 
Notice that $\bbU^\infty \setminus \pi_\infty^{-1}(\La_n(y))$ and $\bbU^\infty_n$ are identical, 
and that the unaltered estimate of Lemma~\ref{lem:a2} may be used to prove \eqref{eq:G_ub_adapt}.

Finally, Lemma~\ref{lem:a3} allows us to deduce  \eqref{eq:hex_adapt} from  \eqref{eq:G_ub_adapt}.\qed

\appendix

\section{Polygon insertion for the hexagonal lattice}\label{sec:madras_join_hex}

Our goal is to prove Theorem~\ref{thm:HW2}
by means of Proposition~\ref{prop:many_poly_hex}.
To do this, we will introduce a join procedure on $\bbH$; record its pertinent properties; and explain perturbations needed to the proof in Sections~\ref{sec:proofsone} and~\ref{sec:proofstwo} of Theorem~\ref{thm:HW1}. Three subsections perform these respective tasks.

\subsection{A counterpart for $\bbH$ of Madras' join procedure}
Here, we present an $\bbH$-version of the procedure reviewed in Section~\ref{sec:MJ}. 

We continue to adopt the  notation $y_{\min}$ and $y_{\max}$. 
Let $\gamma$ be a self-avoiding bridge and $P$ be a polygon on $\bbH$ such that
\begin{equation}\label{eq:madras_cond_hex}
 	\llbracket y_{\min}(P) - 3 ,  y_{\max}(P) + 3 \rrbracket \subseteq  \llbracket y_{\min}(\gamma)  ,   y_{\max}(\gamma) \rrbracket \, . 
\end{equation} 
Let $F(\gamma)$ and $F(P)$ be the faces of $\bbH$ adjacent to at least one edge in $\gamma$ and $P$ respectively. 
Call a face of $F(\gamma)$ good if there are at most two edges of $\gamma$ adjacent to it and, 
in the case where there are two such edges, these two are also adjacent to each other. 
The same applies to $F(P)$. 

Slide $P$ horizontally to its rightmost position such that $F(\gamma)$ and $F(P)$ are at distance four in the dual graph: 
that is, the shortest path of adjacent faces from $F(\gamma)$ to  $F(P)$ contains five faces. 
This may be done because moving $P$ one unit horizontally changes the distance between $F(\gamma)$ and $F(P)$ by at most one.

\begin{figure}
	\begin{center}
	\includegraphics[width=0.55\textwidth]{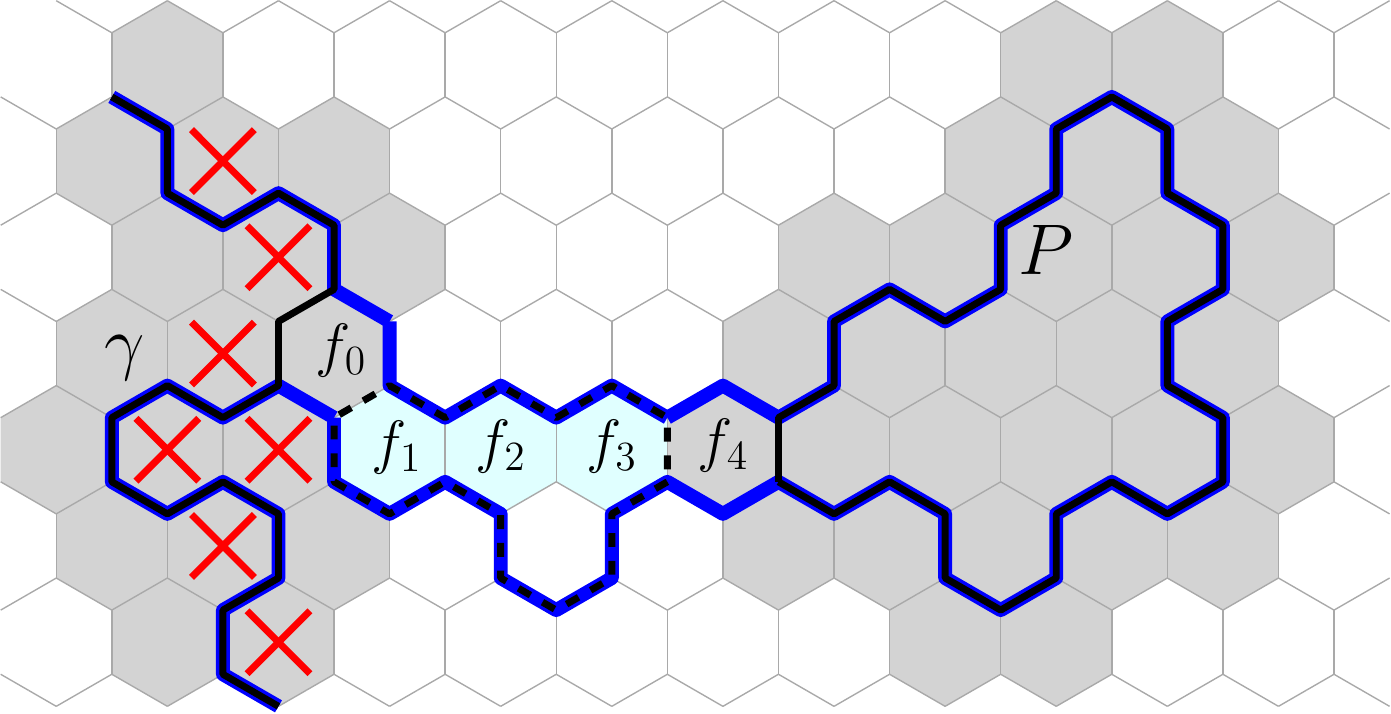}
	\caption{The join $M(\gamma,P)$ is emboldened. The grey faces are part of $F(\gamma)$ and $F(P)$; the bad faces of $\gamma$ are marked.
	The choice of $f_0$ and $f_4$ may be such that both are good. The polygon $Q$ is dotted.}
	\label{fig:MJ_hex_ex}
	\end{center}
\end{figure}

Let $f_0,\dots, f_4$ be a path of adjacent faces with $f_0\in F(\gamma)$ and $f_4 \in F(P)$. 
Notice that $f_0$ and $f_4$ may be chosen to be good. Indeed, given $f_1,f_2,f_3$, there is always a good neighbour of $f_1$ in $F(\gamma)$; and similarly for $f_3$. 
Moreover, by the minimality of the distance, all faces neighbouring $f_2$ are adjacent to neither $\gamma$ nor $P$. 

\begin{figure}[h]
	\begin{center}
	\includegraphics[width=1\textwidth]{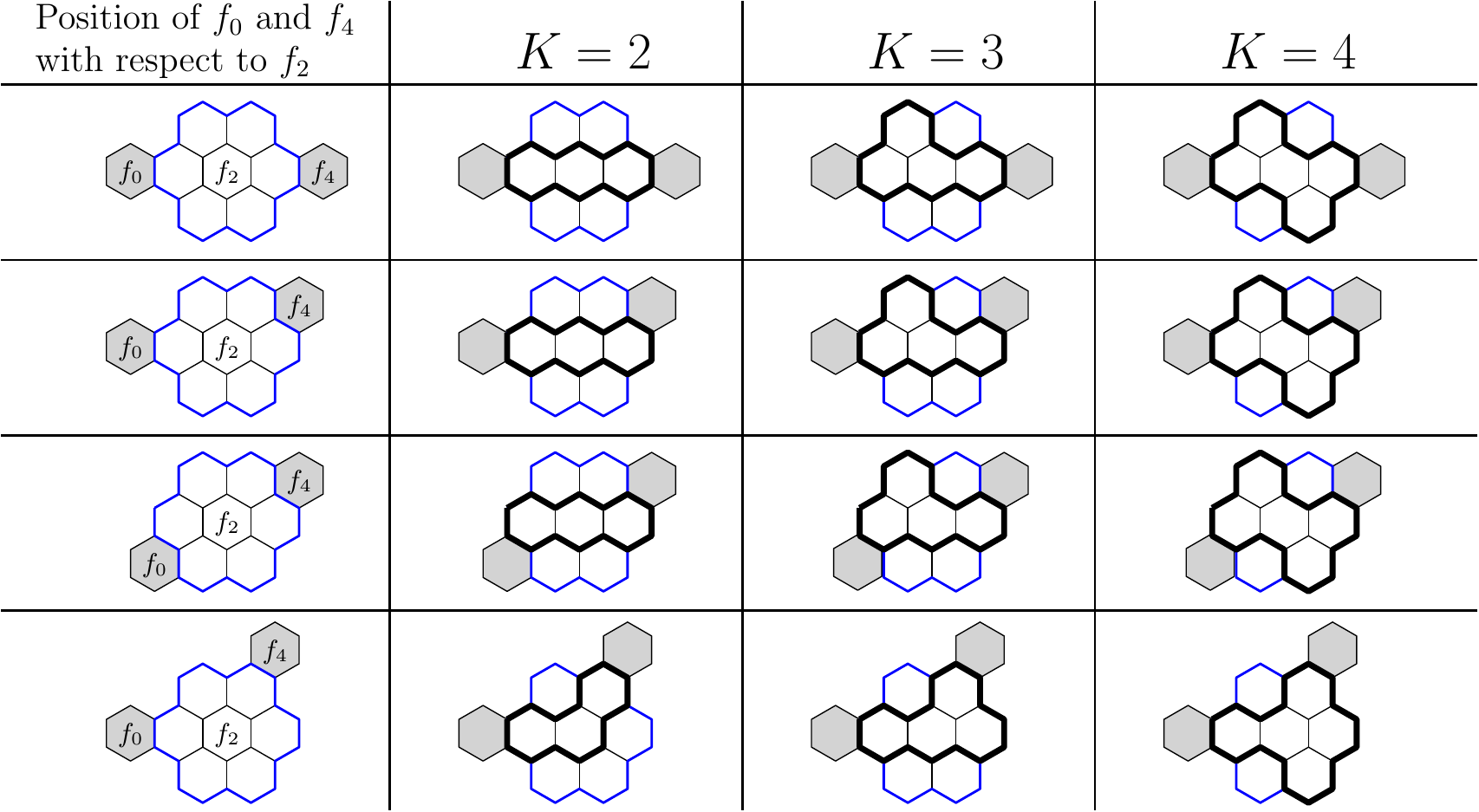}	
	\caption{The choice of $Q$, which is depicted in bold, is determined by $K$ and by the position of $f_0$ and $f_4$ relative to $f_2$.}
	\label{fig:MJ_hex}
	\end{center}
\end{figure}

Let $k_\gamma$ be the number of edges of $\gamma$ adjacent to $f_0$, and let $k_P$   be the number of edges of $P$ adjacent to $f_4$.
Set $K = k_\gamma + k_P$. The join of $\gamma$ and $P$ is then defined by using the exclusive-or operation $\textsf{xor}$: 
\begin{align*}
	M(\gamma,P) := \gamma \,\textsf{xor}\, f_0 \,\textsf{xor}\, Q \,\textsf{xor}\, f_4 \,\textsf{xor}\, P \, ,
\end{align*}
where, in an abuse of notation, $f_0$ and $f_4$ denote the polygons of length six that respectively surround the faces actually denoted by $f_0$ and $f_4$; and where $Q$ 
is a polygon whose interior contains $f_2$ and whose form  depends on the position of $f_0$ and of $f_4$ relative to $f_2$, as well as on $K$, in a manner that is indicated in Figure~\ref{fig:MJ_hex}. An example appears in Figure~\ref{fig:MJ_hex_ex}.

\subsection{Pertinent properties of the join on $\bbH$}
Lemma~\ref{l.madrascondition} records all required features of the Madras join. We need a counterpart for $\bbH$, but it would be unduly repetitive to explicitly state the analogous lemma. 

Rather, we focus on essentials. First note that, in view of the  choice of $f_0,\dots, f_4$ and of $Q$, the join  $M(\gamma,P)$ is a bridge of length 
\begin{align*}
	| M(\gamma,P)| = |\gamma| + |P| + |Q| + 12 - 2 K - 4 = |\gamma| + |P| +  18 \, ;
\end{align*}
so the change $16 \to 18$ is made from $\Z^2$ to $\bbH$. The role of the junction plaquette is now assumed by a {\em junction hexagon}, which we may set equal to $f_3$. The two adjacent vertices of $f_3$ whose intervening edge borders the grey hexagon in $F(P)$ indicated in Figure~\ref{fig:MJ_hex}
have chemical distance along $M(\gamma,P)$ equal to $\vert P \vert +3$.

\subsection{Perturbations to proofs} 
Assisted by fairly clear definitional changes, the $\bbH$-course of the proof of Theorem~\ref{thm:HW2}
offered by Sections~\ref{sec:proofsone}
and~\ref{sec:proofstwo}
runs rather unruffled in its $\Z^2$-tracks.  The change $29 \to 10$ occurs in the lower bound on $\alpha$ due to the use of Proposition~\ref{prop:many_poly_hex} in place of  Proposition~\ref{prop:many_poly_subseq}. Since  Proposition~\ref{prop:many_poly_hex}  holds for all $u \in \N$, every $n \in \N$ satisfies~\eqref{eq:polygonal_existence}; or rather the variant of~(\ref{eq:polygonal_existence}) obtained by the change $4^{-1} \to c_0$. The $\bbH$-counterpart to Proposition~\ref{prop:HW3}(1) thus offers an upper bound on $\vert \HSW_n \vert$ for all $n$. This permits the counterpart to Proposition~\ref{prop:HW3}(2) in which conditions  (\ref{eq:exponent_condition1}) and~(\ref{eq:exponent_condition2}) are used, rather than the strengthened conditions needed in the $\Z^2$ case. The $\bbH$-Proposition~\ref{prop:HW3}(2)  follows directly from~(\ref{eq:SAW_vs_HSW}). 
The bound $\eps < (42)^{-1}$ in Theorem~\ref{thm:HW2} arises by choosing $\alpha = \alpha(\eps) > 10$ small enough that~(\ref{eq:exponent_condition2}) with $\delta =2\eps$ holds. 

In regard to the near injectivity of $\Phi$, the value $L = L(\bbH)$ in Lemma~\ref{l.phikeyprop} equals three, whereas $L(\Z^2) = 12$. Indeed, it may be argued by examining the cases that, when $M(\gamma,P)$ and the junction hexagon are known, there is only one compatible choice of $(\gamma,P)$, whereas before there were possibly as many as four. 

\newcommand{\etalchar}[1]{$^{#1}$}

\end{document}